\documentclass[12pt]{article}
\usepackage{amsmath,amsthm,amssymb,mathtools}
\allowdisplaybreaks
\usepackage{enumitem}
\usepackage{etoolbox}
\usepackage{xcolor}
\usepackage{ulem}
\usepackage[colorlinks,
linkcolor=red,
anchorcolor=blue,
citecolor=magenta
]{hyperref}
\hypersetup{pdftitle={Two-Dimensional Full Compressible Navier--Stokes System with Vacuum at Infinity},pdfauthor={Xin Song},pdfsubject={Global strong solutions, physical small energy and large-time decay}}

\theoremstyle{plain}
\newtheorem{theorem}{Theorem}[section]
\newtheoremstyle{lixinlemma}
  {\topsep}{\topsep}{\itshape}{0pt}
  {\fontfamily{lmss}\bfseries}{}{.5em}{}
\theoremstyle{lixinlemma}
\newtheorem{lemma}[theorem]{Lemma}
\theoremstyle{plain}
\newtheorem{proposition}{Proposition}[section]
\newtheorem{corollary}[theorem]{Corollary}
\newtheorem{definition}{Definition}[section]
\newtheoremstyle{lixinremark}
  {\topsep}{\topsep}{\normalfont}{0pt}
  {\fontfamily{lmss}\bfseries\itshape}{}{.5em}{}
\theoremstyle{lixinremark}
\newtheorem{remark}{Remark}[section]
\theoremstyle{plain}

\numberwithin{equation}{section}

\newcommand{\R}{\mathbb R}
\newcommand{\dd}{\,\mathrm d}
\newcommand{\eps}{\varepsilon}

\newcommand{\Div}{\operatorname{div}}
\newcommand{\supp}{\operatorname{supp}}

\title{The global well-posedness and large time behavior of two-dimensional full compressible Navier--Stokes systems
with the far-field vacuum}

\author{Xinyu Fan$^{a} $,  Xin Song$^{b}$,  Xue Wang$^{c}$ \\
{\normalsize a. School of Mathematics  Sciences, } \\
{\normalsize  Nankai University, Tianjin 300110, P. R. China;}
 \\ {\normalsize b. School of Mathematics and Computer Sciences, }\\ {\normalsize \& Institute of Mathematics and Interdisciplinary Sciences,}\\ {\normalsize  Nanchang University, Nanchang 330031, P. R. China;}\\
{\normalsize c. School of Mathematics  Sciences,} \\ {\normalsize   Tiangong University, Tianjin 300380, P. R. China.}}
\date{}

\makeatletter
\patchcmd{\@title}{\\}{ }{}{}
\patchcmd{\@title}{\\}{ }{}{}
\patchcmd{\@title}{\\}{ }{}{}
\patchcmd{\@maketitle}{\large}{\normalsize}{}{}
\let\XiongOriginalMaketitle\maketitle
\renewcommand{\maketitle}{\setcounter{footnote}{0}\XiongOriginalMaketitle\vspace{-.7cm}}
\renewcommand{\@biblabel}[1]{[#1]}
\makeatother

\newif\ifXiongKeywordStyle
\AfterEndEnvironment{abstract}{%
  \par\begingroup\XiongKeywordStyletrue\footnotesize\let\noindent\indent
}
\pretocmd{\section}{\ifXiongKeywordStyle\par\endgroup\fi}{}{}

\BeforeBeginEnvironment{thebibliography}{%
  \begingroup\footnotesize
  \let\textsc\textnormal
  \let\emph\textnormal
  \let\textbf\textnormal
}
\AfterEndEnvironment{thebibliography}{\endgroup}

\begin{document}
\maketitle
\begingroup
\renewcommand{\thefootnote}{}
\footnotetext[0]{E-mail: fanxinyu17@mails.nankai.edu.cn(X. Fan), songxin2357@email.ncu.edu.cn (X. Song),  xuewa@amss.ac.cn (X.  Wang).}
\endgroup
\begin{abstract}
This paper is concerned with the Cauchy problem for the two-dimensional full
compressible Navier--Stokes equations with constant {viscosities} and heat
conductivity. Under the small initial energy assumption and the compatible condition, we establish the global existence and uniqueness of strong
solutions with far-field vacuum and derive that the system will converge to the steady state within the algebraic rates.  
\end{abstract}

\noindent\textbf{Keywords:} Full compressible Navier--Stokes system,
two-dimensional Cauchy problem, vacuum at infinity,
small initial energy, initial-layer estimates.

\section{Introduction}\label{sec:introduction}
We consider the full compressible Navier--Stokes system
\begin{equation}\label{eq:fns}
\left\{
\begin{aligned}
&\rho_t+\Div(\rho u)=0,\\
&(\rho u)_t+\Div(\rho u\otimes u)-\mu\Delta u
  -(\mu+\lambda)\nabla\Div u+\nabla P=0,\\
&c_v\big[(\rho\theta)_t+\Div(\rho\theta u)\big]
  -\kappa\Delta\theta+P\Div u
  =2\mu|D(u)|^2+\lambda(\Div u)^2,
\end{aligned}\right.
\end{equation}
where $t\ge0$ denotes time, $x=(x_1,x_2)\in\R^2$ is the spatial
variable, and $\rho$, $u=(u^1,u^2)$ and $\theta$ denote the density,
velocity, and absolute temperature, respectively. The pressure and the
deformation tensor are given by
\begin{equation}\label{eq:pressure}
P=\rho\theta,\qquad D(u)=\frac12\big(\nabla u+(\nabla u)^{\mathsf T}\big).
\end{equation}
Without loss of generality, we set the gas constant $R=1$ by rescaling
the temperature and the positive coefficients $c_v$ and $\kappa$.
The constant viscosity coefficients
satisfy the two-dimensional physical restrictions
\begin{equation}\label{eq:physical}
\mu>0,\qquad \mu+\lambda\ge0,
\end{equation}
 The positive constants
$c_v$ and $\kappa$ are, respectively, the specific heat at constant volume
and the heat conductivity coefficient.
The initial data are
\begin{equation}\label{eq:initial}
 (\rho,\rho u,\rho\theta)(x,0)
  =(\rho_0,\rho_0u_0,\rho_0\theta_0)(x),\qquad x\in\R^2.
\end{equation}

The density satisfies
\[
 \lim_{|x|\to\infty}\rho(x,t)=0,\qquad t>0.
\]

The compressible Navier--Stokes equations describe the motion of viscous
compressible fluids; when thermal effects are taken into account, the full
system \eqref{eq:fns} is also known as the Navier--Stokes--Fourier
system~\cite{CaoJiang2026,Wen2025}. Unlike the isentropic equations where
the pressure is a prescribed function of the density, \eqref{eq:fns}
couples the mass and momentum equations to the energy equation, and the
isentropic system cannot in general be recovered from \eqref{eq:fns} by
prescribing a particular temperature~\cite{Wen2025}. In the absence of
vacuum, local well-posedness was established by Nash~\cite{Nash1962},
while global classical solutions near a positive constant equilibrium were
first proved by Matsumura and Nishida~\cite{MatsumuraNishida1980} for
small perturbations in $H^3(\R^3)$ (see also
Danchin~\cite{Danchin2000Critical} for global strong solutions to the
barotropic system in critical homogeneous Besov spaces of $L^2$-type).
When vacuum is allowed, Cho and Kim~\cite{ChoKim2006} proved local
existence of strong solutions to the three-dimensional full system under
compatibility conditions on the initial data for both the momentum and
temperature equations, permitting both interior vacuum and vacuum at
infinity on bounded or unbounded domains.

For the three-dimensional isentropic Cauchy problem with $P=A\rho^\gamma$
($A>0$, $\gamma>1$), Huang, Li and Xin~\cite{HuangLiXin2012} established
the global existence and uniqueness of classical solutions for smooth,
compatible initial data with small total energy and possibly large
oscillations, allowing vacuum (even compactly supported initial density)
and either positive or zero far-field density. For the three-dimensional
full system with positive far-field density and temperature, Huang and
Li~\cite{HuangLi2018} obtained global classical solutions with small
initial energy, interior vacuum and large oscillations. When the far-field
state is vacuum, however, the total energy identity
\[
 \int_{\R^3}\rho\left(\frac{1}{2}|u|^2+c_v\theta\right)(t)\dd x
 =\int_{\R^3}\rho_0\left(\frac{1}{2}|u_0|^2+c_v\theta_0\right)\dd x
\]
does not directly yield a dissipation estimate (see~\cite{Li2020}). Unlike
the non-vacuum far-field case in~\cite{HuangLi2018}, where dissipation
estimates are recovered from the entropy inequality, Wen and
Zhu~\cite{WenZhu2017} used the conservation of mass to prove global
existence of strong solutions in $\R^3$ with vacuum at infinity provided
that the initial mass is sufficiently small. Later, Li~\cite{Li2020}
established global strong solutions to the three-dimensional full Cauchy
problem with vacuum at infinity under a smallness condition on a
scaling-invariant initial quantity with respect to the scaling
\[
 (\rho_0(x),u_0(x),\theta_0(x))
 \longmapsto
 (\rho_0(\ell x),\ell u_0(\ell x),\ell^2\theta_0(\ell x)),
 \qquad \ell>0,
\]
together with the viscosity restriction $2\mu>\lambda$. Recently,
Wen~\cite{Wen2025} proved global existence of strong solutions for both
the three-dimensional isentropic and full systems under scaling-invariant
smallness conditions, removing the restriction $2\mu>\lambda$ for the full
system while retaining only the physical conditions $\mu>0$ and
$2\mu+3\lambda\ge0$. Concerning blowup criteria for strong solutions of the
three-dimensional full system, Feireisl, Wen and
Zhu~\cite{FeireislWenZhu2024} estimated $\int\rho|u|^4\dd x$ using the
coupled terms in the temperature equation without imposing any extra
relation between the viscosity coefficients.

The large-time behavior of global solutions has also been studied near
positive constant equilibria. Using a time-weighted energy method in
$L^p$-type critical Besov spaces, Danchin and
Xu~\cite{DanchinXu2017Decay,DanchinXu2018FullDecay} established optimal
time-decay estimates for both the barotropic system and the full
compressible Navier--Stokes system in dimensions $d\ge3$ for small initial
perturbations satisfying an additional low-frequency smallness condition.

For the two-dimensional Cauchy problem with vacuum at infinity, the
analysis is more subtle since the three-dimensional Sobolev inequality
$\|u\|_{L^6(\R^3)}\le C\|\nabla u\|_{L^2(\R^3)}$ used
in~\cite{HuangLiXin2012} fails in $\R^2$, and the energy bounds do not
directly provide unweighted $L^p$-estimates for the velocity. For periodic
barotropic flows with constant shear viscosity and density-dependent bulk
viscosity $\lambda(\rho)=b\rho^\beta$ ($b>0$), Huang and
Li~\cite{HuangLi2016} established global existence of strong solutions
with large initial data allowing vacuum for $\beta>4/3$, and further
derived a time-independent upper bound for the density as well as the
large-time behavior when $\beta>3/2$ and $1<\gamma<4\beta-3$. In $\R^2$
with vacuum at infinity, Li and Liang~\cite{LiLiang2014} established a
local existence theory for the barotropic equations by introducing
spatially weighted density assumptions and weighted inequalities for the
velocity; under small initial total energy, Fang and
Guo~\cite{FangGuo2016} proved global existence and large-time asymptotic
behavior of strong solutions for constant shear viscosity and
$\lambda(\rho)=\rho^\beta$ ($0\le\beta\le1$). For constant viscosities,
Li and Xin~\cite{LiXin2019} proved the global existence, uniqueness and
large-time decay rates of strong and classical solutions to the
two-dimensional barotropic Cauchy problem with small initial energy and
vacuum at infinity (see also Li, Shang and Tang~\cite{LiShangTang2020} for
global classical solutions to the two-dimensional isentropic equations in
a square with homogeneous Dirichlet boundary conditions and small initial
mass). For the two-dimensional full system \eqref{eq:fns} with constant
coefficients and vacuum at infinity, Wang~\cite{Wang2023} established
local existence and a blowup criterion for strong solutions requiring only
the momentum compatibility condition \eqref{eq:compatibility}.

In this paper, we establish the global existence, uniqueness and
large-time decay rates of strong solutions to the two-dimensional Cauchy
problem \eqref{eq:fns}--\eqref{eq:initial} with constant viscosity and
heat-conductivity coefficients and vacuum at infinity, under only the
physical viscosity restrictions \eqref{eq:physical} and the compatibility
condition \eqref{eq:compatibility}. Here the smallness condition
\eqref{eq:smallness} is imposed solely on the initial total energy $E_0$
in \eqref{eq:energy-definition}, without requiring smallness of the total
mass $m_0$, $\|\nabla\theta_0\|_2$, or $\|g\|_2$. In addition, we
determine the far-field limit of the velocity in
\eqref{eq:velocity-at-infinity} and obtain explicit large-time decay rates
in Theorem~\ref{thm:decay}.

All unmarked spatial norms below are over $\R^2$.
We use the following definition of a strong solution
\cite{Wang2023}.

\begin{definition}\label{def:strong}
 $(\rho,u,\theta)$ is called a strong solution to \eqref{eq:fns}
provided that all derivatives appearing in the system are regular
distributions and the equations are satisfied almost everywhere in
$\R^2\times(0,T)$.
\end{definition}
The initial energy is\begin{equation}\label{eq:energy-definition} E_0=\int_{\R^2}\left(\tfrac12\rho_0|u_0-u_\infty|^2+c_v\rho_0\theta_0\right)\dd x,
\end{equation}
where $$u_\infty=\lim_{|x|\to\infty}u_0(x).$$
Actually, Proposition~\ref{prop:initial-velocity-limit} in Appendix~\ref{app:initial-velocity} proves that this limit exists.
Assume that the initial density has finite positive mass,
\begin{equation}\label{eq:mass}
 m_0=\int_{\R^2}\rho_0\dd x\in(0,\infty).
\end{equation}
Choose $N_0>0$ such that
\begin{equation}\label{eq:masslocal}
\int_{B_{N_0}}\rho_0\dd x\ge\frac12\int_{\R^2}\rho_0\dd x
 =\frac{m_0}{2}.
\end{equation}

\begin{theorem}\label{thm:target}
Let $\eta_0>0$ and set
\[
\bar x=(e+|x|^2)^{1/2}\log^{1+\eta_0}(e+|x|^2).
\]
For given constants $a>1$, $q>2$, $M>0$ and $\bar\rho>0$,
we assume that the initial data
$(\rho_0\ge0,u_0,\theta_0\ge0)$ satisfy
\begin{equation}\label{eq:data}
\begin{gathered}
\bar x^a\rho_0\in L^1\cap H^1\cap W^{1,q},\qquad
\sqrt{\rho_0}u_0,\ \sqrt{\rho_0}\theta_0\in L^2,
\end{gathered}
\end{equation}
and
\begin{equation}\label{eq:data-second}
\begin{gathered}
\nabla u_0\in H^1,\qquad \nabla\theta_0\in L^2,
\end{gathered}
\end{equation}
and the compatibility condition
\begin{equation}\label{eq:compatibility}
-\mu\Delta u_0-(\mu+\lambda)\nabla\Div u_0
   +\nabla(\rho_0\theta_0)=\sqrt{\rho_0}\,g
\end{equation}
with $g\in L^2$ and assume further that
\begin{equation}\label{eq:initial-data-bound}
0\le\rho_0\le\bar\rho,\qquad
\|\bar x^a\rho_0\|_{L^1}+\|\nabla u_0\|_2
       +\|\sqrt{\rho_0}\theta_0\|_2\le M.
\end{equation}
 Then there exists $\eps_0\in(0,1]$, depending only on
$m_0,\bar\rho,M,\mu,\lambda,c_v,\kappa$, such that if
\begin{equation}\label{eq:smallness}
 E_0\le\eps_0,
\end{equation}
the Cauchy problem for \eqref{eq:fns} with initial data \eqref{eq:initial} has a unique global strong
solution $(\rho\ge0,u,\theta\ge0)$. For every finite $T>0$, we have
\[
 0\le\rho(x,t)\le\tfrac32\bar\rho
 \quad\hbox{on }\R^2\times[0,T],
\]
and
\begin{equation}\label{eq:regularity}
\left\{
\begin{aligned}
&\rho\in C([0,T];L^1\cap H^1\cap W^{1,q}),\\
&\bar x_t^a\rho\in L^\infty(0,T;L^1\cap H^1\cap W^{1,q}),\\
&\nabla u\in L^\infty(0,T;H^1)\cap L^2(0,T;W^{1,q}),\\
&\nabla\theta,\ \sqrt t\nabla^2\theta
   \in L^\infty(0,T;L^2)\cap L^2(0,T;H^1),\\
&\sqrt\rho(u-u_\infty),\ \sqrt\rho\theta,\ \sqrt\rho u_t,
   \ \sqrt{t\rho}\,\theta_t,\ \bar x_t^{\alpha/2}\nabla u
   \in L^\infty(0,T;L^2),\\
&\nabla u_t,\ \sqrt\rho\theta_t,\ \sqrt t\,\nabla\theta_t
   \in L^2(0,T;L^2),
\end{aligned}\right.
\end{equation}
where
\[
 \bar x_t=(e+|x-tu_\infty|^2)^{1/2}\log^{1+\eta_0}(e+|x-tu_\infty|^2),\qquad
 \alpha=\min\{a/2,1\}.
\]
For each $t\ge0$,
the velocity has the spatial limit
\begin{equation}\label{eq:velocity-at-infinity}
\begin{gathered}
 \lim_{|x|\to\infty}u(x,t)
 =\frac1{m_0}\int_{\R^2}\rho_0u_0\dd x
  +\frac1{m_0}\int_{\R^2}(\nabla\Delta^{-1}\rho\cdot\nabla)u\dd x,
\end{gathered}
\end{equation}
and
\begin{equation}\label{eq:velocity-at-infinity-second}
\begin{gathered}
 u(\cdot,t)-\lim_{|y|\to\infty}u(y,t)\in L^{16}(\R^2).
\end{gathered}
\end{equation}
\end{theorem}

\begin{remark}\label{rem:local-existence}
Suppose that $\rho_0,\theta_0\ge0$ and that \eqref{eq:mass},
\eqref{eq:data}--\eqref{eq:data-second} and \eqref{eq:compatibility} hold.
By arguments similar to those in \cite[Proposition~3.1, Lemma~3.6 and Section~3.2]{Wang2023}
(see also \cite{LiXin2019}), there exists a time $T_0>0$ such that the Cauchy problem
\eqref{eq:fns}--\eqref{eq:initial} admits a unique local strong solution $(\rho,u,\theta)$
on $\R^2\times[0,T_0]$ satisfying \eqref{eq:regularity}. Moreover, for a radius $R>0$ with
$\int_{B_R}\rho_0\dd x\ge m_0/2$, $T_0$ depends only on
$m_0,R,a,q,\eta_0,\mu,\lambda,c_v,\kappa$, $\|g\|_2$ and the initial norms in
\eqref{eq:data}--\eqref{eq:data-second}.
\end{remark}

\begin{remark}\label{rem:classical-local-agreement}
	If the initial data satisfy additional regularity, the strong
	solution in Theorem~\ref{thm:target} is classical one.
\end{remark}

\begin{remark}\label{rem:smallness-dependence}
In addition to $m_0,\mu,\lambda,c_v$ and $\kappa$, the smallness threshold $\eps_0$ in
\eqref{eq:smallness} depends only on the constants $\bar\rho$ and $M$ in
\eqref{eq:initial-data-bound}, which control
$\|\rho_0\|_\infty$ and
$\|\bar x^a\rho_0\|_1+\|\nabla u_0\|_2+\|\sqrt{\rho_0}\theta_0\|_2$,
respectively, and is independent of all other norms of the initial data.
\end{remark}

With $u_*=m_0^{-1}\int_{\R^2}\rho_0u_0\dd x$ denoting the mean velocity and
$E_0$ as in \eqref{eq:energy-definition}, we define the limiting temperature by
\[
 \theta_*=\frac{E_0-\frac12m_0|u_*-u_\infty|^2}{c_vm_0}
 =\frac{1}{c_vm_0}\int_{\R^2}
   \left(\tfrac12\rho_0|u_0-u_*|^2+c_v\rho_0\theta_0\right)\dd x\ge0.
\]

\begin{theorem}\label{thm:decay}
Let $(\rho,u,\theta)$ be the global strong solution of Theorem~\ref{thm:target}.
There exists $C>0$ such that, for all $t\ge2$, we have
\begin{equation}\label{eq:main-decay}
\begin{aligned}
 &\|\nabla u(t)\|_2+\|P(t)\|_2\le C(1+t)^{-3/8},
       \qquad \|P(t)\|_3\le C(1+t)^{-1/2},\\
 &\|\nabla\theta(t)\|_2
       \le C\log^{1/2}(e+t)(1+t)^{-9/16},\\
 &\|\sqrt\rho(u-u_*)(t)\|_2
       \le C\log^{1/2}(e+t)(1+t)^{-3/8},\\
 &\|\sqrt\rho(\theta-\theta_*)(t)\|_2
       \le C\log(e+t)(1+t)^{-9/16},
\end{aligned}
\end{equation}
where $C=C(m_0,\bar\rho,M,\mu,\lambda,c_v,\kappa)$.
Moreover, as $t\to\infty$, we have
\[
 \frac12\int_{\R^2}\rho|u|^2\dd x\longrightarrow\frac{m_0}{2}|u_*|^2,\qquad
 \frac1{m_0}\int_{\R^2}\rho\theta\dd x\longrightarrow\theta_*.
\]
\end{theorem}

\begin{remark}
By \eqref{eq:main-decay}, \eqref{eq:conserved-limits-second} and
\eqref{eq:unshifted-kinetic-identity}, one also obtains the convergence rate
of the kinetic energy $\frac12\int_{\R^2}\rho|u|^2\dd x$ and the internal
energy $c_v\int_{\R^2}\rho\theta\dd x$ as $t\to\infty$, namely,
\[
 \frac12\int_{\R^2}\rho|u|^2\dd x-\frac{m_0}{2}|u_*|^2
 =c_v\left(m_0\theta_*-\int_{\R^2}\rho\theta\dd x\right)
 =O(t^{-3/4}\log t).
\]
\end{remark}

The Galilean transformation
\begin{equation}\label{eq:moving-frame}
 (\rho,u,\theta)(x,t)\longmapsto
 \bigl(\rho(x+tu_\infty,t),\,u(x+tu_\infty,t)-u_\infty,\,
       \theta(x+tu_\infty,t)\bigr)
\end{equation}
gives the same system with initial data $(\rho_0,u_0-u_\infty,\theta_0)$.
Since $\nabla(u_0-u_\infty)=\nabla u_0$, \eqref{eq:compatibility} and
\eqref{eq:initial-data-bound} hold with the same $g$ and $M$, respectively.
Still denoting the transformed solution and initial data by $(\rho,u,\theta)$
and $(\rho_0,u_0,\theta_0)$, respectively, we can rewrite
\eqref{eq:energy-definition} as
\[
 E_0=\int_{\R^2}\left(\tfrac12\rho_0|u_0|^2+c_v\rho_0\theta_0\right)\dd x.
\]
For a function $f$, write $\dot f=f_t+u\cdot\nabla f$. Then we have
\begin{equation}\label{eq:perturbation-system}
\left\{
\begin{aligned}
 &\rho_t+\Div(\rho u)=0,\\
 &\rho\dot u-\mu\Delta u-(\mu+\lambda)\nabla\Div u+\nabla P=0,\\
 &c_v\rho\dot\theta-\kappa\Delta\theta+P\Div u
        =2\mu|D(u)|^2+\lambda(\Div u)^2.
\end{aligned}\right.
\end{equation}

We now make some comments on the analysis of this paper. By
Remark~\ref{rem:local-existence}, to extend the local strong solution
globally in time, it suffices to establish global a priori estimates for
smooth solutions to \eqref{eq:perturbation-system}. As in the barotropic case \cite{LiXin2019}, the key issue here is
to derive both the time-independent upper bound for the density and the
time-dependent higher-order estimates for $(\rho,u,\theta)$, so some basic
ideas in \cite{LiXin2019} will be adapted here, yet new difficulties arise
for the full compressible Navier--Stokes system
\eqref{eq:perturbation-system} with vacuum at far field. Indeed, unlike
the barotropic case $P=a\rho^\gamma$, the pressure $P=\rho\theta$ is
coupled with the temperature $\theta$. On the short-time interval
$[0,\min\{T,\tau\}]$ (for a small time $\tau>0$), multiplying
\eqref{eq:perturbation-system}$_3$ by $\theta$ and integrating by parts
leads to the integral $\int_{\R^2}\theta|\nabla u|^2\dd x$, which cannot
be bounded directly from the energy estimate for $\int_{\R^2}\rho\theta\dd x$
when vacuum is present. For $t\ge\tau$, one needs sufficient time decay of
the effective viscous flux $F=(2\mu+\lambda)\Div u-P$ so that
$\int_\tau^T\|F\|_\infty\dd t$ is bounded uniformly in $T$. However, the
pressure estimate \eqref{eq:full-pressure-power-balance} with $p=3$
contains the term $3\int_{\R^2}\rho P^2\dot\theta\dd x$, so one must also
estimate $\sqrt\rho\,\dot\theta$.

To overcome these difficulties, on $[0,\min\{T,\tau\}]$ we apply H\"older's
inequality with spatial weights on $\theta$ and $\nabla u$ to estimate
$\int_{\R^2}\theta(2\mu|D(u)|^2+\lambda(\Div u)^2)\dd x$ in
\eqref{eq:initial-sublinear-heating}. From the energy estimates for
$\theta$, $F$ and the vorticity, we obtain \eqref{eq:initial-gradient-time},
which controls the time integral of $\|\nabla u\|_4$ in the commutator
estimate \eqref{eq:initial-commutator-bound}. Integrating
\eqref{eq:signed-potential-identity} along particle paths and using
\eqref{eq:momentum-potential-bound}--\eqref{eq:initial-gradient-time}, we
deduce the short-time upper bound \eqref{eq:initial-density-gain} for the
density. Moreover, by introducing time weights in
Lemmas~\ref{lem:initial-time-weighted}--\ref{lem:short-squared-weight}, we
can choose $\eps_0$ independently of $\|g\|_2$ and $\|\nabla\theta_0\|_2$.
On $[\tau,T]$, to handle the term involving $\dot\theta$ in
\eqref{eq:full-pressure-power-balance}, we use the solution $\phi$ of
\eqref{eq:lift-elliptic-source} so that
$2\mu|D(u)|^2+\lambda(\Div u)^2-P\Div u$ cancels out in the equation for
$\theta-\phi$. We then estimate
$\sqrt\rho\bigl(\dot\theta-m_0^{-1}\int_{\R^2}\rho\dot\theta\dd y\bigr)$ in
Lemma~\ref{lem:temperature-residual-energy} and
$m_0^{-1}\int_{\R^2}\rho\dot\theta\dd x$ in
\eqref{eq:lift-temperature-mean-identity}. Together with the estimates for
$\sqrt\rho\,\dot u$ and $P$, these temperature estimates yield the
time-weighted decay in
Lemmas~\ref{lem:weighted-decay}--\ref{lem:endpoint-weighted-decay} and hence
$\int_\tau^T\|F\|_\infty\dd t\le CE_0^{1/64}$ in
\eqref{eq:late-log-weighted-flux}. Since $P\ge0$, it follows from
\eqref{eq:perturbation-system}$_1$, \eqref{eq:initial-density-gain} and
\eqref{eq:density-characteristics}--\eqref{eq:original-smallness-threshold}
that $\rho\le\frac32\bar\rho$ on $[0,T]$.

Once the uniform upper bound for the density is established, we turn to the
higher-order estimates on $[0,T]$. We first derive the $L^\infty(0,T;L^2)$
bounds for $\sqrt\rho\,\dot u$ and $\nabla\theta$ in
\eqref{eq:higher-final-bounds}, and then obtain the $L^2(0,T;L^\infty)$
bound for $\theta$ in \eqref{eq:higher-final-bounds-second} via
\eqref{eq:continuation-newton} (see
Lemma~\ref{lem:higher-acceleration-temperature}). Using these estimates
together with the Beale--Kato--Majda type inequality
(Lemma~\ref{lem:bkm}), we bound $\nabla\rho$ in
$L^\infty(0,T;L^2\cap L^q)$ and $\nabla u$ in $L^2(0,T;L^\infty)$ through
a logarithmic Gronwall inequality (see Lemma~\ref{lem:higher-density}).
With these bounds at hand, we establish the spatially weighted estimates
for $\rho$ and $\nabla u$ and complete the higher-order estimates for the
velocity and temperature in
Lemmas~\ref{lem:higher-weighted-density}--\ref{lem:higher-thermal-estimates};
in particular, the spatially weighted estimate for $\nabla u$ gives the
far-field limit \eqref{eq:velocity-at-infinity}. Finally, the large-time
behavior in Theorem~\ref{thm:decay} follows from
\eqref{eq:lift-log-anchor}, \eqref{eq:physical-total-energy} and the decay
estimates for $\nabla u$ and $\nabla\theta$.

The rest of this paper is organized as follows. In Section~\ref{sec:two},
we collect some preliminary lemmas and inequalities that will be used
later. We derive the lower-order and higher-order a priori estimates in
Sections~\ref{sec:three} and~\ref{sec:four}, respectively, and complete
the proofs of Theorems~\ref{thm:target} and~\ref{thm:decay} in
Section~\ref{sec:five}. Finally, in Appendix~\ref{app:initial-velocity},
we prove the existence of the far-field limit $u_\infty$ of the initial
velocity.

\section{Preliminaries}\label{sec:two}
In this section, we recall some basic facts and elementary inequalities on
$\R^2$ that will be used in the a priori estimates.

We begin with the Gagliardo--Nirenberg inequalities \cite{Nirenberg1959},
which we use later.

\begin{lemma}[Gagliardo--Nirenberg]\label{lem:gn}
Let $p\in[2,\infty)$, $q\in(1,\infty)$ and $r\in(2,\infty)$.
There exists $C>0$ such that, for $f\in H^1(\R^2)$ and
$g\in L^q(\R^2)\cap D^{1,r}(\R^2)$, we have
\begin{align}
 &\|f\|_{L^p(\R^2)}^p
 \le C\|f\|_{L^2(\R^2)}^2
       \|\nabla f\|_{L^2(\R^2)}^{p-2},
      \label{eq:gn-lp}\\
\intertext{and}
 &\|g\|_{C(\overline{\R^2})}
 \le C\|g\|_{L^q(\R^2)}^{q(r-2)/(2r+q(r-2))}
       \|\nabla g\|_{L^r(\R^2)}^{2r/(2r+q(r-2))},
      \label{eq:gn-infty}
\end{align}
where $C=C(p,q,r)$.
\end{lemma}

The following weighted $L^m$ inequality in $\R^2$ follows from
\cite[Theorem~B.1]{Lions1996}.
\begin{lemma}\label{lem:weighted-poincare}
Let $m\ge2$ and $\zeta>1+m/2$. For every $v\in H^1_{\rm loc}(\R^2)$
with $\nabla v\in L^2$, we have
\begin{equation}\label{eq:weighted-poincare}
 \left(\int_{\R^2}\frac{|v|^m}{e+|x|^2}
              \log^{-\zeta}(e+|x|^2)\dd x\right)^{1/m}
 \le C\|v\|_{L^2(B_1)}+C\|\nabla v\|_2,
\end{equation}
where $C=C(m,\zeta)$.
\end{lemma}

Using Poincar\'e's inequality and \eqref{eq:weighted-poincare} we get the following inequalities.
\begin{lemma}
Suppose $0\le\rho\le 2\bar\rho$, $\int_{\R^2}\rho=m_1>0$ and $\int_{\R^2}|x|\rho\le m_2$.\label{lem:density-weighted}
For $1\le p<\infty$ and $w\in H^1_{\rm loc}$ with $\nabla w\in L^2$,
we have
\begin{equation}\label{eq:density-anchor}
 \|w\|_{L^p(\rho\,\dd x)}
 \le C\left(\int_{\R^2}\rho|w|\,\dd x+\|\nabla w\|_2\right),
\end{equation}
where $C=C(\bar\rho,m_1,m_2,p)$ and $\rho w\in L^1$, and
\predisplaypenalty=0\relax
\begin{equation}\label{eq:spatial-weighted-lp}
 \|(1+|x|^2)^{-\eta/2}w\|_p
 \le C\bigl(\|\sqrt\rho\,w\|_2+\|\nabla w\|_2\bigr),
\end{equation}
where $C=C(\bar\rho,m_1,m_2,\eta,p)$, $\sqrt\rho\,w\in L^2$ and $\eta p>2$.
\end{lemma}
\begin{proof}
For convenience, we set $r=1+4m_2/m_1$. Then we have $\int_{B_r}\rho\dd x\ge m_1-m_2/r\ge m_1/2$.
Writing $w_{B_r}$ for the average of $w$ over $B_r$, we obtain from Poincar\'e's inequality that
\[
\begin{aligned}
 \frac{m_1}{2}|w_{B_r}|
 &\le \int_{\R^2}\rho|w|\dd x
      +\|\rho\|_{L^2(B_r)}\|w-w_{B_r}\|_{L^2(B_r)}\\
 &\le \int_{\R^2}\rho|w|\dd x+Cr(2\bar\rho m_1)^{1/2}\|\nabla w\|_2,
\end{aligned}
\]
and
\begin{equation}\label{eq:local-mean-l2}
 \|w\|_{L^2(B_1)}
 \le \|w-w_{B_r}\|_{L^2(B_r)}+|B_1|^{1/2}|w_{B_r}|
 \le C\left(\int_{\R^2}\rho|w|\dd x+\|\nabla w\|_2\right).
\end{equation}
By H\"older's inequality and \eqref{eq:weighted-poincare} with $m=4p$ and $\zeta=2+2p$, we obtain
\begin{equation}\label{eq:density-lp-local}
\begin{aligned}
 \int_{\R^2}\rho|w|^p\dd x
 &\le \|\rho\|_\infty^{1/4}
       \left(\int_{\R^2}(1+|x|)\rho\dd x\right)^{3/4}
       \left(\int_{\R^2}\frac{|w|^{4p}}{(1+|x|)^3}\dd x\right)^{1/4}\\
 &\le C\left(\|w\|_{L^2(B_1)}+\|\nabla w\|_2\right)^p,
\end{aligned}
\end{equation}
where $(1+|x|)^{-3}\le C(e+|x|^2)^{-1}\log^{-2-2p}(e+|x|^2)$.
Substituting \eqref{eq:local-mean-l2} into \eqref{eq:density-lp-local} and taking the $p$th root gives \eqref{eq:density-anchor}.

For $\eta p>2$, we have
\begin{equation}\label{eq:spatial-weight-comparison}
 \frac{1}{(1+|x|^2)^{\eta p/2}}
 \le \frac{C}{(e+|x|^2)\log^{2+p}(e+|x|^2)}.
\end{equation}
By Cauchy--Schwarz, \eqref{eq:weighted-poincare} with $m=2p$
and \eqref{eq:spatial-weight-comparison}, we get
\begin{equation}\label{eq:spatial-lp-local}
\begin{aligned}
 \int_{\R^2}\frac{|w|^p}{(1+|x|^2)^{\eta p/2}}\dd x
 &\le \left(\int_{\R^2}\frac{\dd x}{(1+|x|^2)^{\eta p/2}}\right)^{1/2}
       \left(\int_{\R^2}\frac{|w|^{2p}}{(1+|x|^2)^{\eta p/2}}\dd x\right)^{1/2}\\
 &\le C\left(\|w\|_{L^2(B_1)}+\|\nabla w\|_2\right)^p.
\end{aligned}
\end{equation}
Using \eqref{eq:local-mean-l2} and $\int_{\R^2}\rho|w|\dd x\le m_1^{1/2}\|\sqrt\rho\,w\|_2$ in \eqref{eq:spatial-lp-local} gives \eqref{eq:spatial-weighted-lp}.
\end{proof}

We next derive the energy identities and the first-moment estimate from \eqref{eq:perturbation-system}.
\begin{lemma}\label{lem:physical-energy-mean-temperature}
Let $(\rho\ge0,u,\theta\ge0)$ be a smooth solution of
\eqref{eq:perturbation-system} on $\R^2\times[0,T]$ in the class
\eqref{eq:regularity}, with initial data \eqref{eq:initial}. Then, for $0\le t\le T$,
$\int_{\R^2}\rho(t)\dd x=m_0$ and
\begin{align}
 &\begin{aligned}
  &\frac12\int_{\R^2}(\rho|u|^2)(x,t)\dd x
    +\int_0^t\bigl(\mu\|\nabla u(s)\|_2^2
       +(\mu+\lambda)\|\Div u(s)\|_2^2\bigr)\dd s\\
  =&\frac12\int_{\R^2}\rho_0|u_0|^2\dd x
       +\int_0^t\int_{\R^2}(P\Div u)(x,s)\dd x\dd s,
 \end{aligned}\label{eq:physical-kinetic-energy}\\
 &\begin{aligned}
  &c_v\int_{\R^2}(\rho\theta)(x,t)\dd x
       +\int_0^t\int_{\R^2}(P\Div u)(x,s)\dd x\dd s\\
  =&c_v\int_{\R^2}\rho_0\theta_0\dd x
       +\int_0^t\bigl(\mu\|\nabla u(s)\|_2^2
          +(\mu+\lambda)\|\Div u(s)\|_2^2\bigr)\dd s,
 \end{aligned}\label{eq:physical-thermal-energy}\\
 &\int_{\R^2}(\tfrac12\rho|u|^2+c_v\rho\theta)=E_0,\qquad
   \int_{\R^2}\rho u=\int_{\R^2}\rho_0u_0.\label{eq:physical-total-energy}\\
 &\int_{\R^2}|x|\rho(t)\dd x\le M+\sqrt{2m_0E_0}\,t.
    \label{eq:physical-first-moment}
\end{align}
\end{lemma}
\begin{proof}
Integrating \eqref{eq:perturbation-system}$_1$ over $\R^2$ and using it to integrate
\eqref{eq:perturbation-system}$_2$ in conservative form, we obtain
\[
 \frac{\mathrm d}{\mathrm dt}\int_{\R^2}\rho\dd x=0,\qquad
 \frac{\mathrm d}{\mathrm dt}\int_{\R^2}\rho u\dd x=0.
\]
Multiplying \eqref{eq:perturbation-system}$_2$ by $u$ and integrating over $\R^2\times(0,t)$ gives
\eqref{eq:physical-kinetic-energy}.
Integrating \eqref{eq:perturbation-system}$_3$ over $\R^2\times(0,t)$ and using
\eqref{eq:perturbation-system}$_1$ and
\begin{equation}\label{eq:viscous-dissipation-identity}
 \int_{\R^2} \bigl(2\mu|D(u)|^2+\lambda(\Div u)^2\bigr)\dd x=\mu\|\nabla u\|_2^2+(\mu+\lambda)\|\Div u\|_2^2
\end{equation}
gives \eqref{eq:physical-thermal-energy}. Adding \eqref{eq:physical-kinetic-energy} and \eqref{eq:physical-thermal-energy} gives the energy identity in \eqref{eq:physical-total-energy}.

Finally, multiplying \eqref{eq:perturbation-system}$_1$ by $|x|$ and integrating
by parts yields
\[
 \frac{\mathrm d}{\mathrm dt}\int_{\R^2}|x|\rho\dd x
 =\int_{\R^2}\rho u\cdot\frac{x}{|x|}\dd x
 \le\sqrt{m_0}\|\sqrt\rho u\|_2\le\sqrt{2m_0E_0}.
\]
Since \(|x|\le\bar x^a\), integration from \(0\) to \(t\) and
\eqref{eq:initial-data-bound} give \eqref{eq:physical-first-moment}.
\end{proof}

\begin{remark}\label{rem:mean-material-acceleration}
By \eqref{eq:perturbation-system}$_1$ and \eqref{eq:physical-total-energy}, we have
\begin{equation}\label{eq:mean-material-acceleration}
 \int_{\R^2}\rho\dot u\dd x
 =\frac{\mathrm d}{\mathrm dt}\int_{\R^2}\rho u\dd x=0.
\end{equation}
\end{remark}

Define the effective viscous flux $F$ and the vorticity $\omega$ by
\begin{equation}\label{eq:flux-definitions}
 F=(2\mu+\lambda)\Div u-P,\qquad \omega=\partial_1u^2-\partial_2u^1.
\end{equation}
We write $\nabla^\perp=(-\partial_2,\partial_1)$. Taking the divergence and the curl of
\eqref{eq:perturbation-system}$_2$ gives
\begin{equation}\label{eq:flux-elliptic}
 \Delta F=\Div(\rho\dot u),\qquad \mu\Delta\omega=\nabla^\perp\cdot(\rho\dot u).
\end{equation}
Combining \eqref{eq:gn-lp} with the div--curl estimate for $u$ and the
$L^p$ estimates for \eqref{eq:flux-elliptic} gives the following estimates
(see \cite[Lemma~2.5]{LiXin2019} and \cite[Lemma~2.6]{Wang2023}).
\begin{lemma}\label{lem:flux-estimates}
For $2\le p<\infty$, a smooth solution of \eqref{eq:perturbation-system}
satisfies
\begin{align}
 &\|\nabla u\|_p
 \le C(\|F\|_p+\|\omega\|_p+\|P\|_p),\label{eq:div-curl-estimate}\\
 &\|\nabla F\|_p+\|\nabla\omega\|_p
 \le C\|\rho\dot u\|_p,\label{eq:flux-gradient-lp}\\
 &\|F\|_p+\|\omega\|_p
 \le C\|\rho\dot u\|_2^{1-2/p}
                  (\|\nabla u\|_2+\|P\|_2)^{2/p},\label{eq:flux-interpolation}
\end{align}
where $C=C(\mu,\lambda,p)$.
\end{lemma}
For the following two lemmas, take the initial data as in Theorem~\ref{thm:target}
and assume
\[
 E_0\le1,\qquad 0\le\rho(x,t)\le2\bar\rho,\qquad \theta(x,t)\ge0
 \quad\text{on }\R^2\times[0,T].
\]
Denote the density mean by
\[
\langle w\rangle_\rho=\frac1{m_0}\int_{\R^2}\rho w\dd x.
\]
For the large-time estimates, we need the following weighted Poincar\'e inequality
with coefficient $C\sqrt{\log(e+t)}$.
We shall apply it with $w=\dot u$ and $w=\theta$ in Section~\ref{sec:three}.
\begin{lemma}\label{lem:lift-log-anchor}
Let $\rho$ be the density of a smooth solution of \eqref{eq:perturbation-system}
in $\R^2\times[0,T]$. Assume further that $2\le p<\infty$ and $w\in H^1_{\rm loc}$
with $\nabla w\in L^2$ and $\rho w\in L^1$. Then we have, for $0\le t\le T$,
\begin{equation}\label{eq:lift-log-anchor}
\|\rho^{1/p}(w-\langle w\rangle_\rho)\|_p
 \le C\sqrt{\log(e+t)}\,\|\nabla w\|_2,
\end{equation}
where $C=C(m_0,\bar\rho,M,p)$ is independent of $t$.
\end{lemma}
\begin{proof}
Fix $t\in[0,T]$ and first let $w\in C_c^\infty(\R^2)$.
We use the Fourier transform and inversion formulas of
\cite[Chapter~6, Theorem~2.4]{SteinShakarchi2003}
\[
 \widehat{w}(\xi)=\int_{\R^2}e^{-2\pi i x\cdot\xi}w(x)\dd x,\qquad
 w(x)=\int_{\R^2}e^{2\pi i x\cdot\xi}\widehat{w}(\xi)\dd\xi.
\]
Decompose $w$ as
$
 w=w_{\mathrm{low}}+w_{\mathrm{middle}}+w_{\mathrm{high}},
$
where
\[
 \widehat{w_{\mathrm{low}}}=\mathbf1_{\{|\xi|\le(1+t)^{-1}\}}\widehat w,\quad
 \widehat{w_{\mathrm{middle}}}=\mathbf1_{\{(1+t)^{-1}<|\xi|\le1\}}\widehat w,\quad
 \widehat{w_{\mathrm{high}}}=\mathbf1_{\{|\xi|>1\}}\widehat w.
\]
By Plancherel's identity and \eqref{eq:gn-lp}, we have
\begin{equation}\label{eq:high-frequency-bound}
\begin{aligned}
 \|w_{\mathrm{high}}\|_{H^1}^2
 =&\|w_{\mathrm{high}}\|_2^2+\|\nabla w_{\mathrm{high}}\|_2^2
 \le2\|\nabla w\|_2^2,
\end{aligned}
\end{equation}
and
\begin{equation}\label{eq:high-frequency-bound-second}
\begin{aligned}
 \|\rho^{1/p}w_{\mathrm{high}}\|_p
 \leq &(2\bar\rho)^{1/p}\|w_{\mathrm{high}}\|_p
 \le C\|\nabla w\|_2.
\end{aligned}
\end{equation}
It follows from Cauchy--Schwarz that
\begin{equation}\label{eq:middle-frequency}
\begin{aligned}
 \|w_{\mathrm{middle}}\|_\infty
 \leq &C\left(\int_{(1+t)^{-1}<|\xi|\le1}|\xi|^{-2}\dd\xi\right)^{1/2}
          \|\nabla w\|_2\\
 \leq &C\sqrt{\log(e+t)}\,\|\nabla w\|_2.
\end{aligned}
\end{equation}
Using $|e^{iz}-1|\le2|z|^{1/p}$ and Cauchy--Schwarz, we obtain
\[
\begin{aligned}
 |w_{\mathrm{low}}(x)-w_{\mathrm{low}}(y)|
 \leq &C|x-y|^{1/p}
       \left(\int_{|\xi|\le(1+t)^{-1}}|\xi|^{2/p-2}\dd\xi\right)^{1/2}
       \|\nabla w\|_2\\
 \leq &C\left(\frac{|x-y|}{1+t}\right)^{1/p}\|\nabla w\|_2,
\end{aligned}
\]
so Jensen's inequality, \eqref{eq:mass} and \eqref{eq:physical-first-moment} with \(E_0\le1\) give
\begin{equation}\label{eq:lift-low-frequency-bound}
\begin{aligned}
 \int_{\R^2}\rho|w_{\mathrm{low}}-\langle w_{\mathrm{low}}\rangle_\rho|^p\dd x
 \leq &\frac1{m_0}\iint_{\R^2\times\R^2}
       \rho(x)\rho(y)|w_{\mathrm{low}}(x)-w_{\mathrm{low}}(y)|^p\dd x\dd y\\
 \leq &\frac{C\|\nabla w\|_2^p}{m_0(1+t)}
       \iint_{\R^2\times\R^2}\rho(x)\rho(y)|x-y|\dd x\dd y\\
 \leq &\frac{2C\|\nabla w\|_2^p}{1+t}
       \int_{\R^2}|x|\rho\dd x\\
 \leq & C\|\nabla w\|_2^p.
\end{aligned}
\end{equation}
Jensen's inequality and \eqref{eq:high-frequency-bound}--\eqref{eq:lift-low-frequency-bound} yield
\[
\begin{aligned}
 \|\rho^{1/p}(w-\langle w\rangle_\rho)\|_p
 \leq &2\|\rho^{1/p}w_{\mathrm{high}}\|_p+2m_0^{1/p}\|w_{\mathrm{middle}}\|_\infty+\|\rho^{1/p}(w_{\mathrm{low}}-\langle w_{\mathrm{low}}\rangle_\rho)\|_p\\
 \leq &C\sqrt{\log(e+t)}\,\|\nabla w\|_2.
\end{aligned}
\]

For every $c\in\R$, we have
\begin{equation}\label{eq:constant-shift-identities}
 \nabla(w+c)=\nabla w,\qquad
 (w+c)-\langle w+c\rangle_\rho=w-\langle w\rangle_\rho.
\end{equation}
Thus \eqref{eq:lift-log-anchor}, proved above for $w\in C_c^\infty(\R^2)$,
also holds for $w\in C_c^\infty(\R^2)+\R$. Choose  $w_n\in C_c^\infty(\R^2)+\R$ such that
\[
 \int_{B_1}(w_n-w)\dd x=0,\qquad
 \|\nabla(w_n-w)\|_2\longrightarrow0.
\]
By Jensen's inequality, Poincar\'e's inequality on $B_1$ and
\eqref{eq:density-lp-local}, we have
\[
\begin{aligned}
 \|\rho^{1/p}\bigl[(w_n-\langle w_n\rangle_\rho)
                       -(w-\langle w\rangle_\rho)\bigr]\|_p
 \leq &2\|\rho^{1/p}(w_n-w)\|_p\\
 \leq &C\bigl(\|w_n-w\|_{L^2(B_1)}+\|\nabla(w_n-w)\|_2\bigr)\\
 \leq &C\|\nabla(w_n-w)\|_2\longrightarrow0.
\end{aligned}
\]
Applying \eqref{eq:lift-log-anchor} to $w_n$ and passing to the limit, we obtain
\[
\begin{aligned}
 \|\rho^{1/p}(w-\langle w\rangle_\rho)\|_p
 =&\lim_{n\to\infty}\|\rho^{1/p}(w_n-\langle w_n\rangle_\rho)\|_p\\
 \leq &C\sqrt{\log(e+t)}\lim_{n\to\infty}\|\nabla w_n\|_2\\
  =&C\sqrt{\log(e+t)}\|\nabla w\|_2.
\end{aligned}
\]
\end{proof}

To construct the weak solution of the elliptic equation in the temperature estimates, we introduce
\begin{equation}\label{eq:actual-anchored-space}
 \mathcal H_\rho
 =\{w\in H^1_{\rm loc}(\R^2):\nabla w\in L^2,\quad
                  \rho w\in L^1,\quad\langle w\rangle_\rho=0\}.
\end{equation}
Equip $\mathcal H_\rho$ with the inner product
\[
 (f,g)_{\mathcal H_\rho}=\int_{\R^2}\nabla f\cdot\nabla g\dd x,
 \qquad f,g\in\mathcal H_\rho.
\]
By Lemma~\ref{lem:lift-log-anchor}, $\mathcal H_\rho$ is a Hilbert space
with this inner product.

We next establish a weighted $L^2$ inequality for a nonnegative
weight $w$. Its application with $w=|\nabla u|^2+P^2$ provides
estimates for the auxiliary elliptic problem.

\begin{lemma}\label{lem:lift-two-measure}
Let $\rho$ be the density of a smooth solution of \eqref{eq:perturbation-system}
in $\R^2\times[0,T]$. Assume further that $0\le w\in L^1\cap L^{3/2}$ with $w\ne0$ and
$(1+|x|^2)^{1/8}w\in L^1$.
Then there exists a positive constant $C=C(m_0,\bar\rho)$ such that,
for $t\in[0,T]$ and every $v$ with
$v-\frac1{m_0}\int_{\R^2}\rho v\dd x\in\mathcal H_\rho$,
\begingroup
\medmuskip=2mu\relax
\thickmuskip=3mu\relax
\begin{equation}\label{eq:lift-two-measure}
\begin{aligned}
\int_{\R^2} w\left|v-\frac1{m_0}\int_{\R^2}\rho v\dd y\right|^2\dd x
 \le& C\|w\|_1\Bigg[1+
 \log\left(e+\frac{\|w\|_{3/2}}{\|w\|_1}\right)\\
&
 +\log\left(e+\frac{\int_{\R^2}(1+|x|^2)^{1/8}w\dd x}{\|w\|_1}
  +\frac{\int_{\R^2}(1+|x|^2)^{1/8}\rho\dd x}{m_0}\right)
 \Bigg]\|\nabla v\|_2^2.
\end{aligned}
\end{equation}
\endgroup
\end{lemma}
\begin{proof}
Fix $t\in[0,T]$. First let $v\in C_c^\infty(\R^2)$ and $L\ge2$.
We decompose
$v=v_{\mathrm{low}}+v_{\mathrm{middle}}+v_{\mathrm{high}}$, where
\[
 \widehat{v_{\mathrm{low}}}=\mathbf1_{\{|\xi|\le L^{-1}\}}\widehat v,\quad
 \widehat{v_{\mathrm{middle}}}=\mathbf1_{\{L^{-1}<|\xi|\le L\}}\widehat v,
 \quad
 \widehat{v_{\mathrm{high}}}=\mathbf1_{\{|\xi|>L\}}\widehat v.
\]
By Plancherel's identity and \eqref{eq:gn-lp} with $p=6$, we have
\[
\begin{aligned}
&\|v_{\mathrm{high}}\|_2\le CL^{-1}\|\nabla v\|_2,
 \qquad \|\nabla v_{\mathrm{high}}\|_2\le\|\nabla v\|_2,
\end{aligned}
\]
and
\[
\begin{aligned}
&\|v_{\mathrm{high}}\|_6^2
 \le C\|v_{\mathrm{high}}\|_2^{2/3}
             \|\nabla v_{\mathrm{high}}\|_2^{4/3}
 \le CL^{-2/3}\|\nabla v\|_2^2.
\end{aligned}
\]
Using H\"older's inequality and $\|\rho\|_2^2\le2\bar\rho m_0$, we obtain
\begin{equation}\label{eq:lift-two-measure-high}
\begin{aligned}
 &\int_{\R^2}w\left|v_{\mathrm{high}}
       -\frac1{m_0}\int_{\R^2}\rho v_{\mathrm{high}}\dd y\right|^2\dd x\\
 \le&2\|w\|_{3/2}\|v_{\mathrm{high}}\|_6^2
      +\frac{2\|w\|_1}{m_0^2}\|\rho\|_2^2\|v_{\mathrm{high}}\|_2^2\\
 \le& C\left(L^{-2/3}\|w\|_{3/2}+L^{-2}\|w\|_1\right)
                                      \|\nabla v\|_2^2.
\end{aligned}
\end{equation}
By Fourier inversion for \(v_{\mathrm{middle}}\) and Cauchy--Schwarz  inequality, we obtain
\[
\begin{aligned}
 \|v_{\mathrm{middle}}\|_\infty
 \le& \int_{L^{-1}<|\xi|\le L}|\widehat v(\xi)|\dd\xi
 \le \left(\int_{L^{-1}<|\xi|\le L}|\xi|^{-2}\dd\xi\right)^{1/2}
                  \||\xi|\widehat v\|_2
 \le C\sqrt{\log L}\|\nabla v\|_2.
\end{aligned}
\]
Since $\int_{\R^2}\rho\dd x=m_0$, we obtain
\begin{equation}\label{eq:lift-two-measure-middle}
\begin{aligned}
 &\int_{\R^2}w\left|v_{\mathrm{middle}}
       -\frac1{m_0}\int_{\R^2}\rho v_{\mathrm{middle}}\dd y\right|^2\dd x\le4\|w\|_1\|v_{\mathrm{middle}}\|_\infty^2
 \le C\|w\|_1\log L\,\|\nabla v\|_2^2.
\end{aligned}
\end{equation}
Using $|e^{iz}-1|\le2|z|^{1/8}$ for $z\in\R$ and
the Cauchy--Schwarz inequality, we obtain
\begin{equation}\label{eq:lift-two-measure-low}
\begin{aligned}
 |v_{\mathrm{low}}(x)-v_{\mathrm{low}}(y)|^2
 \le& C\|\nabla v\|_2^2
  \int_{|\xi|\le L^{-1}}
  \frac{|e^{2\pi i(x-y)\cdot\xi}-1|^2}{|\xi|^2}\dd\xi\\
 \le& C|x-y|^{1/4}\|\nabla v\|_2^2
                       \int_0^{L^{-1}}r^{-3/4}\dd r\\
 \le& CL^{-1/4}|x-y|^{1/4}\|\nabla v\|_2^2.
\end{aligned}
\end{equation}
By Jensen's inequality, \eqref{eq:lift-two-measure-low} and
$|x-y|^{1/4}\le|x|^{1/4}+|y|^{1/4}$, we have
\begin{equation}\label{eq:lift-two-measure-low-weighted}
\begin{aligned}
 &\int_{\R^2}w(x)\left|v_{\mathrm{low}}(x)
       -\frac1{m_0}\int_{\R^2}\rho(y)v_{\mathrm{low}}(y)\dd y\right|^2\dd x\\
 \le&\frac1{m_0}\int_{\R^2}\int_{\R^2}
 w(x)\rho(y)|v_{\mathrm{low}}(x)-v_{\mathrm{low}}(y)|^2\dd y\dd x\\
 \le& CL^{-1/4}\left[
       \int_{\R^2}(1+|x|^2)^{1/8}w\dd x
       +\frac{\|w\|_1}{m_0}\int_{\R^2}(1+|x|^2)^{1/8}\rho\dd x\right]
                                                  \|\nabla v\|_2^2.
\end{aligned}
\end{equation}
Combining \eqref{eq:lift-two-measure-high}, \eqref{eq:lift-two-measure-middle} with \eqref{eq:lift-two-measure-low-weighted}, we obtain
\begingroup
\medmuskip=2mu\relax
\thickmuskip=3mu\relax
\begin{equation}\label{eq:two-measure-frequency}
\begin{aligned}
 \int_{\R^2}w\left|v-\frac1{m_0}\int_{\R^2}\rho v\dd y\right|^2\dd x
 \le& C\|w\|_1\Bigg[1+\log L
           +L^{-2/3}\frac{\|w\|_{3/2}}{\|w\|_1}\\
 &
 +L^{-1/4}\left(
   \frac{\int_{\R^2}(1+|x|^2)^{1/8}w\dd x}{\|w\|_1}
   +\frac{\int_{\R^2}(1+|x|^2)^{1/8}\rho\dd x}{m_0}
                    \right)\Bigg]\|\nabla v\|_2^2.
\end{aligned}
\end{equation}
\endgroup
We choose
\begin{equation}\label{eq:lift-two-measure-scales}
\begin{aligned}
 L=&\left(e+\frac{\|w\|_{3/2}}{\|w\|_1}\right)^{3/2}\times\left(e+
   \frac{\int_{\R^2}(1+|x|^2)^{1/8}w\dd x}{\|w\|_1}
   +\frac{\int_{\R^2}(1+|x|^2)^{1/8}\rho\dd x}{m_0}
          \right)^4.
\end{aligned}
\end{equation}
Substituting \eqref{eq:lift-two-measure-scales} into
\eqref{eq:two-measure-frequency} proves \eqref{eq:lift-two-measure}
for $v\in C_c^\infty(\R^2)$.

Applying \eqref{eq:constant-shift-identities} to $v$, we obtain
\eqref{eq:lift-two-measure} for $v\in C_c^\infty(\R^2)+\R$.
For a general $v$, choose $v_n\in C_c^\infty(\R^2)+\R$ such that
\[
\begin{gathered}
 \int_{B_1}(v_n-v)\dd x=0,\qquad
 \nabla v_n\longrightarrow\nabla v\quad\text{in }L^2,\qquad
 v_n\longrightarrow v\quad\text{in }H^1_{\rm loc}.
\end{gathered}
\]
By \eqref{eq:density-lp-local} with $p=2$, we have
\[
 |\langle v_n-v\rangle_\rho|
 \le m_0^{-1/2}\|\sqrt\rho(v_n-v)\|_2\longrightarrow0.
\]
After passing to a subsequence, we get
\[
 v_n-\langle v_n\rangle_\rho\longrightarrow v-\langle v\rangle_\rho
 \quad\text{a.e. in }\R^2,\qquad
 \|\nabla v_n\|_2^2\longrightarrow\|\nabla v\|_2^2.
\]
Fatou's lemma and \eqref{eq:lift-two-measure} applied to $v_n$ give
\eqref{eq:lift-two-measure} for $v$.
\end{proof}

 \eqref{eq:flux-definitions} and
\eqref{eq:flux-gradient-lp} with $p=2$ give
\[
\begin{aligned}
 \|F\|_2
 &\le(2\mu+\lambda)\|\Div u\|_2
       +\|\rho\|_\infty^{1/2}\|\sqrt\rho\theta\|_2<\infty,\\
 \|\nabla F\|_2
 &\le C\|\rho\dot u\|_2
 \le C\|\rho\|_\infty^{1/2}\|\sqrt\rho\dot u\|_2<\infty.
\end{aligned}
\]
Thus $F\in H^1(\R^2)$.
In the coupled energy estimate for $F$ and $\omega$, we need to control
$\left|\int_{\R^2}F\det\nabla u\dd x\right|$ by
$\frac{\sqrt3}{8\sqrt\pi}\|\nabla F\|_2\|\nabla u\|_2^2$.
This follows from the next inequality with $f=u^1$ and $g=u^2$.
\begin{lemma}\label{lem:planar-jacobian}
If $F\in H^1$ and $f,g\in H^1_{\rm loc}$ with
$\nabla f,\nabla g\in L^2\cap L^4$, then
\begin{equation}\label{eq:lift-planar-divcurl-bound}
\left|\int_{\R^2} F(\partial_1f\partial_2g-\partial_2f\partial_1g)\dd x\right|
       \le \frac{\sqrt3}{4\sqrt\pi}\|\nabla F\|_2\|\nabla f\|_2\|\nabla g\|_2.
\end{equation}
\end{lemma}
\begin{proof}
First let $F\in C_c^\infty(\R^2)$ and choose a ball $B$ with $\supp F\subset B$.
Assume that $\phi\in H_0^1(B)$ solve
\[
 \left\{
 \begin{aligned}
  -\Delta\phi&=\partial_1f\partial_2g-\partial_2f\partial_1g
       &&\text{in }B,\\
  \phi&=0&&\text{on }\partial B.
 \end{aligned}
 \right.
\]
By \cite[Theorem~1.2]{Ge1998}, we have
\[
 \|\nabla\phi\|_{L^2(B)}^2
 \le\frac3{16\pi}\|\nabla f\|_{L^2(B)}^2\|\nabla g\|_{L^2(B)}^2.
\]
Multiplying the equation for $\phi$ by $F$ and integrating over $B$, we obtain
\[
\begin{aligned}
 \left|\int_{\R^2}F(\partial_1f\partial_2g-\partial_2f\partial_1g)\dd x\right|
 =&\left|\int_B\nabla F\cdot\nabla\phi\dd x\right|\\
 \le&\|\nabla F\|_{L^2(B)}\|\nabla\phi\|_{L^2(B)}\\
 \le&\frac{\sqrt3}{4\sqrt\pi}\|\nabla F\|_2\|\nabla f\|_2\|\nabla g\|_2.
\end{aligned}
\]
For $F\in H^1$, we can choose $F_n\in C_c^\infty(\R^2)$ such that $F_n\to F$ in $H^1$.
Since $\nabla f,\nabla g\in L^4$, H\"older's inequality gives
\[
\begin{aligned}
 &\left|\int_{\R^2}(F_n-F)(\partial_1f\partial_2g-\partial_2f\partial_1g)\dd x\right|
 \le\|F_n-F\|_2\|\nabla f\|_4\|\nabla g\|_4\longrightarrow0.
\end{aligned}
\]
Using \eqref{eq:lift-planar-divcurl-bound}, already proved for
$F_n\in C_c^\infty(\R^2)$, and $\|\nabla F_n-\nabla F\|_2\to0$, we obtain
\[
\begin{aligned}
 &\left|\int_{\R^2}F(\partial_1f\partial_2g-\partial_2f\partial_1g)\dd x\right|\\
 =&\lim_{n\to\infty}
 \left|\int_{\R^2}F_n(\partial_1f\partial_2g-\partial_2f\partial_1g)\dd x\right|\\
 \le&\frac{\sqrt3}{4\sqrt\pi}\|\nabla f\|_2\|\nabla g\|_2
       \lim_{n\to\infty}\|\nabla F_n\|_2\\
 =&\frac{\sqrt3}{4\sqrt\pi}\|\nabla F\|_2\|\nabla f\|_2\|\nabla g\|_2.
\end{aligned}
\]
\end{proof}

To estimate $\nabla\rho$ in $L^2\cap L^q$, we use the following
Beale--Kato--Majda type inequality (see \cite[Lemma~2.7]{LiXin2019}).
\begin{lemma}\label{lem:bkm}
Let $q>2$ and $\nabla u\in L^2\cap W^{1,q}$. Then
\begin{equation}\label{eq:higher-bkm}
 \|\nabla u\|_\infty
 \le C(\|\Div u\|_\infty+\|\omega\|_\infty)
           \log(e+\|\nabla^2u\|_q)+C\|\nabla u\|_2+C,
\end{equation}
where $C=C(q)$.
\end{lemma}

To estimate $\|\theta\|_\infty$, we use its density mean and the following
estimate for the logarithmic potential of $\Delta\theta$.
\begin{lemma}\label{lem:continuation-newton}
Let $f\in L^{4/3}$, $\int_{\R^2} f\dd x=0$ and
$\int_{\R^2}(1+|x|^2)^{1/8}|f|\dd x<\infty$. Then $\nabla\Delta^{-1}f\in L^2$,
\begin{equation}\label{eq:newton-potential-bound}
\begin{aligned}
 \Delta^{-1}f(x)&=\frac1{2\pi}\int_{\R^2}\log|x-y|\,f(y)\dd y\in C_0(\R^2),
\end{aligned}
\end{equation}
and
\begin{equation}\label{eq:newton-potential-bound-second}
\begin{aligned}
 \|\Delta^{-1}f\|_\infty
 &\le\frac1\pi\left(\|f\|_{4/3}
                        +\int_{\R^2}(1+|x|^2)^{1/8}|f|\dd x\right).
\end{aligned}
\end{equation}
If $v\in H^1_{\rm loc}$, $\nabla v\in L^2$ and $\Delta v=f$, then
$v-\Delta^{-1}f$ is constant.
\end{lemma}
\begin{proof}
Since $\int_{\R^2}f(y)\dd y=0$, we have
\begin{equation}\label{eq:newton-kernel-split}
\begin{aligned}
 2\pi\Delta^{-1}f(x)=&\int_{\R^2}\bigl[\log|x-y|-\log(1+|x|)\bigr]f(y)\dd y\\
 =&\int_{|x-y|<1}\log|x-y|\,f(y)\dd y
   +\int_{\R^2}\log\frac{\max\{1,|x-y|\}}{1+|x|}\,f(y)\dd y.
\end{aligned}
\end{equation}
Observe that
\begin{equation}\label{eq:newton-kernel-estimates}
 \|\log|x|\,\mathbf1_{\{|x|<1\}}\|_4
   =(3\pi/2)^{1/4}\le2,
\end{equation}
and
\begin{equation}\label{eq:newton-log-pointwise}
 \left|\log\frac{\max\{1,|x-y|\}}{1+|x|}\right|
   \le\log(2+|y|)\le2(1+|y|^2)^{1/8}.
\end{equation}
By H\"older's inequality and \eqref{eq:newton-kernel-split}--\eqref{eq:newton-log-pointwise}, we obtain
\[
\begin{aligned}
 2\pi\|\Delta^{-1}f\|_\infty
 \le&(3\pi/2)^{1/4}\|f\|_{4/3}
       +2\int_{\R^2}(1+|y|^2)^{1/8}|f(y)|\dd y,
\end{aligned}
\]
which gives \eqref{eq:newton-potential-bound-second}.
By H\"older's inequality, \eqref{eq:newton-kernel-estimates} and
translation continuity in $L^{4/3}(\R^2)$, we obtain
\[
\begin{aligned}
 &\sup_{x\in\R^2}\left|
   \int_{|y|<1}\log|y|\,[f(x+h-y)-f(x-y)]\dd y\right|\\
 \le&(3\pi/2)^{1/4}\|f(\cdot+h)-f\|_{4/3}
       \longrightarrow0\qquad(h\to0),
\end{aligned}
\]
and
\[
\begin{aligned}
 &\left|\int_{|x-y|<1}\log|x-y|\,f(y)\dd y\right|\\
 \le&(3\pi/2)^{1/4}
       \|f\|_{L^{4/3}(\{|y|>|x|-1\})}
       \longrightarrow0\qquad(|x|\to\infty).
\end{aligned}
\]
Combining Lebesgue's theorem with
\[
 \left|\log\frac{\max\{1,|x-y|\}}{1+|x|}\,f(y)\right|
 \le2(1+|y|^2)^{1/8}|f(y)|\in L^1(\R^2),
\]
we obtain
\[
 x\mapsto
 \int_{\R^2}\log\frac{\max\{1,|x-y|\}}{1+|x|}\,f(y)\dd y
 \in C_0(\R^2).
\]
By \eqref{eq:newton-kernel-split}, we have $\Delta^{-1}f\in C_0(\R^2)$.

Since $\int_{\R^2}f\dd x=0$, using $|e^{iz}-1|\le2|z|^{1/4}$ gives
\[
\begin{aligned}
 |\widehat f(\xi)|
 =&\left|\int_{\R^2}(e^{-2\pi i x\cdot\xi}-1)f(x)\dd x\right|\le2(2\pi)^{1/4}|\xi|^{1/4}
       \int_{\R^2}|x|^{1/4}|f(x)|\dd x.
\end{aligned}
\]
By Hausdorff--Young inequality \cite[(1)]{Beckner1975} with $p=4/3$, we have
\[
 \|\widehat f\|_4\le\|f\|_{4/3}.
\]
Applying Plancherel's identity to \(\nabla\Delta^{-1}f\) and using
H\"older's inequality, we obtain
\[
\begin{aligned}
 \|\nabla\Delta^{-1}f\|_2^2
 =&\frac1{4\pi^2}
   \left(\int_{|\xi|\le1}|\xi|^{-2}|\widehat f(\xi)|^2\dd\xi
         +\int_{|\xi|>1}|\xi|^{-2}|\widehat f(\xi)|^2\dd\xi\right)\\
 \le&\frac{\sqrt{2\pi}}{\pi^2}
       \left(\int_{\R^2}|x|^{1/4}|f(x)|\dd x\right)^2
       \int_{|\xi|\le1}|\xi|^{-3/2}\dd\xi+\frac1{4\pi^2}
       \left(\int_{|\xi|>1}|\xi|^{-4}\dd\xi\right)^{1/2}
       \|\widehat f\|_4^2\\
 \le&\frac{4\sqrt{2\pi}}{\pi}
       \left(\int_{\R^2}|x|^{1/4}|f(x)|\dd x\right)^2
       +\frac1{4\pi^{3/2}}\|f\|_{4/3}^2<\infty.
\end{aligned}
\]
Finally, $\Delta(v-\Delta^{-1}f)=0$ in $\mathcal D'(\R^2)$ and
$\nabla(v-\Delta^{-1}f)\in L^2$.
Choose $0\le\eta_\varepsilon\in C_c^\infty(B_\varepsilon)$ with
$\int_{\R^2}\eta_\varepsilon\dd x=1$. Then we can deduce
\[
 \Delta[\eta_\varepsilon*(v-\Delta^{-1}f)]
 =\eta_\varepsilon*\Delta(v-\Delta^{-1}f)=0.
\]
By the mean-value property \cite[Theorem~2.1]{GilbargTrudinger2001} and
Cauchy--Schwarz, we obtain, for $i=1,2$,
\[
\begin{aligned}
 |[\eta_\epsilon*\partial_i(v-\Delta^{-1}f)](x)|
 =&\frac1{\pi r^2}
       \left|\int_{B_r(x)}[\eta_\epsilon*\partial_i(v-\Delta^{-1}f)](y)\dd y\right|\\
 \le&\frac1{\sqrt\pi\,r}
       \|\eta_\epsilon*\partial_i(v-\Delta^{-1}f)\|_2\\
 \le&\frac1{\sqrt\pi\,r}
       \|\partial_i(v-\Delta^{-1}f)\|_2
       \longrightarrow0\qquad(r\to\infty).
\end{aligned}
\]
Hence we have
\[
 0=\eta_\epsilon*\nabla(v-\Delta^{-1}f)
   \longrightarrow\nabla(v-\Delta^{-1}f)
   \quad\text{in }L^2\quad(\varepsilon\rightarrow0).
\]
Thus $\nabla(v-\Delta^{-1}f)=0$, which implies that $v-\Delta^{-1}f$ is constant.
\end{proof}

The following estimate follows from
\cite[Chap~I, (4.3)]{Sobolev1938}. We use it to estimate
$F,\omega$ in $L^4$ and $\nabla u$ in $L^3$.
\begin{lemma}\label{lem:riesz-potential}
Let $1<p<2$ and $f\in L^p(\R^2)$. Then
\begin{equation}\label{eq:riesz-potential-lp}
 \|I_1f\|_{2p/(2-p)}\le C\|f\|_p,
\end{equation}
where $C=C(p)$ and
$I_1f(x)=\displaystyle\int_{\R^2}\frac{f(y)}{|x-y|}\dd y$.
\end{lemma}

\section{A Priori Estimates: Lower Order Estimates}\label{sec:three}\label{subsec:initial-rough}\label{subsec:initial-signed}\label{subsec:late-spatial}\label{subsec:material-tests}\label{subsec:thermal-lift}\label{subsec:restart-short}\label{subsec:weighted-closure}\label{subsec:threshold}
In this section, let $(\rho\ge0,u,\theta\ge0)$ be a smooth solution of
\eqref{eq:perturbation-system} on $\R^2\times[0,T]$ with smooth initial
data satisfying the transformed assumptions of Theorem~\ref{thm:target}
and with the regularity in \eqref{eq:regularity}.

Throughout this section, $C$ denotes a generic positive constant depending
only on $m_0,\bar\rho,M,\mu,\lambda,c_v,\kappa$, but not on $T$, where
$\bar\rho$ and $M$ are fixed in \eqref{eq:initial-data-bound} and $M$ need not be small.
Constants denoted by $C_T$ may also depend on $T$.
Set \(\sigma(t)=\min\{1,t\}\) for \(t\ge0\). The following proposition gives the bootstrap estimates for the density,
the velocity gradient and the temperature.

\begin{proposition}\label{prop:lower-order-closure}
There exists $\eps_0\in(0,1]$, depending only on $m_0,\bar\rho,M,\mu,\lambda,c_v,\kappa$,
such that, if $E_0\le\eps_0$ and
\begin{equation}\label{eq:initial-bootstrap-assumption}
\left\{
\begin{aligned}
 &\sup_{0\le t\le T}\|\rho(t)\|_\infty\le2\bar\rho,\\
 &\sup_{0\le t\le T}
   \bigl(\|\nabla u(t)\|_2^2+\|\sqrt\rho\,\theta(t)\|_2^2
                +\sigma(t)\|\nabla\theta(t)\|_2^2\bigr)\le 2K,
\end{aligned}\right.
\end{equation}
where
\begin{equation}\label{eq:bootstrap-K}
 K=M^2\max\left\{\frac{2(2\mu+\lambda)}{\mu},\
  1+\frac{2\bar\rho(4\mu+\lambda)}{\mu(2\mu+\lambda)^2},\
  \frac{c_v}{2\kappa}
      +\frac{\bar\rho[2(2\mu+\lambda)c_v+\kappa]}{\mu(2\mu+\lambda)^2c_v}\right\},
\end{equation}
then
\begin{equation}\label{eq:initial-bootstrap-improvement}
\left\{
\begin{aligned}
 &\sup_{0\le t\le T}\|\rho(t)\|_\infty\le\tfrac32\bar\rho,\\
 &\sup_{0\le t\le T}
   \bigl(\|\nabla u(t)\|_2^2+\|\sqrt\rho\,\theta(t)\|_2^2
                +\sigma(t)\|\nabla\theta(t)\|_2^2\bigr)\le K.
\end{aligned}\right.
\end{equation}
\end{proposition}
If \(E_0=0\), \eqref{eq:physical-kinetic-energy} and \eqref{eq:physical-total-energy}
give \(P=0\) and \(u=0\), since \(m_0>0\). Then \eqref{eq:perturbation-system}$_1$
and \eqref{eq:perturbation-system}$_3$ give \(\rho=\rho_0\) and \(\Delta\theta=0\).
Since \(\nabla\theta\in L^2\) and \(\rho\theta=0\), we have \(\theta=0\).
Thus \eqref{eq:initial-bootstrap-improvement} and the estimates in Sections~\ref{sec:three}
and~\ref{sec:four} hold for \(E_0=0\), with \(\phi=0\) in \eqref{eq:lift-elliptic-source}.
In the proofs of these estimates we take \(E_0>0\).

We first estimate $\nabla u$ and $\sqrt\rho\theta$
on the short-time interval. For the short-time velocity estimate in the
barotropic case, see \cite[Lemma~3.6]{LiXin2019}.
\begin{lemma}\label{lem:initial-rough-control}
Under the assumptions of Proposition~\ref{prop:lower-order-closure},
there exists $\tau\in(0,1]$, depending only on
$m_0,\bar\rho,M,\mu,\lambda,c_v,\kappa$, such that
\begin{equation}\label{eq:initial-rough-lemma}
 \sup_{0\le t\le\min\{T,\tau\}}
 \left(\|\nabla u(t)\|_2^2+\|\sqrt\rho\,\theta(t)\|_2^2
        +\int_0^t\|\nabla\theta(s)\|_2^2\dd s\right)\le K-m_0^2,
\end{equation}
and
\begin{equation}\label{eq:initial-rough-lemma-second}
 \int_0^{\min\{T,\tau\}}\|\sqrt\rho\,\dot u(t)\|_2^2\dd t\le C,
\end{equation}
where $C=C(m_0,\bar\rho,M,\mu,\lambda,c_v,\kappa)$.
\end{lemma}
\begin{proof}
Since \(\bar x^a\ge1\), \eqref{eq:mass} and \eqref{eq:initial-data-bound} give
\[
 \|\nabla u_0\|_2+\|\sqrt{\rho_0}\theta_0\|_2
 \le M-\|\bar x^a\rho_0\|_1
 \le M-\int_{\R^2}\rho_0\dd x=M-m_0
\]
and
\[
 \|P_0\|_2^2=\int_{\R^2}\rho_0^2\theta_0^2\dd x
 \le\bar\rho\|\sqrt{\rho_0}\theta_0\|_2^2 .
\]
Using \eqref{eq:flux-definitions} and \eqref{eq:bootstrap-K}, we obtain 
\begin{equation}\label{eq:initial-layer-constants}
\begin{aligned}
 &\left[\frac{2\|F_0\|_2^2}{\mu(2\mu+\lambda)}+2\|\omega_0\|_2^2
 +\left(\frac K{M^2}-\frac{2\bar\rho}{\mu(2\mu+\lambda)}\right)
                  \|\sqrt{\rho_0}\theta_0\|_2^2\right]^{1/2}\\
 \le&\sqrt{\frac{2(2\mu+\lambda)}{\mu}}\,\|\nabla u_0\|_2
          +\frac{\sqrt K}{M}\|\sqrt{\rho_0}\theta_0\|_2
 \le\frac{\sqrt K}{M}(M-m_0).
\end{aligned}
\end{equation}
Let $\Gamma=(\Gamma_{ij})_{i,j=1}^2$ be the matrix of fundamental solutions
to $\mu\Delta+(\mu+\lambda)\nabla\Div$ on $\R^2$. For $\mu+\lambda>0$, we can get from
\cite[(1.5)--(1.6)]{AndoJiKangKimYu2018} that
\begin{equation}\label{eq:lift-lame-fundamental-solution}
 \Gamma_{ij}(x)
 =\frac{\lambda+3\mu}{4\pi\mu(\lambda+2\mu)}\delta_{ij}\log|x|
  -\frac{\lambda+\mu}{4\pi\mu(\lambda+2\mu)}\frac{x_ix_j}{|x|^2},
 \qquad x\ne0.
\end{equation}
For $\mu+\lambda=0$, this is $\Gamma=(2\pi\mu)^{-1}(\log|x|)I$.
Hence \(|\partial_i\Gamma(x)|+|x|\,|\nabla\partial_i\Gamma(x)|\le \frac{C}{|x|}\).

We show that, for \(i=1,2\),
\begin{equation}\label{eq:lift-lame-gradient-representation}
 \partial_i u
 =\frac{\partial_i\nabla\Delta^{-1}P}{2\mu+\lambda}
     +(\partial_i\Gamma)*(\rho\dot u).
\end{equation}
By \eqref{eq:perturbation-system}$_2$ and
\((\mu\Delta+(\mu+\lambda)\nabla\Div)\Gamma=\delta_0 I\), we have
\begin{equation}\label{eq:lift-lame-differentiated}
\begin{aligned}
 (\mu\Delta+(\mu+\lambda)\nabla\Div)
 \left(\frac{\partial_i\nabla\Delta^{-1}P}{2\mu+\lambda}
           +(\partial_i\Gamma)*(\rho\dot u)\right)
 =&\partial_i\nabla P+\partial_i(\rho\dot u)\\
     =&(\mu\Delta+(\mu+\lambda)\nabla\Div)\partial_i u .
\end{aligned}
\end{equation}
The Fourier symbol of $\mu\Delta+(\mu+\lambda)\nabla\Div$ satisfies
\[
 \det\!\left[-4\pi^2\bigl(\mu|\xi|^2 I+(\mu+\lambda)\xi\otimes\xi\bigr)\right]
 =16\pi^4\mu(2\mu+\lambda)|\xi|^4>0,\qquad \xi\ne0 .
\]
Taking the Fourier transform of \eqref{eq:lift-lame-differentiated}, we obtain
\begin{equation}\label{eq:lift-lame-fourier-support}
 \operatorname{supp}\!\left(
  \widehat{\partial_i u}-\frac{\widehat{\partial_i\nabla\Delta^{-1}P}}{2\mu+\lambda}
  -\widehat{(\partial_i\Gamma)*(\rho\dot u)}\right)\subset\{0\}.
\end{equation}
By \eqref{eq:lift-lame-fourier-support} and \cite[Theorems~6.24(d) and~6.25]{Rudin1991}, there exists \(N\ge0\)
and \(c_\alpha\in\mathbb C^2\) such that
\begin{equation}\label{eq:lift-lame-dirac-expansion}
 \widehat{\partial_i u}-\frac{\widehat{\partial_i\nabla\Delta^{-1}P}}{2\mu+\lambda}
  -\widehat{(\partial_i\Gamma)*(\rho\dot u)}
 =\sum_{|\alpha|\le N}c_\alpha\,\partial_\xi^\alpha\delta_0 .
\end{equation}
Taking the inverse Fourier transform of \eqref{eq:lift-lame-dirac-expansion} and using
\(\widehat{(-2\pi\mathrm{i}x)^\alpha}=\partial_\xi^\alpha\delta_0\)
(see \cite[\S7.16, (5)]{Rudin1991}), we obtain
\[
 \partial_i u-\frac{\partial_i\nabla\Delta^{-1}P}{2\mu+\lambda}
       -(\partial_i\Gamma)*(\rho\dot u)
 =\sum_{|\alpha|\le N}c_\alpha(-2\pi\mathrm{i}x)^\alpha
 \in\mathbb C[x_1,x_2] .
\]
By \eqref{eq:initial-bootstrap-assumption}, we have
\[
 \|\partial_i\nabla\Delta^{-1}P\|_2^2
 \le\|P\|_2^2\le 2\bar\rho\,\|\sqrt\rho\,\theta\|_2^2 .
\]
Using \eqref{eq:riesz-potential-lp} with \(p=4/3\) and \eqref{eq:lift-lame-fundamental-solution}, we obtain
\[
\begin{aligned}
 \|(\partial_i\Gamma)*(\rho\dot u)\|_4
 \le C\|I_1(|\rho\dot u|)\|_4
 \le C\|\rho\dot u\|_{4/3}
 \le C(2\bar\rho m_0)^{1/4}\|\sqrt\rho\,\dot u\|_2 .
\end{aligned}
\]
Since \(\partial_i u\in L^2\) and \((L^2+L^4)\cap\mathbb C[x_1,x_2]=\{0\}\), 
 \eqref{eq:lift-lame-gradient-representation} holds.

By H\"older inequality, we arrive at
\begin{equation}\label{eq:initial-pressure-weight}
\||x|^{1/16}P\|_2
 \le(\int_{\R^2}\rho\theta^{16/7})^{7/16}(\int_{\R^2}|x|\rho^9)^{1/16}
 \le C\|\sqrt\rho\theta\|_2^{3/4}(\|\sqrt\rho\theta\|_2+\|\nabla\theta\|_2)^{1/4}.
\end{equation}
Using \eqref{eq:mean-material-acceleration} in \eqref{eq:lift-lame-gradient-representation}, we obtain, for
\(|x|>1\),
\begin{equation}\label{eq:force-kernel-cancellation}
 |\bigl((\partial_i\Gamma)*(\rho\dot u)\bigr)(x)|
 \le C|x|^{-5/4}\||y|^{1/4}\rho\dot u\|_1
     +C|x|^{-1/4}I_1(|y|^{1/4}|\rho\dot u|)(x).
\end{equation}
Applying H\"older's inequality and \eqref{eq:riesz-potential-lp} with \(p=16/15\)
to \eqref{eq:force-kernel-cancellation}, we obtain, for \(0\le t\le\min\{T,1\}\),
\[
\begin{aligned}
 \||x|^{1/16}((\partial_i\Gamma)*(\rho\dot u))\|_{L^2(\R^2\setminus B_1)}
 \le C\bigl(\||x|^{1/4}\rho\dot u\|_1+\||x|^{1/4}\rho\dot u\|_{16/15}\bigr)
 \le C\|\sqrt\rho\,\dot u\|_2,
\end{aligned}
\]
where \(|x|^{-3/16}\in L^{16}(\R^2\setminus B_1)\) and the last inequality uses
\eqref{eq:physical-first-moment} and \eqref{eq:initial-bootstrap-assumption}.

Since \((1+|x|^2)^\beta\in A_2(\R^2)\) for \(0<\beta<1\), the weighted
\(L^2\) estimate for the Riesz transforms \cite[Example~7.1.7 and Theorem~7.4.6]{Grafakos2014} gives
\begin{equation}\label{eq:weighted-riesz-pressure}
 \|(1+|x|^2)^{\beta/2}\partial_i\nabla\Delta^{-1}P\|_2
 \le C\|(1+|x|^2)^{\beta/2}P\|_2.
\end{equation}
Taking \(\beta=1/16\) and using \eqref{eq:lift-lame-gradient-representation} and
\eqref{eq:initial-pressure-weight}--\eqref{eq:weighted-riesz-pressure}, we obtain
\begin{equation}\label{eq:initial-weighted-gradient}
\begin{aligned}
 \|(1+|x|^2)^{1/32}\nabla u\|_2
 &\le C\bigl(\|\nabla u\|_2+\|\sqrt\rho\theta\|_2
                         +\|\sqrt\rho\dot u\|_2\bigr)\\
 &\quad+C\|\sqrt\rho\theta\|_2^{3/4}
        (\|\sqrt\rho\theta\|_2+\|\nabla\theta\|_2)^{1/4}.
\end{aligned}
\end{equation}

Taking \(p=4\) in \eqref{eq:density-anchor}, \eqref{eq:div-curl-estimate} and
\eqref{eq:flux-interpolation}, we deduce from \eqref{eq:physical-first-moment}
and \eqref{eq:initial-bootstrap-assumption} that
\begin{equation}\label{eq:initial-heating-interpolation} 
 \|\nabla u\|_4
 \le C\bigl(\|\sqrt\rho\,\dot u\|_2^{1/2}+\|\rho^{1/4}\theta\|_4\bigr)
 \le C\bigl(1+\|\sqrt\rho\,\dot u\|_2+\|\nabla\theta\|_2\bigr).
\end{equation}
By H\"older's inequality, \eqref{eq:initial-bootstrap-assumption} and
\eqref{eq:initial-weighted-gradient}--\eqref{eq:initial-heating-interpolation}, we have
\begin{equation}\label{eq:initial-heating-interpolation-second}
\begin{aligned}
 \|(1+|x|^2)^{1/256}\nabla u\|_{512/255}^{2}
 \le&\|(1+|x|^2)^{1/32}\nabla u\|_2^{1/4}
       \|\nabla u\|_2^{111/64}\|\nabla u\|_4^{1/64}\\
 \le&C\|\nabla u\|_2^{1/4}
       \bigl(1+\|\sqrt\rho\,\dot u\|_2^2+\|\nabla\theta\|_2^2\bigr)^{17/128}.
\end{aligned}
\end{equation}
Applying \eqref{eq:spatial-weighted-lp} with \((\eta,p)=(1/64,256)\) and
\eqref{eq:initial-heating-interpolation-second}, we obtain
\begin{equation}\label{eq:initial-sublinear-heating}
\begin{aligned}
 &\int_{\R^2}\theta\bigl(2\mu|D(u)|^2+\lambda(\Div u)^2\bigr)\dd x\\
 \le&C\|(1+|x|^2)^{-1/128}\theta\|_{256}
        \|(1+|x|^2)^{1/256}\nabla u\|_{512/255}^{2}\\
 \le&C\bigl(\|\nabla u\|_2+\|\sqrt\rho\,\theta\|_2\bigr)^{1/4}
       \bigl(1+\|\sqrt\rho\,\dot u\|_2^2+\|\nabla\theta\|_2^2\bigr)^{81/128}.
\end{aligned}
\end{equation}
It follows from \eqref{eq:density-anchor} and \eqref{eq:initial-bootstrap-assumption} that
\begin{equation}\label{eq:initial-sublinear-heating-second}
\begin{aligned}
 \int_{\R^2}\rho\theta^2|\Div u|\,\dd x
 \le&C\|\nabla u\|_2\|\sqrt\rho\,\theta\|_2^{2/3}
                \|\rho^{1/8}\theta\|_8^{4/3}\\
 \le&C\bigl(\|\nabla u\|_2+\|\sqrt\rho\,\theta\|_2\bigr)^{1/4}
        (1+\|\nabla\theta\|_2^2)^{2/3}.
\end{aligned}
\end{equation}
Multiplying \eqref{eq:perturbation-system}$_3$ by \(\theta\) and integrating over \(\R^2\),
we use \eqref{eq:perturbation-system}$_1$ and apply Young's inequality to
\eqref{eq:initial-sublinear-heating}--\eqref{eq:initial-sublinear-heating-second}. For every \(\zeta>0\), this gives
\begin{equation}\label{eq:initial-thermal-testing}
\frac{c_v}{2}\frac{\mathrm d}{\mathrm d t}\|\sqrt\rho\,\theta\|_2^2+\kappa\|\nabla\theta\|_2^2
 \le\zeta\bigl(\|\sqrt\rho\,\dot u\|_2^2+\|\nabla\theta\|_2^2\bigr)+C_\zeta.
\end{equation}
Multiplying $\rho\dot u=\nabla F+\mu\nabla^\perp\omega$ by $\dot u$
and integrating by parts, we use \eqref{eq:perturbation-system}$_1$ and
\eqref{eq:perturbation-system}$_3$ to obtain
\begin{equation}\label{eq:initial-flux-identity}
\begin{aligned}
 &\frac{\mathrm d}{\mathrm d t}
   \left(\frac{\|F\|_2^2}{2(2\mu+\lambda)}+\frac{\mu}{2}\|\omega\|_2^2\right)
   +\|\sqrt\rho\,\dot u\|_2^2\\
 =&-\frac{\mu}{2}\int_{\R^2}(\Div u)\omega^2\dd x
    +\frac{1}{2(2\mu+\lambda)}\int_{\R^2}(\Div u)F^2\dd x
    -\int_{\R^2}F\partial_i u^j\partial_j u^i\dd x\\
 &+\frac{1+c_v^{-1}}{2\mu+\lambda}\int_{\R^2}FP\Div u\dd x
    +\frac{\kappa}{(2\mu+\lambda)c_v}\int_{\R^2}\nabla F\cdot\nabla\theta\dd x\\
 &-\frac{1}{(2\mu+\lambda)c_v}
       \int_{\R^2}F\bigl(2\mu|D(u)|^2+\lambda(\Div u)^2\bigr)\dd x .
\end{aligned}
\end{equation}
By \eqref{eq:div-curl-estimate} with \(p=3\), the cubic terms in
\eqref{eq:initial-flux-identity} are bounded by
\(C(\|F\|_3^3+\|\omega\|_3^3+\|P\|_3^3)\).
Equations \eqref{eq:gn-lp} and \eqref{eq:flux-gradient-lp}, with \(p=3\) and \(p=2\)
respectively, and Young's inequality give under \eqref{eq:initial-bootstrap-assumption}
\begin{equation}\label{eq:initial-flux-products}
 C\bigl(\|F\|_3^3+\|\omega\|_3^3+\|P\|_3^3\bigr)
 \le\frac{1}{4}\|\sqrt\rho\,\dot u\|_2^2
    +C\left(\frac{\|F\|_2^2}{2(2\mu+\lambda)}
             +\frac{\mu}{2}\|\omega\|_2^2\right)^2+C\|P\|_3^3.
\end{equation}
Multiplying \eqref{eq:flux-elliptic}$_1$ by \(F\) and integrating by parts,
we deduce from \eqref{eq:initial-bootstrap-assumption} that
\[
 \|\nabla F\|_2^2
 =\int_{\R^2}\rho\dot u\cdot\nabla F\dd x
 \le\sqrt{2\bar\rho}\,\|\sqrt\rho\,\dot u\|_2\|\nabla F\|_2.
\]
It follows from Cauchy--Schwarz and Young's inequality that
\begin{equation}\label{eq:initial-flux-temperature-cross}
\begin{aligned}
 \frac{\kappa}{(2\mu+\lambda)c_v}\left|\int_{\R^2}\nabla F\cdot\nabla\theta\dd x\right|
 \le&\frac{\kappa\sqrt{2\bar\rho}}{(2\mu+\lambda)c_v}
      \|\sqrt\rho\,\dot u\|_2\|\nabla\theta\|_2\\
 \le&\frac{1}{4}\|\sqrt\rho\,\dot u\|_2^2
    +\frac{2\kappa^2\bar\rho}{(2\mu+\lambda)^2c_v^2}\|\nabla\theta\|_2^2.
\end{aligned}
\end{equation}
By H\"older's inequality, \eqref{eq:density-anchor} with \(p=6\),
and \eqref{eq:initial-bootstrap-assumption}, we get
\begin{equation}\label{eq:initial-cubic-pressure}
\begin{aligned}
 \|P\|_3^3
 \le&4\bar\rho^2\int_{\R^2}\rho\theta^3\dd x
 \le C\|\sqrt\rho\,\theta\|_2^{3/2}
       \bigl(\|\sqrt\rho\,\theta\|_2+\|\nabla\theta\|_2\bigr)^{3/2}\le\|\nabla\theta\|_2^2+C.
\end{aligned}
\end{equation}
By \eqref{eq:initial-cubic-pressure} and Young's inequality,
the term \(C\|P\|_3^3\) in \eqref{eq:initial-flux-products} is bounded by \(\|\nabla\theta\|_2^2+C\).
Together with \eqref{eq:initial-flux-temperature-cross}, this gives from \eqref{eq:initial-flux-identity}
\begin{equation}\label{eq:initial-flux-budget}
\frac{d}{dt}\left(\frac{\|F\|_2^2}{2(2\mu+\lambda)}+\frac\mu2\|\omega\|_2^2\right)+\tfrac12 \|\sqrt\rho\dot u\|_2^2\le\left(1+\frac{2\kappa^2\bar\rho}{(2\mu+\lambda)^2c_v^2}\right)\|\nabla\theta\|_2^2+C.
\end{equation}

Choose
\begin{equation}\label{eq:initial-absorption-choice}
 \zeta=\frac{\kappa(2\mu+\lambda)^2c_v^2}
              {16\bigl(2(2\mu+\lambda)^2c_v^2+2\kappa^2\bar\rho\bigr)}.
\end{equation}
Then we have
\[
 \left(\frac4\kappa+\frac{4\kappa\bar\rho}{(2\mu+\lambda)^2c_v^2}\right)\zeta
 =\frac18,\qquad \zeta\le\frac\kappa4.
\]
Multiplying \eqref{eq:initial-thermal-testing} by $\frac4\kappa+\frac{4\kappa\bar\rho}{(2\mu+\lambda)^2c_v^2}$ and adding
\eqref{eq:initial-flux-budget}, we obtain
\begin{equation}\label{eq:initial-rough-budget}
\begin{aligned}
 &\frac{d}{dt}\left(
 \frac{\|F\|_2^2}{2(2\mu+\lambda)}+\frac\mu2\|\omega\|_2^2
              +\left(\frac{2c_v}{\kappa}+\frac{2\kappa\bar\rho}{(2\mu+\lambda)^2c_v}\right)\|\sqrt\rho\theta\|_2^2\right)
+\tfrac14\|\sqrt\rho\dot u\|_2^2+\|\nabla\theta\|_2^2\le C.
\end{aligned}
\end{equation}
Integrating \eqref{eq:initial-rough-budget} over \([0,t]\) and using
\eqref{eq:initial-layer-constants}, we obtain, for \(0\le t\le\min\{T,1\}\),
\begin{equation}\label{eq:initial-flux-and-dissipation}
 \sup_{0\le s\le t}
 \bigl(\|F(s)\|_2^2+\|\omega(s)\|_2^2+\|\sqrt\rho\,\theta(s)\|_2^2\bigr)\le C,
\end{equation}
and
\begin{equation}\label{eq:initial-flux-and-dissipation-second}
 \int_0^t\bigl(\|\sqrt\rho\,\dot u(s)\|_2^2+\|\nabla\theta(s)\|_2^2\bigr)\dd s\le C,
\end{equation}
which gives \eqref{eq:initial-rough-lemma-second}.

By \eqref{eq:flux-interpolation} with \(p=3\), \eqref{eq:initial-sublinear-heating}--\eqref{eq:initial-sublinear-heating-second}
and \eqref{eq:initial-cubic-pressure}, H\"older's inequality and
\eqref{eq:initial-flux-and-dissipation}--\eqref{eq:initial-flux-and-dissipation-second} give
\[
\begin{aligned}
 &\int_0^t\left(\|F\|_3^3+\|\omega\|_3^3+\|P\|_3^3
      +\int_{\R^2}\bigl[\theta(2\mu|D(u)|^2+\lambda(\Div u)^2)
                     +\rho\theta^2|\Div u|\bigr]\dd x\right)\dd s\\
 \le&Ct^{1/4}\left(\int_0^t
           (1+\|\sqrt\rho\,\dot u\|_2^2+\|\nabla\theta\|_2^2)\dd s\right)^{3/4}\\
 \le& Ct^{1/4}.
\end{aligned}
\]
Similar to \eqref{eq:initial-flux-temperature-cross}, we obtain
\[
 \frac{\kappa}{(2\mu+\lambda)c_v}
       \left|\int_{\R^2}\nabla F\cdot\nabla\theta\dd x\right|
 \le\|\sqrt\rho\,\dot u\|_2^2+
   \frac{\bar\rho\kappa^2}{2(2\mu+\lambda)^2c_v^2}\|\nabla\theta\|_2^2 .
\]
Since \(2\mu+\lambda\ge\mu\) and
\(\|P\|_2^2\le2\bar\rho\|\sqrt\rho\,\theta\|_2^2\), we conclude from  \eqref{eq:flux-definitions} that
\[
\begin{aligned}
 &\|\nabla u\|_2^2
  =\frac{\|F+P\|_2^2}{(2\mu+\lambda)^2}+\|\omega\|_2^2
 \le\frac{2\|F\|_2^2}{\mu(2\mu+\lambda)}+2\|\omega\|_2^2
        +\frac{4\bar\rho}{(2\mu+\lambda)^2}\|\sqrt\rho\,\theta\|_2^2 .
\end{aligned}
\]
Observe from \eqref{eq:bootstrap-K} that
\[
 \frac K{M^2}-\frac{2\bar\rho}{\mu(2\mu+\lambda)}
 \ge 1+\frac{4\bar\rho}{(2\mu+\lambda)^2},
\]
and
\[
 \frac{2\kappa}{c_v}\left(\frac K{M^2}-\frac{2\bar\rho}{\mu(2\mu+\lambda)}\right)
 \ge 1+\frac{2\bar\rho\kappa^2}{\mu(2\mu+\lambda)^2c_v^2}.
\]
Multiplying \eqref{eq:perturbation-system}$_3$ by
\(\frac{2}{c_v}\bigl(\frac K{M^2}-\frac{2\bar\rho}{\mu(2\mu+\lambda)}\bigr)\theta\),
integrating over \(\R^2\), adding \(4/\mu\) times \eqref{eq:initial-flux-identity},
and integrating over \([0,t]\), we obtain from
\eqref{eq:div-curl-estimate} with \(p=3\) and \eqref{eq:initial-layer-constants} that
\[
\begin{aligned}
 &\|\nabla u(t)\|_2^2+\|\sqrt\rho\,\theta(t)\|_2^2
                   +\int_0^t\|\nabla\theta(s)\|_2^2\dd s
 \le\frac K{M^2}(M-m_0)^2+Ct^{1/4}
 \le K-2m_0^2+Ct^{1/4},
\end{aligned}
\]
where \(M\ge m_0\) by \eqref{eq:initial-data-bound} and \(\frac{K}{M^2}\ge2\) by
\eqref{eq:bootstrap-K}. Taking \(C\tau^{1/4}\le m_0^2\)
gives \eqref{eq:initial-rough-lemma}.
\end{proof}

By \eqref{eq:regularity} and \eqref{eq:gn-lp}, \(u\in L^\infty(0,T;L^\infty)\).
Since \(q>2\), Sobolev embedding gives \(u\in L^1(0,T;W^{1,\infty}(\R^2))\).
Let \(X\) be the flow of \(u\) (see~\cite[pp.~1192, (ODE)]{AmbrosioCrippa2014}) with
\begin{equation}\label{eq:particle-flow}
 \partial_tX(t,x)=u(t,X(t,x))\quad(0<t\le T),\qquad X(0,x)=x,\quad x\in\R^2.
\end{equation}
For each \(t\in[0,T]\), \(X(t,\cdot)\) is a bi-Lipschitz homeomorphism of \(\R^2\)
(see~\cite[pp.~1195--1198]{AmbrosioCrippa2014}).
We integrate \eqref{eq:perturbation-system}$_1$ along \(X\) to improve the density estimate from $2\bar\rho$ to
$3\bar\rho/2$.
\begin{lemma}\label{lem:initial-density-control}
Under the assumptions of Proposition~\ref{prop:lower-order-closure}, we have
\[
 \sup_{0\le t\le\min\{T,\tau\}}\|\rho(t)\|_\infty
 \le\tfrac32\bar\rho.
\]
\end{lemma}
\begin{proof}
Applying \(\Delta^{-1}\Div\) to \eqref{eq:fns}$_2$ and using
\eqref{eq:flux-definitions}, we obtain
\[
 \partial_t\Delta^{-1}\Div(\rho u)
 +\sum_{i,j}\Delta^{-1}\partial_i\partial_j(\rho u_i u_j)
 =(2\mu+\lambda)\Div u-P=F,
\]
which  yields from \eqref{eq:perturbation-system}$_1$ and \eqref{eq:flux-definitions} that
\begin{equation}\label{eq:signed-potential-identity}
\begin{aligned}
 &(\partial_t+u\cdot\nabla)
   \bigl[(2\mu+\lambda)\log\rho+\Delta^{-1}\Div(\rho u)\bigr]\\
 =&-F-P+(\partial_t+u\cdot\nabla)\Delta^{-1}\Div(\rho u)\\
 =&-P-\sum_{i,j}\bigl[\Delta^{-1}\partial_i\partial_j(\rho u_i u_j)
                   -u_i\Delta^{-1}\partial_i\partial_j(\rho u_j)\bigr].
\end{aligned}
\end{equation}
By \eqref{eq:initial-bootstrap-assumption} and \eqref{eq:physical-total-energy}, we have
\[
\begin{aligned}
 \|\rho u\|_1&\le\sqrt{2m_0E_0},\qquad
 \|\rho u\|_2\le\sqrt{4\bar\rho E_0},
\end{aligned}
\]
Taking \(p=4,6\) in \eqref{eq:density-anchor} and using
\eqref{eq:physical-first-moment} and \eqref{eq:initial-rough-lemma} yields
\[
 \|\rho u\|_4+\|\rho u\|_6
 \le C\left(\|\rho u\|_1+\|\nabla u\|_2\right)\le C.
\]
By H\"older's inequality, we obtain
\begin{equation}\label{eq:momentum-potential-bound}
\begin{aligned}
 |\Delta^{-1}\Div(\rho u)|
 \le&\frac{1}{2\pi}\int_{\R^2}
   \frac{|\rho u|}{|x-y|}\dd x\\
 \le&\frac{1}{2\pi}\left(
   \int_{|x-y|\le h}\frac{|\rho u|}{|x-y|}\dd x
   +h^{-1}\int_{|x-y|>h}|\rho u|\dd x\right)\\
 \le& C\bigl(h^{1/2}\|\rho u\|_4+h^{-1}\|\rho u\|_1\bigr)\\
 \le& C\|\rho u\|_4^{2/3}\|\rho u\|_1^{1/3}\\
 \le& CE_0^{1/6},
\end{aligned}
\end{equation}
where \(h=(\|\rho u\|_1/\|\rho u\|_4)^{2/3}\) for \(\rho u\ne0\).
The commutator estimate \cite[Theorem~B.1(2)]{FridMarroquinNariyoshi2024} gives
\[
\begin{aligned}
 &\left\|\nabla\sum_{i,j}
   \bigl[\Delta^{-1}\partial_i\partial_j(\rho u_i u_j)
         -u_i\Delta^{-1}\partial_i\partial_j(\rho u_j)\bigr]
   \right\|_{4p/(4+p)}
 \le C\|\nabla u\|_4\|\rho u\|_p,\qquad p\in\{2,6\}.
\end{aligned}
\]
It follows from \eqref{eq:gn-infty} with \((q,r)=(4,12/5)\) and
Sobolev's inequality \cite[Chap.~I, (4.3)]{Sobolev1938} that
\begin{equation}\label{eq:initial-commutator-bound}
\begin{aligned}
 \left\|\sum_{i,j}\bigl[\Delta^{-1}\partial_i\partial_j(\rho u_i u_j)
                   -u_i\Delta^{-1}\partial_i\partial_j(\rho u_j)\bigr]\right\|_\infty
 \le&C\|\nabla u\|_4\|\rho u\|_2^{1/4}\|\rho u\|_6^{3/4}\\
 \le& CE_0^{1/8}\|\nabla u\|_4,
\end{aligned}
\end{equation}
By \eqref{eq:density-anchor}, \eqref{eq:physical-first-moment},
\eqref{eq:div-curl-estimate}, \eqref{eq:flux-interpolation},
\eqref{eq:initial-rough-lemma}--\eqref{eq:initial-rough-lemma-second}
and H\"older's inequality with \(\tau\le1\), we have
\begin{equation}\label{eq:initial-gradient-time}
\begin{aligned}
 \int_0^{\min\{T,\tau\}}\|\nabla u\|_4\dd t
 \le& C\int_0^{\min\{T,\tau\}}\bigl(\|\sqrt\rho\dot u\|_2^{1/2}
            +\|\sqrt\rho\theta\|_2+\|\nabla\theta\|_2\bigr)\dd t\\
 \le&C\left(\int_0^{\min\{T,\tau\}}
                   \|\sqrt\rho\dot u\|_2^2\dd t\right)^{1/4}\\
 &+C\left(\sup_{0\le t\le\min\{T,\tau\}}
                   \|\sqrt\rho\theta(t)\|_2^2
          +\int_0^{\min\{T,\tau\}}\|\nabla\theta\|_2^2\dd t\right)^{1/2}\\
 \le&C.
\end{aligned}
\end{equation}
Solving \eqref{eq:perturbation-system}$_1$ along the flow in \eqref{eq:particle-flow}, we conclude
\[
 \rho(t,X(t,x))=\rho_0(x)\exp\!\left(-\int_0^t\Div u(s,X(s,x))\dd s\right).
\]
For \(\rho_0(x)>0\), integrating \eqref{eq:signed-potential-identity} along \(X(\cdot,x)\) over \([0,t]\) and using
\(P\ge0\) and \eqref{eq:momentum-potential-bound}--\eqref{eq:initial-gradient-time}, we obtain
\begin{equation}\label{eq:initial-density-gain}
\begin{aligned}
 \log\frac{\rho(t,X(t,x))}{\rho_0(x)}
 \le&\frac{\|\Delta^{-1}\Div(\rho_0u_0)\|_\infty
             +\|\Delta^{-1}\Div(\rho u)(t)\|_\infty}{2\mu+\lambda}\\
 &\quad+\frac1{2\mu+\lambda}\int_0^t
 \left\|\sum_{i,j}\bigl[\Delta^{-1}\partial_i\partial_j(\rho u_i u_j)
                   -u_i\Delta^{-1}\partial_i\partial_j(\rho u_j)\bigr](s)\right\|_\infty\dd s\\
 \le& CE_0^{1/8}.
\end{aligned}
\end{equation}
For \(C\) in \eqref{eq:initial-density-gain}, we can choose \(\eps_0\) so that
\[
 CE_0^{1/8}\le C\eps_0^{1/8}\le\log\frac{3}{2}.
\]
\eqref{eq:initial-data-bound} and
\eqref{eq:initial-density-gain} give
\[
 \|\rho(t)\|_\infty=\|\rho(t,X(t,\cdot))\|_\infty
 \le\tfrac32\|\rho_0\|_\infty\le\tfrac32\bar\rho,
 \qquad 0\le t\le\min\{T,\tau\}.\qedhere
\]
\end{proof}

We next use \eqref{eq:initial-rough-lemma} to control
the time-weighted temperature gradient in \eqref{eq:initial-bootstrap-improvement}
for short times.
\begin{lemma}\label{lem:initial-time-weighted}
Under the assumptions of Proposition~\ref{prop:lower-order-closure},
there is a choice of \(\tau\) in Lemma~\ref{lem:initial-rough-control}
such that
\begin{equation}\label{eq:initial-weighted-lemma}
\begin{aligned}
 &\int_0^{\min\{T,\tau\}}\bigl(\|\nabla u\|_2^2+
                    \|\sqrt\rho\theta\|_2^2\bigr)\dd t\le CE_0^{2/3},
\end{aligned}
\end{equation}
and
\begin{equation}\label{eq:initial-weighted-lemma-second}
\begin{aligned}
 &\sup_{0\le t\le \min\{T,\tau\}}t\bigl(\|\nabla u(t)\|_2^2+
           \|\sqrt\rho\theta(t)\|_2^2+\|P(t)\|_2^2\bigr)\\
 &\qquad+\int_0^{\min\{T,\tau\}}t\bigl(\|\sqrt\rho\dot u\|_2^2+
                             \|\nabla\theta\|_2^2\bigr)\dd t
                                                    \le CE_0^{1/12},
\end{aligned}
\end{equation}
where \(C=C(m_0,\bar\rho,M,\mu,\lambda,c_v,\kappa)\). Moreover,
\begin{equation}\label{eq:initial-weighted-temperature-gradient}
 \sup_{0\le t\le\min\{T,\tau\}}
 \bigl(\|\nabla u(t)\|_2^2+\|\sqrt\rho\,\theta(t)\|_2^2
                +t\|\nabla\theta(t)\|_2^2\bigr)\le K.
\end{equation}
\end{lemma}
\begin{proof}
Let \(0\le t\le\min\{T,\tau\}\).
By \eqref{eq:density-anchor}, \eqref{eq:physical-total-energy}
and \eqref{eq:physical-first-moment}, we have
\begin{equation}\label{eq:physical-initial-interpolation}
 \|\sqrt\rho\theta\|_2^2
 \le\left(\int_{\R^2}\rho\theta\dd x\right)^{2/3}
       \left(\int_{\R^2}\rho\theta^4\dd x\right)^{1/3}
 \le CE_0^{2/3}(E_0+\|\nabla\theta\|_2)^{4/3}.
\end{equation}
Using \eqref{eq:initial-rough-lemma},
\eqref{eq:physical-initial-interpolation} and H\"older's inequality, we obtain
\begin{equation}\label{eq:initial-temperature-integral}
\begin{aligned}
 \int_0^{\min\{T,\tau\}}\|\sqrt\rho\theta\|_2^2\dd t
 \le&CE_0^2\min\{T,\tau\}+CE_0^{2/3}
       \left(\int_0^{\min\{T,\tau\}}\|\nabla\theta\|_2^2\dd t\right)^{2/3}\\
 \le& CE_0^{2/3}.
\end{aligned}
\end{equation}
Applying Young's inequality to \eqref{eq:physical-kinetic-energy} and combining
 \eqref{eq:physical-total-energy} with \eqref{eq:initial-bootstrap-assumption}, one has
\begin{equation}\label{eq:small-mechanical-integral}
 \int_0^{\min\{T,\tau\}}\|\nabla u\|_2^2\dd t
 \le\frac{2E_0}{\mu}
       +\frac{4\bar\rho}{\mu^2}\int_0^{\min\{T,\tau\}}\|\sqrt\rho\theta\|_2^2\dd t.
\end{equation}
Together with \eqref{eq:initial-temperature-integral}, this gives \eqref{eq:initial-weighted-lemma}.

Applying H\"older's inequality in time to
\eqref{eq:initial-sublinear-heating}--\eqref{eq:initial-sublinear-heating-second}
and using \eqref{eq:initial-rough-lemma}--\eqref{eq:initial-rough-lemma-second}
and \eqref{eq:initial-weighted-lemma}, we obtain
\begin{equation}\label{eq:small-heating-integral}
\begin{aligned}
 &\int_0^{\min\{T,\tau\}}\!\int_{\R^2}\theta
           \bigl(2\mu|D(u)|^2+\lambda(\Div u)^2\bigr)\dd x\dd t\\
 \le& C\left(\int_0^{\min\{T,\tau\}}
            (\|\nabla u\|_2+\|\sqrt\rho\theta\|_2)^2\dd t\right)^{1/8}\left(\int_0^{\min\{T,\tau\}}
       (1+\|\sqrt\rho\dot u\|_2^2+\|\nabla\theta\|_2^2)\dd t\right)^{81/128}\\
       \le& CE_0^{1/12},
\end{aligned}
\end{equation}
and
\begin{equation}\label{eq:small-heating-integral-second}
\begin{aligned}
 &\int_0^{\min\{T,\tau\}}\!\int_{\R^2}\rho\theta^2|\Div u|\dd x\dd t\\
 \le& C\left(\int_0^{\min\{T,\tau\}}
            (\|\nabla u\|_2+\|\sqrt\rho\theta\|_2)^2\dd t\right)^{1/8}
       \left(\int_0^{\min\{T,\tau\}}(1+\|\nabla\theta\|_2^2)\dd t\right)^{2/3}\\
 \le& CE_0^{1/12}.
\end{aligned}
\end{equation}
Multiplying \eqref{eq:perturbation-system}$_3$ by $s\theta$ and integrating over $\R^2\times(0,t)$,
we obtain from \eqref{eq:perturbation-system}$_1$ that
\begin{equation}\label{eq:time-weighted-thermal-identity}
\begin{aligned}
 &\frac{c_v}{2}t\|\sqrt\rho\theta(t)\|_2^2
       +\kappa\int_0^t s\|\nabla\theta(s)\|_2^2\dd s\\
 =&\frac{c_v}{2}\int_0^t\|\sqrt\rho\theta(s)\|_2^2\dd s
       +\int_0^t s\int_{\R^2}(\theta \bigl(2\mu|D(u)|^2+\lambda(\Div u)^2\bigr)-\rho\theta^2\Div u)\dd x\dd s.
\end{aligned}
\end{equation}
Substituting \eqref{eq:initial-temperature-integral} and
\eqref{eq:small-heating-integral}--\eqref{eq:small-heating-integral-second}
into \eqref{eq:time-weighted-thermal-identity}, we obtain
\[
\begin{aligned}
 \frac{c_v}{2}t\|\sqrt\rho\theta(t)\|_2^2
       +\kappa\int_0^t s\|\nabla\theta(s)\|_2^2\dd s
 \le& CE_0^{2/3}+CtE_0^{1/12}\le CE_0^{1/12},
\end{aligned}
\]
which deduce that
\begin{equation}\label{eq:small-thermal-budget}
\sup_{0\le t\le \min\{T,\tau\}}t\|\sqrt\rho\theta(t)\|_2^2
       +\int_0^{\min\{T,\tau\}}t\|\nabla\theta\|_2^2\dd t
       \le CE_0^{1/12}.
\end{equation}
Integrating \eqref{eq:initial-cubic-pressure} over
\([0,\min\{T,\tau\}]\) and using
\eqref{eq:initial-rough-lemma}, \eqref{eq:initial-temperature-integral} and H\"older's inequality in time, we obtain
\begin{equation}\label{eq:initial-pressure-time-integral-second}
\begin{aligned}
 \int_0^{\min\{T,\tau\}}\|P\|_3^3\dd t
 \le&C\int_0^{\min\{T,\tau\}}\|\sqrt\rho\theta\|_2^3\dd t\\
 &+C\left(\int_0^{\min\{T,\tau\}}\|\sqrt\rho\theta\|_2^6\dd t\right)^{1/4}
      \left(\int_0^{\min\{T,\tau\}}\|\nabla\theta\|_2^2\dd t\right)^{3/4}\\
 \le& CE_0^{1/6}.
\end{aligned}
\end{equation}
By \eqref{eq:flux-definitions}, \eqref{eq:initial-bootstrap-assumption},
\eqref{eq:initial-flux-and-dissipation} and \eqref{eq:initial-weighted-lemma}, we arrive at
\begin{equation}\label{eq:initial-flux-time-integrals}
\begin{aligned}
 \int_0^{\min\{T,\tau\}}\bigl(\|F\|_2^2+\|\omega\|_2^2+\|F\|_2^4+\|\omega\|_2^4\bigr)\dd t
 \le& C\int_0^{\min\{T,\tau\}}\bigl(\|F\|_2^2+\|\omega\|_2^2\bigr)\dd t\\
 \le& CE_0^{2/3}.
\end{aligned}
\end{equation}
Substituting \eqref{eq:initial-flux-products} and \eqref{eq:initial-flux-temperature-cross}
into \eqref{eq:initial-flux-identity} gives
\begin{equation}\label{eq:initial-flux-cubic-pressure}
\begin{aligned}
 &\frac{d}{dt}\left(\frac{\|F\|_2^2}{2(2\mu+\lambda)}
                         +\frac\mu2\|\omega\|_2^2\right)
                         +\tfrac12\|\sqrt\rho\dot u\|_2^2
 \le C\bigl(\|F\|_2^4+\|\omega\|_2^4+\|P\|_3^3\bigr)
       +\frac{2\kappa^2\bar\rho}{(2\mu+\lambda)^2c_v^2}\|\nabla\theta\|_2^2.
\end{aligned}
\end{equation}
Multiplying \eqref{eq:initial-flux-cubic-pressure} by $s$ and integrating over $[0,t]$, we obtain
\begin{equation}\label{eq:initial-weighted-flux}
\begin{aligned}
 &t\left(\frac{\|F(t)\|_2^2}{2(2\mu+\lambda)}
                         +\frac\mu2\|\omega(t)\|_2^2\right)
       +\tfrac12\int_0^t s\|\sqrt\rho\dot u\|_2^2\dd s\\
 \le& C\int_0^t\bigl(\|F\|_2^2+\|\omega\|_2^2\bigr)\dd s
       +C\int_0^t s\bigl(\|F\|_2^4+\|\omega\|_2^4+\|P\|_3^3\bigr)\dd s\\
 &+\frac{2\kappa^2\bar\rho}{(2\mu+\lambda)^2c_v^2}
                    \int_0^t s\|\nabla\theta\|_2^2\dd s.
\end{aligned}
\end{equation}
Substituting \eqref{eq:small-thermal-budget}, \eqref{eq:initial-pressure-time-integral-second} and \eqref{eq:initial-flux-time-integrals} into \eqref{eq:initial-weighted-flux}, we obtain
\begin{equation}\label{eq:initial-weighted-velocity}
 \sup_{0\le t\le \min\{T,\tau\}}t\bigl(\|F(t)\|_2^2+\|\omega(t)\|_2^2\bigr)
 +\int_0^{\min\{T,\tau\}}t\|\sqrt\rho\dot u\|_2^2\dd t\le CE_0^{1/12}.
\end{equation}
Using \eqref{eq:div-curl-estimate} with \(p=2\) and \eqref{eq:initial-bootstrap-assumption},
we deduce from
\eqref{eq:small-thermal-budget} and \eqref{eq:initial-weighted-velocity} that \eqref{eq:initial-weighted-lemma-second} holds .

We next estimate \(t\|\nabla\theta(t)\|_2^2\).
Applying \(\partial_t+\Div(u\,\cdot)\) to \eqref{eq:perturbation-system}$_2$
and using \eqref{eq:perturbation-system}$_1$, we get
\begin{equation}\label{eq:lift-material-momentum}
\begin{aligned}
 \rho\ddot u^j-\mu\Delta\dot u^j-(\mu+\lambda)\partial_j\Div\dot u
 =&-\partial_j(\rho\dot\theta)+\partial_i(P\partial_ju^i)\\
 &+\mu\partial_i\bigl[(\Div u)\partial_iu^j
             -\partial_iu^k\partial_ku^j-\partial_ku^i\partial_ku^j\bigr]\\
 &+(\mu+\lambda)\partial_j
             \bigl[(\Div u)^2-\partial_ku^l\partial_lu^k\bigr]
             -(\mu+\lambda)\partial_i\bigl[(\Div u)\partial_ju^i\bigr].
\end{aligned}
\end{equation}
Multiplying \eqref{eq:lift-material-momentum} by \(\dot u^j\) and integrating by parts, we have
\[
 2\left|\int_{\R^2}\rho\dot\theta\,\Div\dot u\dd x\right|
 \le \frac{\mu}{2}\|\nabla\dot u\|_2^2
        +\frac{8\bar\rho}{\mu}\|\sqrt\rho\,\dot\theta\|_2^2,
\]
which implies from Young's inequality that
\begin{equation}\label{eq:material-acceleration-products}
 \frac{\dd}{\dd t}\|\sqrt\rho\,\dot u\|_2^2+\mu\|\nabla\dot u\|_2^2
 \le \frac{8\bar\rho}{\mu}\|\sqrt\rho\,\dot\theta\|_2^2
       +C\bigl(\|\nabla u\|_4^4+\|P\|_4^4\bigr).
\end{equation}
Multiplying \eqref{eq:perturbation-system}$_3$ by \(\dot\theta\) and integrating
by parts, we obtain
\begin{equation}\label{eq:higher-cross-identity}
\begin{aligned}
 &\frac{\dd}{\dd t}\left[
       \frac\kappa2\|\nabla\theta\|_2^2
       -\int_{\R^2}\theta\bigl(2\mu|D(u)|^2+\lambda(\Div u)^2\bigr)\dd x\right]
       +c_v\|\sqrt\rho\dot\theta\|_2^2\\
 =&\kappa\int_{\R^2}\Delta\theta\,u\cdot\nabla\theta\dd x
          -\int_{\R^2} P\Div u\,\dot\theta\dd x\\
 &-\int_{\R^2}\theta(\partial_t+u\cdot\nabla+\Div u)
            \bigl(2\mu|D(u)|^2+\lambda(\Div u)^2\bigr)\dd x.
\end{aligned}
\end{equation}
It follows from Young's inequality and \eqref{eq:initial-bootstrap-assumption} that
\begin{equation}\label{eq:heat-pressure-work}
 \left|\int_{\R^2}P\Div u\,\dot\theta\dd x\right|
 \le\frac{c_v}{4}\|\sqrt\rho\,\dot\theta\|_2^2
       +C\int_{\R^2}\theta^2|\nabla u|^2\dd x
\end{equation}
and
\begin{equation}\label{eq:heat-viscous-derivative}
\begin{aligned}
 &\left|\int_{\R^2}\theta(\partial_t+u\cdot\nabla+\Div u)
       \bigl(2\mu|D(u)|^2+\lambda(\Div u)^2\bigr)\dd x\right|\\
 \le& C\int_{\R^2}\theta
       \bigl(|\nabla u||\nabla\dot u|+|\nabla u|^3\bigr)\dd x\\
 \le& \zeta\|\nabla\dot u\|_2^2
       +C_\zeta\int_{\R^2}\theta^2|\nabla u|^2\dd x
       +C\|\nabla u\|_4^4 .
\end{aligned}
\end{equation}
Integration by parts and \eqref{eq:gn-lp} give
\begin{equation}\label{eq:heat-transport-product}
 \left|\int_{\R^2}\Delta\theta\,u\cdot\nabla\theta\dd x\right|
 \le C\|\nabla u\|_2\|\nabla\theta\|_4^2
 \le C\|\nabla u\|_2\|\nabla\theta\|_2\|\nabla^2\theta\|_2 .
\end{equation}
By \eqref{eq:perturbation-system}$_3$ and \eqref{eq:initial-bootstrap-assumption}, we get
\[
 \|\nabla^2\theta\|_2
 \le \frac{c_v\sqrt{2\bar\rho}}{\kappa}\|\sqrt\rho\,\dot\theta\|_2
       +C(\|\nabla u\|_4^2+\|P\|_4^2),
\]
which yields from \eqref{eq:higher-cross-identity}--\eqref{eq:heat-transport-product} that, for each \(\zeta>0\),
\begin{equation}\label{eq:heat-source-products}
\begin{aligned}
 &\frac{\dd}{\dd t}\left[\frac{\kappa}{2}\|\nabla\theta\|_2^2
      -\int_{\R^2}\theta\bigl(2\mu|D(u)|^2+\lambda(\Div u)^2\bigr)\dd x\right]
       +\frac{c_v}{2}\|\sqrt\rho\,\dot\theta\|_2^2\\
 \le& \zeta\|\nabla\dot u\|_2^2
       +C_\zeta\left(\|\nabla\theta\|_2^2+\|\nabla u\|_4^4+\|P\|_4^4
                     +\int_{\R^2}\theta^2|\nabla u|^2\dd x\right).
\end{aligned}
\end{equation}

Using H\"older's inequality and \eqref{eq:density-anchor} with \(p=8\),
we deduce from \eqref{eq:physical-first-moment} that
\[
 \|P\|_4^4\le (2\bar\rho)^3
   \|\sqrt\rho\,\theta\|_2^{4/3}\|\rho^{1/8}\theta\|_8^{8/3}
 \le C\bigl(1+\|\nabla\theta\|_2^{8/3}\bigr).
\]
It follows from \eqref{eq:div-curl-estimate}, \eqref{eq:flux-interpolation} and
\eqref{eq:initial-rough-lemma} that
\begin{equation}\label{eq:initial-quartic-products}
 \|\nabla u\|_4^4+\|P\|_4^4
 \le C\bigl(1+\|\sqrt\rho\,\dot u\|_2^2+\|\nabla\theta\|_2^{8/3}\bigr).
\end{equation}
Taking \((\eta,p)=(1/128,512)\) in \eqref{eq:spatial-weighted-lp}
and using \eqref{eq:initial-heating-interpolation-second}, we obtain
\begin{equation}\label{eq:initial-squared-heating}
\begin{aligned}
 \int_{\R^2}\theta^2|\nabla u|^2\dd x
 \le&\|(1+|x|^2)^{-1/256}\theta\|_{512}^2
       \|(1+|x|^2)^{1/256}\nabla u\|_{512/255}^2\\
 \le&C\bigl(1+\|\sqrt\rho\,\dot u\|_2^2+\|\nabla\theta\|_2^2\bigr)^{145/128}.
\end{aligned}
\end{equation}

In \eqref{eq:heat-source-products}, choose
\begin{equation}\label{eq:short-thermal-absorption-choice}
 \zeta=\frac{\mu}{2\bigl(1+32\bar\rho/(\mu c_v)\bigr)}.
\end{equation}
It follows from \eqref{eq:initial-sublinear-heating} and Young's inequality that
\[
 0\le\int_{\R^2}\theta\bigl(2\mu|D(u)|^2+\lambda(\Div u)^2\bigr)\dd x
 \le \frac{\|\sqrt\rho\,\dot u\|_2^2}{2(1+32\bar\rho/(\mu c_v))}
       +\frac{\kappa}{4}\|\nabla\theta\|_2^2+C .
\]
Multiplying \eqref{eq:heat-source-products} by \(s(1+32\bar\rho/(\mu c_v))\)
and \eqref{eq:material-acceleration-products} by \(s\), adding and integrating over \([0,t]\), we obtain from \eqref{eq:initial-quartic-products} and
\eqref{eq:initial-squared-heating} that
\[
\begin{aligned}
 &t\bigl(\|\sqrt\rho\,\dot u(t)\|_2^2+\|\nabla\theta(t)\|_2^2\bigr)
    +\int_0^t s\bigl(\|\nabla\dot u\|_2^2+\|\sqrt\rho\,\dot\theta\|_2^2\bigr)\dd s\\
 \le& C\int_0^t \bigl(1+\|\sqrt\rho\,\dot u\|_2^2+\|\nabla\theta\|_2^2\bigr)\dd s+C\int_0^t s\bigl(1+\|\sqrt\rho\,\dot u\|_2^2+\|\nabla\theta\|_2^2\bigr)^{4/3}\dd s .
\end{aligned}
\]
By \eqref{eq:initial-rough-lemma}--\eqref{eq:initial-rough-lemma-second} and H\"older's inequality, we have
\[
 \int_0^t\bigl(1+\|\sqrt\rho\,\dot u\|_2^2+\|\nabla\theta\|_2^2\bigr)^{1/3}\dd s
 \le t^{2/3}\left(\int_0^t
       \bigl(1+\|\sqrt\rho\,\dot u\|_2^2+\|\nabla\theta\|_2^2\bigr)\dd s\right)^{1/3}\le C .
\]
Using Gronwall's inequality, we obtain
\begin{equation}\label{eq:initial-weighted-material-derivatives}
\begin{aligned}
 &\sup_{0\le s\le t}s\bigl(\|\sqrt\rho\,\dot u(s)\|_2^2+\|\nabla\theta(s)\|_2^2\bigr)
       +\int_0^t s\bigl(\|\nabla\dot u\|_2^2+\|\sqrt\rho\,\dot\theta\|_2^2\bigr)\dd s\\
 \le& C\int_0^t\bigl(1+\|\sqrt\rho\,\dot u\|_2^2+\|\nabla\theta\|_2^2\bigr)\dd s\le C .
\end{aligned}
\end{equation}

By \eqref{eq:initial-sublinear-heating} and \eqref{eq:initial-weighted-material-derivatives},
we have
\[
 t\int_{\R^2}\theta\bigl(2\mu|D(u)|^2+\lambda(\Div u)^2\bigr)\dd x
 \le C t^{47/128},
\]
and
\[
 \int_0^t s\bigl(1+\|\sqrt\rho\,\dot u\|_2^2+\|\nabla\theta\|_2^2\bigr)^{4/3}\dd s
 \le C t^{2/3}\int_0^t\bigl(1+\|\sqrt\rho\,\dot u\|_2^2+\|\nabla\theta\|_2^2\bigr)\dd s
 \le C t^{2/3}.
\]
Multiplying \eqref{eq:heat-source-products} by \(s\) and integrating over \([0,t]\), we obtain,
for each \(\zeta>0\),
\begin{equation}\label{eq:initial-weighted-gradient-comparison}
 \frac{\kappa}{2}t\|\nabla\theta(t)\|_2^2
 \le \frac{\kappa}{2}\int_0^t\|\nabla\theta(s)\|_2^2\dd s
       +C\zeta+Ct^{47/128}+C_\zeta t^{2/3}.
\end{equation}
For \(C\) in \eqref{eq:initial-weighted-gradient-comparison}, we take
\(\zeta=\min\{1,\kappa m_0^2/(4C)\}\) and choose \(\tau\) so that
\[
 C\tau^{47/128}+C_\zeta\tau^{2/3}\le\kappa m_0^2/4 .
\]
Then we have
\[
 t\|\nabla\theta(t)\|_2^2\le\int_0^t\|\nabla\theta(s)\|_2^2\dd s+m_0^2 .
\]
Combining this inequality with \eqref{eq:initial-rough-lemma} gives
\eqref{eq:initial-weighted-temperature-gradient}.
\end{proof}

Using \eqref{eq:initial-weighted-lemma-second}, we now establish the following
time-weighted estimates for $\sqrt\rho\,\dot u$ and $\nabla\theta$ on
$[0,\min\{T,\tau\}]$.
\begin{lemma}\label{lem:short-squared-weight}
Under the assumptions of Proposition~\ref{prop:lower-order-closure},
 we have
\begin{equation}\label{eq:short-squared-weight}
\begin{aligned}
 &\sup_{0\le t\le \min\{T,\tau\}}\sigma^2(t)
   \bigl(\|\sqrt\rho\,\dot u(t)\|_2^2+\|\nabla\theta(t)\|_2^2\bigr)+\int_0^{\min\{T,\tau\}}\sigma^2(t)
   \bigl(\|\nabla\dot u\|_2^2+\|\sqrt\rho\,\dot\theta\|_2^2\bigr)\dd t
 \le CE_0^{1/12},
\end{aligned}
\end{equation}
where \(\tau\) as in Lemma~\ref{lem:initial-time-weighted}.
\end{lemma}
\begin{proof}
Let \(0\le t\le\min\{T,\tau\}\).
Taking \(p=4\) in \eqref{eq:density-anchor},
\eqref{eq:div-curl-estimate} and \eqref{eq:flux-interpolation}, we obtain
from \eqref{eq:initial-bootstrap-assumption}, \eqref{eq:physical-first-moment}
and \eqref{eq:initial-weighted-gradient} that
\begin{equation}\label{eq:short-spatial-norms}
\begin{aligned}
 &\|\rho^{1/4}\theta\|_4+\|P\|_4+\|\nabla u\|_4+\|(1+|x|^2)^{1/32}\nabla u\|_2\\
 \le& C\bigl(\|\nabla u\|_2+\|\sqrt\rho\,\theta\|_2+\|\sqrt\rho\,\dot u\|_2+\|\nabla\theta\|_2\bigr).
\end{aligned}
\end{equation}
Using the first inequality in \eqref{eq:initial-heating-interpolation-second}, we get
\begin{equation}\label{eq:short-gradient-product}
\begin{aligned}
 &\|(1+|x|^2)^{1/256}\nabla u\|_{512/255}^2
 \le C\|\nabla u\|_2^{3/2}\bigl(\|\nabla u\|_2+\|\sqrt\rho\,\theta\|_2+\|\sqrt\rho\,\dot u\|_2+\|\nabla\theta\|_2\bigr)^{1/2}.
\end{aligned}
\end{equation}
For \(z\in H^1_{\rm loc}\) with \(\sqrt\rho\,z,\nabla z\in L^2\),
H\"older's inequality, \eqref{eq:density-anchor} and \eqref{eq:spatial-weighted-lp} give
\begin{equation}\label{eq:short-spatial-multiplier}
\begin{aligned}
 &\int_{\R^2}(|\nabla u|^2+P^2)|z|^2\dd x\\
 \le&\|(1+|x|^2)^{1/256}\nabla u\|_{512/255}^2
       \|(1+|x|^2)^{-1/256}z\|_{512}^2
       +2\bar\rho\|\rho^{1/4}\theta\|_4^2\|\rho^{1/4}z\|_4^2\\
 \le& C\bigl(\|\nabla u\|_2^2+\|\sqrt\rho\,\theta\|_2^2+\|\sqrt\rho\,\dot u\|_2^2+\|\nabla\theta\|_2^2\bigr)
          \bigl(\|\sqrt\rho\,z\|_2+\|\nabla z\|_2\bigr)^2,
\end{aligned}
\end{equation}
where the last inequality uses \eqref{eq:short-spatial-norms} and \eqref{eq:short-gradient-product}.
Combining the first inequality in \eqref{eq:initial-sublinear-heating} with
\eqref{eq:spatial-weighted-lp} and \eqref{eq:short-gradient-product}, we deduce from
Young's inequality and \eqref{eq:initial-bootstrap-assumption} that
\begin{equation}\label{eq:short-viscous-bound}
\begin{aligned}
 0\le&\int_{\R^2}\theta\bigl(2\mu|D(u)|^2+\lambda(\Div u)^2\bigr)\dd x\\
 \le& C\|\nabla u\|_2^{3/2}\bigl(\|\nabla u\|_2+\|\sqrt\rho\,\theta\|_2+\|\sqrt\rho\,\dot u\|_2+\|\nabla\theta\|_2\bigr)^{3/2}\\
 \le&\frac{\mu c_v}{2(\mu c_v+32\bar\rho)}\|\sqrt\rho\,\dot u\|_2^2
     +\frac{\kappa}{4}\|\nabla\theta\|_2^2
     +C\bigl(\|\nabla u\|_2^2+\|\sqrt\rho\,\theta\|_2^2\bigr).
\end{aligned}
\end{equation}
By \eqref{eq:initial-bootstrap-assumption} and \eqref{eq:initial-weighted-material-derivatives}, we have
\[
 s\bigl(\|\nabla u\|_2^2+\|\sqrt\rho\,\theta\|_2^2+\|\sqrt\rho\,\dot u\|_2^2
            +\|\nabla\theta\|_2^2\bigr)\le C,\qquad 0\le s\le t,
\]
which, together with \eqref{eq:short-spatial-norms} and
\eqref{eq:short-spatial-multiplier} with \(z=\theta\), gives
\begin{equation}\label{eq:short-weighted-products}
\begin{aligned}
 &s^2\left(\|\nabla u\|_4^4+\|P\|_4^4+\int_{\R^2}\theta^2|\nabla u|^2\dd x\right)\\
 \le& Cs\bigl(\|\nabla u\|_2^2+\|\sqrt\rho\,\theta\|_2^2+\|\sqrt\rho\,\dot u\|_2^2+\|\nabla\theta\|_2^2\bigr).
\end{aligned}
\end{equation}
With \(\zeta\) given by \eqref{eq:short-thermal-absorption-choice}, multiplying \eqref{eq:heat-source-products} by \(s^2(1+32\bar\rho/(\mu c_v))\)
and \eqref{eq:material-acceleration-products} by \(s^2\), adding and integrating over \([0,t]\), we obtain from \eqref{eq:short-viscous-bound}
and \eqref{eq:short-weighted-products} that
\[
\begin{aligned}
 &t^2\bigl(\|\sqrt\rho\,\dot u(t)\|_2^2+\|\nabla\theta(t)\|_2^2\bigr)
      +\int_0^t s^2\bigl(\|\nabla\dot u\|_2^2+\|\sqrt\rho\,\dot\theta\|_2^2\bigr)\dd s\\
 \le& Ct^2\bigl(\|\nabla u(t)\|_2^2+\|\sqrt\rho\,\theta(t)\|_2^2\bigr)\\
 &+C\int_0^t s\bigl(\|\nabla u\|_2^2+\|\sqrt\rho\,\theta\|_2^2+\|\sqrt\rho\,\dot u\|_2^2+\|\nabla\theta\|_2^2\bigr)\dd s\\
 \le& CE_0^{1/12},
\end{aligned}
\]
where the last inequality follows from \eqref{eq:initial-weighted-lemma-second} and \(t\le1\).
Taking the supremum and using \(\sigma(t)=t\) on \([0,\min\{T,\tau\}]\) gives
\eqref{eq:short-squared-weight}.
\end{proof}

To control the density along particle paths for $t\ge\tau$, we first
estimate the effective viscous flux in $L^\infty$ in terms of the material derivative of the velocity.
\begin{lemma}\label{lem:acceleration-and-flux}
Under the assumptions of Proposition~\ref{prop:lower-order-closure}, if \(0\le t\le T\)
and \(\nabla\dot u(t)\in L^2\), then
\begin{equation}\label{eq:lift-flux-supremum}
 \|F\|_\infty
 \le C\log^{1/6}(e+t)\|\sqrt\rho\,\dot u\|_2^{2/3}\|\nabla\dot u\|_2^{1/3},
\end{equation}
where \(C=C(m_0,\bar\rho,M,\mu,\lambda,c_v,\kappa)\).
\end{lemma}
\begin{proof}
Substituting \eqref{eq:lift-lame-gradient-representation} into \eqref{eq:flux-definitions}, we  obtain
\[
 F=(2\mu+\lambda)\sum_{i,j=1}^2(\partial_i\Gamma_{ij})*(\rho\dot u^j).
\]
Using \eqref{eq:riesz-potential-lp} with \(p=4/3\), H\"older's inequality, \eqref{eq:mass} and \eqref{eq:initial-bootstrap-assumption}, we have
\[
 \|F\|_4
 \le C\|I_1(|\rho\dot u|)\|_4
 \le C\|\rho\dot u\|_{4/3}
 \le C\|\sqrt\rho\|_4\|\sqrt\rho\,\dot u\|_2
 \le C\|\sqrt\rho\,\dot u\|_2.
\]
Taking \(w=\dot u\) and \(p=8\) in \eqref{eq:lift-log-anchor}, we deduce from \eqref{eq:mean-material-acceleration} and \eqref{eq:initial-bootstrap-assumption} that
\[
 \|\rho\dot u\|_8
 \le (2\bar\rho)^{7/8}\|\rho^{1/8}(\dot u-\langle\dot u\rangle_\rho)\|_8
 \le C\sqrt{\log(e+t)}\,\|\nabla\dot u\|_2.
\]
By \eqref{eq:gn-infty} with \(q=r=4\), \eqref{eq:flux-gradient-lp} with \(p=4\)
and \eqref{eq:initial-bootstrap-assumption}, we have
\[
\begin{aligned}
 \|F\|_\infty
 \le& C\|F\|_4^{1/2}\|\nabla F\|_4^{1/2}
 \le C\|F\|_4^{1/2}\|\rho\dot u\|_4^{1/2}\\
 \le& C\|F\|_4^{1/2}\|\rho\dot u\|_2^{1/6}\|\rho\dot u\|_8^{1/3}\\
 \le& C\|\sqrt\rho\,\dot u\|_2^{2/3}\|\rho\dot u\|_8^{1/3}\\
 \le& C\log^{1/6}(e+t)\|\sqrt\rho\,\dot u\|_2^{2/3}\|\nabla\dot u\|_2^{1/3},
\end{aligned}
\]
which gives \eqref{eq:lift-flux-supremum}.
\end{proof}

To handle the terms $2\mu|D(u)|^2+\lambda(\Div u)^2$ and $P\Div u$ in
\eqref{eq:perturbation-system}$_3$, we establish the following estimates for
$|\nabla u|^2+P^2$.
\begin{lemma}\label{lem:velocity-pressure-multiplier}
Under the assumptions of Proposition~\ref{prop:lower-order-closure}, for \(0\le t\le T\), we have
\begin{equation}\label{eq:lift-source-moment}
\begin{aligned}
 &\bigl\||\nabla u|^2+P^2\bigr\|_{3/2}\le C(\|\sqrt\rho\dot u\|_2^2+\|P\|_3^2),
\end{aligned}
\end{equation}
and
\begin{equation}\label{eq:lift-source-moment-second}
\begin{aligned}
 &\int_{\R^2}(1+|x|^2)^{1/8}(|\nabla u|^2+P^2)\dd x\\
 \le&C(1+t)^{1/2}\|\sqrt\rho\,\dot u\|_2^2
  +C(1+t)^{1/4}\bigl(\langle\theta\rangle_\rho^2+\log(e+t)\|\nabla\theta\|_2^2\bigr).
\end{aligned}
\end{equation}
If \(T\ge\tau\), then we also have, for \(\tau\le t\le T\) and \(z\in\mathcal H_\rho\),
\begin{equation}\label{eq:lift-source-multiplier}
\begin{aligned}
 \int_{\R^2} (|\nabla u|^2+P^2)z^2\dd x
 \le&C\bigg[\log(e+t)(\|\nabla u\|_2^2+\|P\|_2^2)^{1/2}\\
 &+(\|\nabla u\|_2^2+\|P\|_2^2)
         (\|\sqrt\rho\,\dot u\|_2^2+\|P\|_3^3)^{1/4}\bigg]\|\nabla z\|_2^2.
\end{aligned}
\end{equation}
\end{lemma}
\begin{proof}
By \eqref{eq:flux-definitions}, \eqref{eq:lift-lame-gradient-representation}
and \eqref{eq:div-curl-estimate} with \(p=3\), we deduce from
\eqref{eq:riesz-potential-lp} with \(p=6/5\) and H\"older's inequality that
\[
\begin{aligned}
 \|\nabla u\|_3
 &\le C\bigl(\|I_1(|\rho\dot u|)\|_3+\|P\|_3\bigr)\\
 &\le C\bigl(\|\rho\dot u\|_{6/5}+\|P\|_3\bigr)\\
 &\le C\bigl(\|\sqrt\rho\|_3\|\sqrt\rho\,\dot u\|_2+\|P\|_3\bigr),
\end{aligned}
\]
which gives \eqref{eq:lift-source-moment}.
For \(0<\alpha\le1\), H\"older's inequality and \eqref{eq:physical-first-moment} give
\begin{equation}\label{eq:lift-first-moment-second}
 \int_{\R^2}(1+|x|^2)^{\alpha/2}\rho(t)\dd x
 \le m_0^{1-\alpha}\bigl(m_0+M+\sqrt{2m_0E_0}\,t\bigr)^\alpha
 \le C(1+t)^\alpha.
\end{equation}
Taking \(\alpha=1/4\) and \(\alpha=3/8\) in \eqref{eq:lift-first-moment-second} and using
\eqref{eq:initial-bootstrap-assumption}, we obtain
\begin{equation}\label{eq:density-pressure-weight-moments}
 0<\|\rho^2(1+|x|^2)^{1/8}\|_1
 \le 2\bar\rho\int_{\R^2}(1+|x|^2)^{1/8}\rho\dd x
 \le C(1+t)^{1/4},
\end{equation}
and
\begin{equation}\label{eq:density-pressure-weight-moments-second}
 \|\rho^2(1+|x|^2)^{1/8}\|_{3/2}
 \le (2\bar\rho)^{4/3}
       \left(\int_{\R^2}(1+|x|^2)^{3/16}\rho\dd x\right)^{2/3}
 \le C(1+t)^{1/4}.
\end{equation}
Similarly, taking \(\alpha=1/2\) in \eqref{eq:lift-first-moment-second} and using \eqref{eq:initial-bootstrap-assumption} gives
\begin{equation}\label{eq:density-pressure-weight-moments-third}
 \int_{\R^2}\rho^2(1+|x|^2)^{1/4}\dd x
 \le 2\bar\rho\int_{\R^2}(1+|x|^2)^{1/4}\rho\dd x
 \le C(1+t)^{1/2}.
\end{equation}
Taking \(v=\theta\) in \eqref{eq:lift-two-measure} with \(w=\rho^2(1+|x|^2)^{1/8}\),
we deduce from \eqref{eq:lift-first-moment-second},
\eqref{eq:density-pressure-weight-moments}--\eqref{eq:density-pressure-weight-moments-third} and the boundedness of \(s\log(e+1/s)\) on \((0,1]\) that
\begin{equation}\label{eq:pressure-weighted-moment}
 \int_{\R^2}(1+|x|^2)^{1/8}P^2\dd x
 \le C(1+t)^{1/4}\bigl(\langle\theta\rangle_\rho^2
                                  +\log(e+t)\|\nabla\theta\|_2^2\bigr).
\end{equation}

For \(|x-y|\le2\), we have \((1+|x|^2)^{1/16}\le C(1+|y|^2)^{1/16}\).
By \eqref{eq:lift-lame-fundamental-solution}, Young's inequality,
H\"older's inequality and \eqref{eq:density-pressure-weight-moments-third}, we obtain
\begin{equation}\label{eq:force-near-weighted}
\begin{aligned}
 &\left\|(1+|x|^2)^{1/16}
       \int_{|x-y|\le2}(\partial_i\Gamma)(x-y)(\rho\dot u)(y)\dd y\right\|_2\\
 \le& C\||x|^{-1}\mathbf1_{B_2}\|_{4/3}\|(1+|x|^2)^{1/16}\rho\dot u\|_{4/3}\\
 \le& C\|(1+|x|^2)^{1/16}\sqrt\rho\|_4
                    \|\sqrt\rho\dot u\|_2
 \le C(1+t)^{1/8}\|\sqrt\rho\dot u\|_2.
\end{aligned}
\end{equation}
Since
\(|(\partial_i\Gamma)(x-y)-(\partial_i\Gamma)(x)|\le C|y|/|x|^2\) for \(|x|>2|y|+4\), we have
\[
\begin{aligned}
 &\int_{\R^2}(1+|x|^2)^{1/8}
   \left|\mathbf1_{\{|x-y|>2\}}(\partial_i\Gamma)(x-y)
          -\mathbf1_{\{|x|>2\}}(\partial_i\Gamma)(x)\right|^2\dd x\\
 \le& C(1+|y|^2)^{1/8}\log(e+|y|)
 \le C(1+|y|^2)^{1/4}.
\end{aligned}
\]
By Minkowski's inequality, \eqref{eq:mean-material-acceleration} and \eqref{eq:lift-first-moment-second}, we obtain
\begin{equation}\label{eq:force-far-weighted}
\begin{aligned}
 &\left\|(1+|x|^2)^{1/16}
       \int_{|x-y|>2}(\partial_i\Gamma)(x-y)(\rho\dot u)(y)\dd y\right\|_2\\
 \le& C\int_{\R^2}(1+|y|^2)^{1/8}\rho|\dot u|\dd y
 \le C(1+t)^{1/4}\|\sqrt\rho\dot u\|_2.
\end{aligned}
\end{equation}
Substituting \eqref{eq:pressure-weighted-moment} and \eqref{eq:force-near-weighted}--\eqref{eq:force-far-weighted}
into \eqref{eq:lift-lame-gradient-representation}, and applying \eqref{eq:weighted-riesz-pressure}
with \(\beta=1/8\), we obtain \eqref{eq:lift-source-moment-second}.

Let \(\tau\le t\le T\) and \(w=|\nabla u|^2+P^2\). By \eqref{eq:initial-bootstrap-assumption}, we have
\begin{equation}\label{eq:lift-source-weight-admissibility}
 \|w\|_1=\|\nabla u\|_2^2+\|P\|_2^2\le 2\max\{1,2\bar\rho\}K.
\end{equation}
For \(w\not\equiv0\), \eqref{eq:lift-two-measure}, \eqref{eq:physical-total-energy},
\eqref{eq:initial-bootstrap-assumption}, \eqref{eq:lift-source-moment}--\eqref{eq:lift-source-moment-second}
and \eqref{eq:lift-first-moment-second} with \(\alpha=1/4\) give
\[
\begin{aligned}
 \int_{\R^2}wz^2\dd x
 \le& C\|w\|_1\bigl[\log(e+t)+\log(e+\|w\|_1^{-1})
       +\log(e+\|\sqrt\rho\,\dot u\|_2^2+\|P\|_3^3)\bigr]\|\nabla z\|_2^2\\
 \le& C\bigl[\log(e+t)\|w\|_1^{1/2}
       +\|w\|_1(\|\sqrt\rho\,\dot u\|_2^2+\|P\|_3^3)^{1/4}\bigr]\|\nabla z\|_2^2,
\end{aligned}
\]
where the last inequality uses \eqref{eq:lift-source-weight-admissibility}
and \(\log(e+s)\le C(1+s^{1/4})\) for \(s\ge0\).
This proves \eqref{eq:lift-source-multiplier}.
\end{proof}

The following lemma gives the energy estimates for $\sqrt\rho\,\dot u$, $P$,
$F$ and $\omega$.
\begin{lemma}\label{lem:mechanical-and-pressure-energy}
Under the assumptions of Proposition~\ref{prop:lower-order-closure}, for \(0\le t\le T\), we have
\begin{equation}\label{eq:lift-mechanical-interpolation}
 \|\nabla u\|_2^2+\|P\|_2^2
 \le CE_0^{1/2}\bigl(\|\sqrt\rho\,\dot u\|_2^2+\|P\|_3^3\bigr)^{1/2}.
\end{equation}
There exists \(C_0=C_0(m_0,\bar\rho,\mu,\lambda)\ge\max\{1,8\bar\rho/\mu\}\) such that
for \(0\le t\le T\),
\begin{equation}\label{eq:lift-acceleration-differential}
\begin{aligned}
 &\|\sqrt{\rho(t)}\,\dot u(t)\|_2^2+\mu\int_0^t\|\nabla\dot u\|_2^2\dd r\\
 \le&\int_{\{\rho_0>0\}}|g|^2\dd x
      +C_0\int_0^t\bigl(\|\sqrt\rho\,\dot u\|_2^4+\|P\|_4^4+
                         \|\sqrt\rho\,\dot\theta\|_2^2\bigr)\dd r,
\end{aligned}
\end{equation}
and
\begin{equation}\label{eq:lift-pressure-differential}
\begin{aligned}
 &\|P(t)\|_3^3+\frac{1}{2\mu+\lambda}\int_0^t\|P\|_4^4\dd r
 \le\|P_0\|_3^3
      +C_0\int_0^t\bigl(\|\sqrt\rho\,\dot u\|_2^4+\|\sqrt\rho\,\dot\theta\|_2^2\bigr)\dd r.
\end{aligned}
\end{equation}
Furthermore, we also have
\begin{equation}\label{eq:lift-short-flux-energy}
\begin{aligned}
 &\frac{\|F(t)\|_2^2}{2(2\mu+\lambda)}+\frac{\mu}{2}\|\omega(t)\|_2^2
     +\frac14\int_0^t\|\sqrt\rho\,\dot u\|_2^2\dd r\\
 \le&\frac{\|F_0\|_2^2}{2(2\mu+\lambda)}+\frac{\mu}{2}\|\omega_0\|_2^2
     +C\int_0^t\bigl(\|\nabla\theta\|_2^2+E_0\|P\|_3^3\bigr)\dd r,
\end{aligned}
\end{equation}
where \(C=C(\bar\rho,\mu,\lambda,c_v,\kappa)\).
\end{lemma}
\begin{proof}
Multiplying \eqref{eq:perturbation-system}$_2$ by \(u\) and integrating by parts, we obtain
\begin{equation}\label{eq:mechanical-energy-testing}
 \mu\|\nabla u\|_2^2+(\mu+\lambda)\|\Div u\|_2^2
 =\int_{\R^2}P\Div u\dd x-\int_{\R^2}\rho u\cdot\dot u\dd x.
\end{equation}
By H\"older's inequality and \eqref{eq:physical-total-energy}, we have
\[
 \|P\|_2^2\le\|P\|_1^{1/2}\|P\|_3^{3/2}
 \le (E_0/c_v)^{1/2}\|P\|_3^{3/2},
 \qquad \|\sqrt\rho\,u\|_2\le(2E_0)^{1/2}.
\]
Applying Young's inequality in \eqref{eq:mechanical-energy-testing} and using the preceding estimates, we obtain
\[
 \|\nabla u\|_2^2+\|P\|_2^2
 \le C\bigl(\|P\|_2^2+\|\sqrt\rho\,u\|_2\|\sqrt\rho\,\dot u\|_2\bigr)
 \le CE_0^{1/2}\bigl(\|\sqrt\rho\,\dot u\|_2+\|P\|_3^{3/2}\bigr),
\]
which gives \eqref{eq:lift-mechanical-interpolation}.

Combining \eqref{eq:flux-definitions} and \eqref{eq:lift-lame-gradient-representation}
with \eqref{eq:div-curl-estimate} at \(p=4\), we use H\"older's inequality, \eqref{eq:mass},
\eqref{eq:initial-bootstrap-assumption} and \eqref{eq:riesz-potential-lp} with \(p=4/3\) to obtain
\begin{equation}\label{eq:quartic-velocity-pressure}
\begin{aligned}
 \|\nabla u\|_4^4+\|P\nabla u\|_2^2
 \le& C\bigl(\|I_1(|\rho\dot u|)\|_4^4+\|P\|_4^4\bigr)\\
 \le& C\bigl(\|\sqrt\rho\,\dot u\|_2^4+\|P\|_4^4\bigr).
\end{aligned}
\end{equation}
Integrating \eqref{eq:material-acceleration-products} over \([0,t]\) and using
\eqref{eq:compatibility}, \eqref{eq:perturbation-system}$_2$ and \eqref{eq:quartic-velocity-pressure},
we obtain \eqref{eq:lift-acceleration-differential}.

Next, using \eqref{eq:perturbation-system}$_1$ to write \(\dot P+P\Div u=\rho\dot\theta\),
multiplying this equation by \(3P^2\) and integrating over \(\R^2\), we obtain from \eqref{eq:flux-definitions} that
\[
 \frac{\mathrm d}{\mathrm dt}\|P\|_3^3+\frac{2}{2\mu+\lambda}\|P\|_4^4
 =-\frac{2}{2\mu+\lambda}\int_{\R^2}P^3F\dd x
   +3\int_{\R^2}\rho P^2\dot\theta\dd x.
\]
By \eqref{eq:flux-definitions}, \eqref{eq:lift-lame-gradient-representation},
\eqref{eq:riesz-potential-lp} with \(p=4/3\), \eqref{eq:mass},
\eqref{eq:initial-bootstrap-assumption} and H\"older's and Young's inequalities, we deduce
\[
\begin{aligned}
 \frac{2}{2\mu+\lambda}\left|\int_{\R^2}P^3F\dd x\right|
      +3\left|\int_{\R^2}\rho P^2\dot\theta\dd x\right|
 \le& C\|P\|_4^3\|I_1(|\rho\dot u|)\|_4
      +C\|P\|_4^2\|\sqrt\rho\,\dot\theta\|_2\\
 \le& C\|P\|_4^3\|\sqrt\rho\,\dot u\|_2
      +C\|P\|_4^2\|\sqrt\rho\,\dot\theta\|_2\\
 \le& \frac{1}{2\mu+\lambda}\|P\|_4^4
      +C\bigl(\|\sqrt\rho\,\dot u\|_2^4+\|\sqrt\rho\,\dot\theta\|_2^2\bigr),
\end{aligned}
\]
which gives
\begin{equation}\label{eq:cubic-pressure-differential}
 \frac{\dd}{\dd t}\|P\|_3^3+\frac{1}{2\mu+\lambda}\|P\|_4^4
 \le C_0\bigl(\|\sqrt\rho\,\dot u\|_2^4+\|\sqrt\rho\,\dot\theta\|_2^2\bigr).
\end{equation}
Integrating \eqref{eq:cubic-pressure-differential} over \([0,t]\) gives \eqref{eq:lift-pressure-differential}.

Finally, substituting
\[
 \partial_i u^j\partial_j u^i=(\Div u)^2-2\det\nabla u
\]
and
\[
 2\mu|D(u)|^2+\lambda(\Div u)^2
 =(2\mu+\lambda)(\Div u)^2+\mu\omega^2-4\mu\det\nabla u
\]
into \eqref{eq:initial-flux-identity}, we arrive at
\begin{equation}\label{eq:regrouped-planar-flux}
\begin{aligned}
 &\frac{\dd}{\dd t}\left(\frac{\|F\|_2^2}{2(2\mu+\lambda)}
                         +\frac\mu2\|\omega\|_2^2\right)
                         +\|\sqrt\rho\,\dot u\|_2^2\\
 =&-\frac{c_v^{-1}+1/2}{2\mu+\lambda}\int_{\R^2}(\Div u)F^2\dd x
       -\mu(c_v^{-1}+1/2)\int_{\R^2}(\Div u)\omega^2\dd x
       +\frac{\mu}{c_v(2\mu+\lambda)}\int_{\R^2} P\omega^2\dd x\\
 &+2\left(1+\frac{2\mu}{c_v(2\mu+\lambda)}\right)
                         \int_{\R^2} F\det\nabla u\dd x
       +\frac{\kappa}{c_v(2\mu+\lambda)}\int_{\R^2}\nabla\theta\cdot\nabla F\dd x.
\end{aligned}
\end{equation}
By \eqref{eq:gn-lp} and \eqref{eq:flux-gradient-lp}, we have
\[
 \|F\|_4^2\le C\|F\|_2\|\sqrt\rho\,\dot u\|_2,
 \qquad
 \|\omega\|_4^2\le C\|\nabla u\|_2\|\sqrt\rho\,\dot u\|_2,
\]
and by \eqref{eq:lift-planar-divcurl-bound} with \(f=u^1\) and \(g=u^2\), we obtain
\[
 \left|\int_{\R^2} F\det\nabla u\dd x\right|
 \le\frac{\sqrt3}{8\sqrt\pi}\|\nabla F\|_2\|\nabla u\|_2^2.
\]
Applying Young's inequality in \eqref{eq:regrouped-planar-flux} and using the preceding estimates gives
\begin{equation}\label{eq:regrouped-flux-differential}
\frac{\dd}{\dd t}\left(\frac{\|F\|_2^2}{2(2\mu+\lambda)}+\frac\mu2\|\omega\|_2^2\right)+\tfrac12\|\sqrt\rho\,\dot u\|_2^2
 \le C\bigl[(\|\nabla u\|_2^2+\|P\|_2^2)^2+\|\nabla\theta\|_2^2\bigr].
\end{equation}
Substituting \eqref{eq:lift-mechanical-interpolation} into \eqref{eq:regrouped-flux-differential}, we obtain
\[
 \frac{\dd}{\dd t}\left(\frac{\|F\|_2^2}{2(2\mu+\lambda)}+\frac\mu2\|\omega\|_2^2\right)
 +\tfrac12\|\sqrt\rho\,\dot u\|_2^2
 \le CE_0\bigl(\|\sqrt\rho\,\dot u\|_2^2+\|P\|_3^3\bigr)+C\|\nabla\theta\|_2^2.
\]
Taking \(C\eps_0\le1/4\) and integrating over \([0,t]\) gives \eqref{eq:lift-short-flux-energy}.
\end{proof}

To eliminate the terms $2\mu|D(u)|^2+\lambda(\Div u)^2-P\Div u$ from the
equation for $\theta-\phi$, we define $\phi$ by the following elliptic
problem.
\begin{lemma}\label{lem:elliptic-temperature-correction}
Under the assumptions of Proposition~\ref{prop:lower-order-closure}, for each \(0\le t\le T\), the problem
\begin{equation}\label{eq:lift-elliptic-source}
\left\{
\begin{aligned}
 &-\kappa\Delta\phi=2\mu|D(u)|^2+\lambda(\Div u)^2-P\Div u-\frac{\rho}{m_0}
       \int_{\R^2}\bigl[2\mu|D(u)|^2+\lambda(\Div u)^2-P\Div u\bigr]\dd y,\\
 &\int_{\R^2}\rho\phi\dd x=0.
\end{aligned}
\right.
\end{equation}
has a unique weak solution \(\phi\in\mathcal H_\rho\). Moreover, we have 
\begin{equation}\label{eq:lift-gradient-testing}
 \kappa\|\nabla\phi\|_2^2
 \le C\|\nabla u\|_2
       \left(\int_{\R^2}(|\nabla u|^2+P^2)\phi^2\dd x\right)^{1/2},
\end{equation}
\begin{equation}\label{eq:lift-hessian}
 \|\nabla^2\phi\|_2^2
 \le C\bigl(\|\sqrt\rho\,\dot u\|_2^4+\|P\|_4^4+E_0\|\sqrt\rho\,\dot u\|_2^2\bigr),
\end{equation}
and
\begin{equation}\label{eq:lift-gradient-product}
 \||\nabla u|\nabla\phi\|_2^2
 \le C\|\nabla\phi\|_2
       \bigl(\|\sqrt\rho\,\dot u\|_2^4+\|P\|_4^4+E_0\|\sqrt\rho\,\dot u\|_2^2\bigr).
\end{equation}
\end{lemma}
\begin{proof}
For \(z_0\in\mathcal H_\rho\), H\"older's inequality gives
\begin{equation}\label{eq:lift-elliptic-testing}
 \left|\int_{\R^2}\bigl(2\mu|D(u)|^2+\lambda(\Div u)^2-P\Div u\bigr)z_0\dd x\right|
 \le C\|\nabla u\|_2
       \left(\int_{\R^2}(|\nabla u|^2+P^2)z_0^2\dd x\right)^{1/2}.
\end{equation}
Using \eqref{eq:short-spatial-multiplier} and \eqref{eq:lift-log-anchor} with \(p=2\) for \(t\le\tau\),
and \eqref{eq:lift-source-multiplier} and \eqref{eq:lift-source-weight-admissibility} for \(t>\tau\),
we deduce from \eqref{eq:initial-bootstrap-assumption} and \eqref{eq:lift-elliptic-testing} that
\[
\begin{aligned}
 &\left|\int_{\R^2}\bigl(2\mu|D(u)|^2+\lambda(\Div u)^2-P\Div u\bigr)z_0\dd x\right|\\
 \le& C\bigl(\log^{1/2}(e+t)+\|\sqrt\rho\,\dot u\|_2+\|\nabla\theta\|_2+\|P\|_3^{3/8}\bigr)
       \|\nabla z_0\|_2.
\end{aligned}
\]
The Lax--Milgram theorem gives a unique \(\phi\in\mathcal H_\rho\) satisfying
\begin{equation}\label{eq:lift-variational}
 \kappa\int_{\R^2}\nabla\phi\cdot\nabla z_0\dd x
 =\int_{\R^2}\bigl(2\mu|D(u)|^2+\lambda(\Div u)^2-P\Div u\bigr)z_0\dd x,
 \qquad z_0\in\mathcal H_\rho .
\end{equation}
Taking \(z_0=\psi-\langle\psi\rangle_\rho\) with \(\psi\in C_c^\infty(\R^2)\) in \eqref{eq:lift-variational}
yields \eqref{eq:lift-elliptic-source}, while taking \(z_0=\phi\) and using \eqref{eq:lift-elliptic-testing}
gives \eqref{eq:lift-gradient-testing}.

For any two weak solutions \(\phi_1,\phi_2\in\mathcal H_\rho\) of \eqref{eq:lift-elliptic-source},
we have
\[
 \int_{\R^2}\nabla(\phi_1-\phi_2)\cdot\nabla\psi\dd x=0,
 \qquad \psi\in C_c^\infty(\R^2).
\]
Choosing \(w_n\in C_c^\infty(\R^2)+\R\) as in the proof of Lemma~\ref{lem:lift-log-anchor}
such that \(\nabla w_n\to\nabla(\phi_1-\phi_2)\) in \(L^2\) and passing to the limit, we obtain
\[
 \|\nabla(\phi_1-\phi_2)\|_2^2
 =\lim_{n\to\infty}\int_{\R^2}\nabla(\phi_1-\phi_2)\cdot\nabla w_n\dd x=0.
\]
Hence \(\nabla(\phi_1-\phi_2)=0\), which together with \(\langle\phi_1-\phi_2\rangle_\rho=0\)
implies that \(\phi_1=\phi_2\).

Next, by \eqref{eq:viscous-dissipation-identity}, \eqref{eq:mechanical-energy-testing}
and \eqref{eq:physical-total-energy}, we have
\[
 \left|\int_{\R^2}\bigl(2\mu|D(u)|^2+\lambda(\Div u)^2-P\Div u\bigr)\dd x\right|
 =\left|\int_{\R^2}\rho u\cdot\dot u\dd x\right|
 \le (2E_0)^{1/2}\|\sqrt\rho\,\dot u\|_2.
\]
Applying Plancherel's identity to \eqref{eq:lift-elliptic-source} and using \eqref{eq:initial-bootstrap-assumption},
\eqref{eq:quartic-velocity-pressure} and the preceding estimate, we obtain
\[
\begin{aligned}
 \|\nabla^2\phi\|_2^2
 \le& C\bigl(\|\nabla u\|_4^4+\|P\nabla u\|_2^2
        +\|\rho\|_2^2 E_0\|\sqrt\rho\,\dot u\|_2^2\bigr)\\
 \le& C\bigl(\|\sqrt\rho\,\dot u\|_2^4+\|P\|_4^4+E_0\|\sqrt\rho\,\dot u\|_2^2\bigr),
\end{aligned}
\]
which proves \eqref{eq:lift-hessian}.
Finally, using H\"older's inequality, \eqref{eq:gn-lp}, \eqref{eq:quartic-velocity-pressure} and \eqref{eq:lift-hessian}, we get
\[
\begin{aligned}
 \||\nabla u|\nabla\phi\|_2^2
 \le& \|\nabla u\|_4^2\|\nabla\phi\|_4^2
 \le C\|\nabla u\|_4^2\|\nabla\phi\|_2\|\nabla^2\phi\|_2\\
 \le& C\|\nabla\phi\|_2\bigl(\|\sqrt\rho\,\dot u\|_2^4+\|P\|_4^4+E_0\|\sqrt\rho\,\dot u\|_2^2\bigr),
\end{aligned}
\]
which yields \eqref{eq:lift-gradient-product}.
\end{proof}

We now estimate $\|\sqrt\rho\,\dot\phi\|_2$ on $[\tau,T]$, which will be used
in the energy estimate for $\nabla(\theta-\phi)$.
\begin{lemma}\label{lem:temperature-correction-rate}
Under the assumptions of Proposition~\ref{prop:lower-order-closure}, let \(T\ge\tau\) and
let \(\phi\) be the solution of \eqref{eq:lift-elliptic-source}.
Then, for \(\tau\le t\le T\), we have
\begin{equation}\label{eq:lift-long-material-bound}
\begin{aligned}
 \|\sqrt\rho\dot\phi\|_2^2
 &\le C\log^2(e+t)\bigg[(\|\nabla u\|_2^2+\|P\|_2^2)^{1/2}\\
 &\quad+(\|\nabla u\|_2^2+\|P\|_2^2)
       (\|\sqrt\rho\dot u\|_2^2+\|P\|_3^3)^{1/4}\bigg]\\
 &\quad\times\bigl(\|\nabla\dot u\|_2^2+\|\sqrt\rho\dot\theta\|_2^2
                          +\|\sqrt\rho\dot u\|_2^4+\|P\|_4^4\bigr)\\
 &\quad+C\log(e+t)\|\nabla\phi\|_2
       \bigl(\|\sqrt\rho\dot u\|_2^4+\|P\|_4^4
                          +E_0\|\sqrt\rho\dot u\|_2^2\bigr).
\end{aligned}
\end{equation}
\end{lemma}
\begin{proof}
Noting that \(\dot P+P\Div u=\rho\dot\theta\) from
\eqref{eq:perturbation-system}$_1$, we obtain
\begin{equation}\label{eq:lift-laplacian-commutator}
\begin{aligned}
 &(\partial_t+u\cdot\nabla+\Div u)\bigl(2\mu|D(u)|^2+\lambda(\Div u)^2-P\Div u\bigr)\\
 =&4\mu D(u):\nabla\dot u+(2\lambda\Div u-P)\Div\dot u
          -\rho(\Div u)\dot\theta\\
 &+(\Div u)\bigl(2\mu|D(u)|^2+\lambda(\Div u)^2\bigr)
       -2\mu(\partial_iu^j+\partial_ju^i)\partial_iu^k\partial_ku^j\\
 &+(P-2\lambda\Div u)\partial_iu^j\partial_ju^i,
\end{aligned}
\end{equation}
and a direct calculation gives
\begin{equation}\label{eq:lift-laplacian-commutator-second}
\begin{aligned}
 &(\partial_t+u\cdot\nabla+\Div u)\Delta\phi
 =\Delta\dot\phi+\partial_i\bigl[(\Div u)\partial_i\phi
           -(\partial_iu^j+\partial_ju^i)\partial_j\phi\bigr].
\end{aligned}
\end{equation}
By \eqref{eq:perturbation-system}$_1$ and \eqref{eq:lift-elliptic-source}$_2$, we have
\[
 \int_{\R^2}\rho\dot\phi\dd x
 =\frac{\mathrm d}{\mathrm dt}\int_{\R^2}\rho\phi\dd x=0.
\]
Applying \(\partial_t+\Div(u\,\cdot)\) to \eqref{eq:lift-elliptic-source}$_1$,
we deduce from \eqref{eq:perturbation-system}$_1$ that
\begin{equation}\label{eq:lift-material-elliptic-equation}
\begin{aligned}
 -\kappa(\partial_t+\Div(u\,\cdot))\Delta\phi
 =&(\partial_t+\Div(u\,\cdot))\bigl(2\mu|D(u)|^2+\lambda(\Div u)^2-P\Div u\bigr)\\
 &-\frac{\rho}{m_0}\frac{\dd}{\dd t}
       \int_{\R^2}\bigl[2\mu|D(u)|^2+\lambda(\Div u)^2-P\Div u\bigr]\dd y .
\end{aligned}
\end{equation}
Multiplying \eqref{eq:lift-material-elliptic-equation} by \(z_0\in\mathcal H_\rho\)
and integrating by parts, we obtain from \eqref{eq:lift-laplacian-commutator-second}
and \(\int_{\R^2}\rho z_0\dd x=0\) that
\[
\begin{aligned}
 \kappa\int_{\R^2}\nabla\dot\phi\cdot\nabla z_0\dd x
 =&\int_{\R^2}z_0(\partial_t+u\cdot\nabla+\Div u)\bigl(2\mu|D(u)|^2+\lambda(\Div u)^2-P\Div u\bigr)\dd x\\
 &-\kappa\int_{\R^2}\bigl[(\Div u)\partial_i\phi-(\partial_i u^j+\partial_j u^i)\partial_j\phi\bigr]\partial_i z_0\dd x.
\end{aligned}
\]
Using \eqref{eq:initial-bootstrap-assumption}, \eqref{eq:lift-laplacian-commutator} and H\"older's inequality, we obtain, for every \(z_0\in\mathcal H_\rho\),
\begin{equation}\label{eq:lift-material-testing}
\begin{aligned}
 \kappa\left|\int_{\R^2}\nabla\dot\phi\cdot\nabla z_0\dd x\right|
 \le& C\int_{\R^2}\bigl(|\nabla u|\,|\nabla\dot u|+P|\nabla\dot u|
       +\rho|\Div u|\,|\dot\theta|+|\nabla u|^3+P|\nabla u|^2\bigr)|z_0|\dd x\\
 &+C\int_{\R^2}|\nabla u|\,|\nabla\phi|\,|\nabla z_0|\dd x\\
 \le& C\bigl(\|\nabla\dot u\|_2+\|\sqrt\rho\,\dot\theta\|_2+\|\nabla u\|_4^2\bigr)
       \left(\int_{\R^2}(|\nabla u|^2+P^2)z_0^2\dd x\right)^{1/2}\\
 &+C\||\nabla u|\nabla\phi\|_2\|\nabla z_0\|_2.
\end{aligned}
\end{equation}
Taking \(z_0=\dot\phi\) in \eqref{eq:lift-material-testing} and applying
\eqref{eq:lift-source-multiplier}, \eqref{eq:quartic-velocity-pressure} and
\eqref{eq:lift-gradient-product}, we obtain
\[
\begin{aligned}
 \|\nabla\dot\phi\|_2^2
 \le& C\log(e+t)\Bigl[(\|\nabla u\|_2^2+\|P\|_2^2)^{1/2}
      +(\|\nabla u\|_2^2+\|P\|_2^2)(\|\sqrt\rho\,\dot u\|_2^2+\|P\|_3^3)^{1/4}\Bigr]\\
 &\quad\times\bigl(\|\nabla\dot u\|_2^2+\|\sqrt\rho\,\dot\theta\|_2^2
      +\|\sqrt\rho\,\dot u\|_2^4+\|P\|_4^4\bigr)\\
 &\quad+C\|\nabla\phi\|_2\bigl(\|\sqrt\rho\,\dot u\|_2^4+\|P\|_4^4
      +E_0\|\sqrt\rho\,\dot u\|_2^2\bigr).
\end{aligned}
\]
Since \(\langle\dot\phi\rangle_\rho=0\), applying \eqref{eq:lift-log-anchor} with \(p=2\) gives
\[
 \|\sqrt\rho\,\dot\phi\|_2^2\le C\log(e+t)\|\nabla\dot\phi\|_2^2,
\]
which together with the preceding inequality yields \eqref{eq:lift-long-material-bound}.
\end{proof}

From \eqref{eq:perturbation-system}$_3$ and \eqref{eq:lift-elliptic-source},
we derive the following estimates for $\theta-\phi$.
\begin{lemma}\label{lem:temperature-residual-energy}
Under the assumptions of Proposition~\ref{prop:lower-order-closure}, let \(\phi\) be
the solution of \eqref{eq:lift-elliptic-source}.
For almost every \(t\in[0,T]\), we have
\begin{equation}\label{eq:lift-residual-elliptic}
 \frac{\|\nabla(\theta-\phi)\|_2^2}{\log(e+t)}
 +\|\nabla^2(\theta-\phi)\|_2^2+\|\nabla(\theta-\phi)\|_4^2
 \le C\Bigl\|\sqrt\rho\Bigl(\dot\theta-\frac{\dd}{\dd t}\langle\theta\rangle_\rho\Bigr)\Bigr\|_2^2,
\end{equation}
where \(C=C(m_0,\bar\rho,M,c_v,\kappa)\).
For \(0\le t\le T\), we also have
\begin{equation}\label{eq:lift-convection}
\begin{aligned}
 &\frac{\kappa}{2}\|\nabla(\theta-\phi)(t)\|_2^2
      +\frac{3c_v}{4}\int_0^t\Bigl\|\sqrt\rho\Bigl(\dot\theta-\frac{\dd}{\dd r}\langle\theta\rangle_\rho\Bigr)\Bigr\|_2^2\dd r\\
 \le&\frac{\kappa}{2}\|\nabla(\theta_0-\phi(0))\|_2^2+C\int_0^t\|\nabla u\|_2
      \Bigl\|\sqrt\rho\Bigl(\dot\theta-\frac{\dd}{\dd r}\langle\theta\rangle_\rho\Bigr)\Bigr\|_2^2\dd r\\
 &+c_v\int_0^t\|\sqrt\rho\,\dot\phi\|_2^2\dd r,
\end{aligned}
\end{equation}
where \(C=C(m_0,\bar\rho,c_v,\kappa)\).
Moreover, we have
\begin{equation}\label{eq:short-weighted-lift}
 \sup_{0\le t\le\min\{T,\tau\}}t^2
 \bigl(\|\sqrt\rho\,\dot u(t)\|_2^2+\|P(t)\|_3^3
                    +\|\nabla(\theta-\phi)(t)\|_2^2\bigr)\le CE_0^{1/12},
\end{equation}
where \(C=C(m_0,\bar\rho,M,\mu,\lambda,c_v,\kappa)\).
\end{lemma}
\begin{proof}
Differentiating \eqref{eq:physical-thermal-energy} and using
\eqref{eq:viscous-dissipation-identity}, \eqref{eq:mechanical-energy-testing}
and \eqref{eq:mean-material-acceleration}, we obtain
\[
\begin{aligned}
 c_vm_0\frac{\mathrm d}{\mathrm dt}\langle\theta\rangle_\rho
 =&\int_{\R^2}\bigl(2\mu|D(u)|^2+\lambda(\Div u)^2-P\Div u\bigr)\dd x\\
 =&-\int_{\R^2}\rho\bigl(u-\langle u\rangle_\rho\bigr)\cdot\dot u\dd x.
\end{aligned}
\]
Applying the Cauchy--Schwarz inequality, \eqref{eq:physical-total-energy} and \eqref{eq:lift-log-anchor}, we arrive at
\begin{equation}\label{eq:lift-temperature-mean-identity}
m_0\left|\frac{\mathrm d}{\mathrm d t}\langle\theta\rangle_\rho\right|^2
 \le C\min\{E_0,\log(e+t)\|\nabla u\|_2^2\}
                                      \|\sqrt\rho\,\dot u\|_2^2.
\end{equation}
By \eqref{eq:perturbation-system}$_1$, we have
\begin{equation}\label{eq:lift-mean-rate-zero}
 \int_{\R^2}\rho\bigl(\dot\theta-\frac{\dd}{\dd t}\langle\theta\rangle_\rho\bigr)\dd x=0
\end{equation}
and
\begin{equation}\label{eq:lift-mean-rate}
 \|\sqrt\rho\,\dot\theta\|_2^2
 =\int_{\R^2}\rho\bigl(\dot\theta-\frac{\dd}{\dd t}\langle\theta\rangle_\rho\bigr)^2\dd x
     +m_0\left|\frac{\dd}{\dd t}\langle\theta\rangle_\rho\right|^2.
\end{equation}
Using \eqref{eq:lift-log-anchor} and \eqref{eq:physical-total-energy}, we also get
\begin{equation}\label{eq:weighted-temperature-mean}
 \|\sqrt\rho\,\theta\|_2^2
 \le C\log(e+t)\|\nabla\theta\|_2^2
       +m_0\langle\theta\rangle_\rho^2
 \le C\log(e+t)\|\nabla\theta\|_2^2+CE_0^2.
\end{equation}

Combining \eqref{eq:perturbation-system}$_3$ with \eqref{eq:lift-elliptic-source} gives
\begin{equation}\label{eq:lift-residual-equation}
\kappa\Delta(\theta-\phi)=c_v\rho\bigl(\dot\theta-\frac{\dd}{\dd t}\langle\theta\rangle_\rho\bigr).
\end{equation}
Multiplying \eqref{eq:lift-residual-equation} by \(\theta-\phi-\langle\theta\rangle_\rho\),
integrating by parts and using \eqref{eq:lift-log-anchor}, we conclude that
\[
 \kappa\|\nabla(\theta-\phi)\|_2^2
 \le C\sqrt{\log(e+t)}\,\Bigl\|\sqrt\rho\Bigl(\dot\theta-\frac{\dd}{\dd t}\langle\theta\rangle_\rho\Bigr)\Bigr\|_2\|\nabla(\theta-\phi)\|_2.
\]
By \eqref{eq:lift-residual-equation}, Plancherel's identity and
\eqref{eq:riesz-potential-lp} with \(p=4/3\), we have
\[
\begin{aligned}
 &\|\nabla^2(\theta-\phi)\|_2+\|\nabla(\theta-\phi)\|_4\\
 \le& C\frac{c_v}{\kappa}\Biggl(\Bigl\|\rho\Bigl(\dot\theta-\frac{\dd}{\dd t}\langle\theta\rangle_\rho\Bigr)\Bigr\|_2
       +\Bigl\|\rho\Bigl(\dot\theta-\frac{\dd}{\dd t}\langle\theta\rangle_\rho\Bigr)\Bigr\|_{4/3}\Biggr)\\
 \le& C\frac{c_v}{\kappa}\bigl(\sqrt{2\bar\rho}+(2m_0\bar\rho)^{1/4}\bigr)
        \Bigl\|\sqrt\rho\Bigl(\dot\theta-\frac{\dd}{\dd t}\langle\theta\rangle_\rho\Bigr)\Bigr\|_2,
\end{aligned}
\]
where the last inequality uses H\"older's inequality, \eqref{eq:mass} and
\eqref{eq:initial-bootstrap-assumption}. This proves \eqref{eq:lift-residual-elliptic}.

Multiplying \eqref{eq:lift-residual-equation} by \(\dot\theta-\dot\phi\) and integrating by parts,
we deduce from \eqref{eq:lift-mean-rate-zero} that
\begin{equation}\label{eq:lift-residual-energy}
\begin{aligned}
 &\frac{\kappa}{2}\frac{\dd}{\dd t}\|\nabla(\theta-\phi)\|_2^2
       +c_v\int_{\R^2}\rho\bigl(\dot\theta-\frac{\dd}{\dd t}\langle\theta\rangle_\rho\bigr)^2\dd x\\
 =&c_v\int_{\R^2}\rho\bigl(\dot\theta-\frac{\dd}{\dd t}\langle\theta\rangle_\rho\bigr)\dot\phi\dd x
       -\kappa\int_{\R^2}\partial_i(\theta-\phi)\,\partial_i u^j\,
                                      \partial_j(\theta-\phi)\dd x\\
 &+\frac{\kappa}{2}\int_{\R^2}(\Div u)|\nabla(\theta-\phi)|^2\dd x.
\end{aligned}
\end{equation}
Using \eqref{eq:lift-residual-equation} and \eqref{eq:riesz-potential-lp} with \(p=4/3\), we obtain
\[
\|\nabla(\theta-\phi)\|_4^2
 \le C\frac{c_v^2}{\kappa^2}(2m_0\bar\rho)^{1/2}
       \int_{\R^2}\rho\bigl(\dot\theta-\frac{\dd}{\dd t}\langle\theta\rangle_\rho\bigr)^2\dd x.
\]
Young's inequality also gives
\[
 c_v\left|\int_{\R^2}\rho\bigl(\dot\theta-\frac{\dd}{\dd t}\langle\theta\rangle_\rho\bigr)\dot\phi\dd x\right|
 \le\frac{c_v}{4}\int_{\R^2}\rho\bigl(\dot\theta-\frac{\dd}{\dd t}\langle\theta\rangle_\rho\bigr)^2\dd x
       +c_v\|\sqrt\rho\,\dot\phi\|_2^2.
\]
Substituting these two estimates into \eqref{eq:lift-residual-energy} and integrating over \([0,t]\), we obtain \eqref{eq:lift-convection}.

Finally, for \(0\le t\le\min\{T,\tau\}\), taking \(z=\phi\) in \eqref{eq:short-spatial-multiplier}
and using \eqref{eq:lift-log-anchor}, \eqref{eq:lift-gradient-testing},
\eqref{eq:initial-cubic-pressure} and \eqref{eq:initial-bootstrap-assumption}, we conclude that
\[
 \|\nabla\phi\|_2^2+\|P\|_3^3
 \le C\bigl(\|\nabla u\|_2^2+\|\sqrt\rho\,\theta\|_2^2
             +\|\sqrt\rho\,\dot u\|_2^2+\|\nabla\theta\|_2^2\bigr).
\]
Multiplying this inequality by \(t^2\) and using \(\|\nabla(\theta-\phi)\|_2^2\le2\|\nabla\theta\|_2^2+2\|\nabla\phi\|_2^2\),
\eqref{eq:initial-weighted-lemma-second} and \eqref{eq:short-squared-weight}, we have
\[
\begin{aligned}
 &\sup_{0\le t\le\min\{T,\tau\}}t^2
  \bigl(\|\sqrt\rho\,\dot u(t)\|_2^2+\|P(t)\|_3^3
        +\|\nabla(\theta-\phi)(t)\|_2^2\bigr)\\
 \le& C\sup_{0\le t\le\min\{T,\tau\}}t^2
  \bigl(\|\nabla u(t)\|_2^2+\|\sqrt\rho\,\theta(t)\|_2^2
        +\|\sqrt\rho\,\dot u(t)\|_2^2+\|\nabla\theta(t)\|_2^2\bigr)\\
 \le& CE_0^{1/12},
\end{aligned}
\]
which proves \eqref{eq:short-weighted-lift}.
\end{proof}

Combining
Lemmas~\ref{lem:mechanical-and-pressure-energy}--\ref{lem:temperature-residual-energy},
we establish the following time-weighted estimates on $[\tau,T]$
(see \cite[Lemma~3.4]{LiXin2019} for the barotropic case).
\begin{lemma}\label{lem:weighted-decay}
Under the assumptions of Proposition~\ref{prop:lower-order-closure},
let \(T\ge\tau\). Then
\begin{equation}\label{eq:weighted-decay-lemma}
\begin{aligned}
 &\sup_{\tau\le t\le T}(1+t)^{3/2}
  \bigl(\|\sqrt\rho\dot u\|_2^2+\|P\|_3^3+\|\nabla(\theta-\phi)\|_2^2\bigr)\\
 &\quad+\int_{\tau}^T (1+t)^{3/2}
  \bigl(\|\nabla\dot u\|_2^2+\|\sqrt\rho\dot\theta\|_2^2+
                        \|P\|_4^4\bigr)\dd t\le CE_0^{1/32},
\end{aligned}
\end{equation}
where $C=C(m_0,\bar\rho,M,\mu,\lambda,c_v,\kappa)$.
Moreover, for \(\tau\le t\le T\),
\begin{equation}\label{eq:decay-mechanical-second}
 \|\nabla u(t)\|_2^2+\|P(t)\|_2^2
 \le CE_0^{33/64}(1+t)^{-3/4},
\end{equation}
and
\begin{equation}\label{eq:decay-thermal}
 \|\nabla\theta(t)\|_2^2\le CE_0^{1/32}(1+t)^{-3/2}
          +CE_0^{99/128}\log(e+t)(1+t)^{-9/8},
\end{equation}
where $C=C(m_0,\bar\rho,M,\mu,\lambda,c_v,\kappa)$.
For each $2\le p\le3$, we also have
\begin{equation}\label{eq:decay-p-family}
 \sup_{\tau\le t\le T}(1+t)^{3(p-1)/4}
       \bigl(\|\nabla u\|_p^p+\|P\|_p^p\bigr)
       \le CE_0^{(95-31p)/64},
\end{equation}
where $C=C(m_0,\bar\rho,M,\mu,\lambda,c_v,\kappa,p)$.
\end{lemma}
\begin{proof}
With \(C_0\) given in Lemma~\ref{lem:mechanical-and-pressure-energy}, we choose
\begin{equation}\label{eq:combined-lower-energy}
 L=\max\left\{1,4C_0(2\mu+\lambda),\frac{16(C_0+1)}{c_v},\sqrt{\frac{2}{\kappa}}\right\}.
\end{equation}
For \(\tau\le t\le T\), we set
\begin{equation}\label{eq:decay-weighted-energy}
 Y(t)=(1+t)^{3/2}\left(\|\sqrt\rho\,\dot u\|_2^2+L\|P\|_3^3
                    +\frac{\kappa L^2}{2}\|\nabla(\theta-\phi)\|_2^2\right).
\end{equation}
Taking \(p=2\) and \(w=\dot u\) in \eqref{eq:lift-log-anchor} and using
\eqref{eq:mean-material-acceleration}, we have
\begin{equation}\label{eq:decay-acceleration-anchor}
 \|\sqrt\rho\,\dot u\|_2^2\le C\log(e+t)\|\nabla\dot u\|_2^2.
\end{equation}
It follows from \eqref{eq:lift-mechanical-interpolation} and \eqref{eq:decay-weighted-energy} that
\begin{equation}\label{eq:decay-mechanical-energy}
 \|\nabla u\|_2^2+\|P\|_2^2
 \le CE_0^{1/2}(1+t)^{-3/4}Y(t)^{1/2}.
\end{equation}

Multiplying \(\dot P+P\Div u=\rho\dot\theta\) by \(3LP^2\), integrating over \(\R^2\)
and adding \eqref{eq:material-acceleration-products}, we deduce from \eqref{eq:flux-definitions},
\eqref{eq:div-curl-estimate}, \eqref{eq:flux-interpolation} with \(p=4\),
\eqref{eq:initial-bootstrap-assumption} and Young's inequality that
\begin{equation}\label{eq:decay-mechanical-differential}
\begin{aligned}
 &\frac{\dd}{\dd t}\left(\|\sqrt\rho\,\dot u\|_2^2+L\|P\|_3^3\right)
       +\mu\|\nabla\dot u\|_2^2+\frac{L}{2(2\mu+\lambda)}\|P\|_4^4\\
 \le&C_0(1+L)\|\sqrt\rho\,\dot\theta\|_2^2
        +C(\|\nabla u\|_2^2+\|P\|_2^2)\|\sqrt\rho\,\dot u\|_2^2.
\end{aligned}
\end{equation}
Taking \(z=\phi\) in \eqref{eq:lift-source-multiplier} and substituting the result into
\eqref{eq:lift-gradient-testing}, we obtain
\begin{equation}\label{eq:decay-correction-gradient}
\begin{aligned}
 \|\nabla\phi\|_2^2
 \le C(\|\nabla u\|_2^2+\|P\|_2^2)^{3/2}
 \times\left[\log(e+t)+(\|\nabla u\|_2^2+\|P\|_2^2)^{1/2}
                   (\|\sqrt\rho\,\dot u\|_2^2+\|P\|_3^3)^{1/4}\right].
\end{aligned}
\end{equation}
Substituting \eqref{eq:decay-weighted-energy}--\eqref{eq:decay-mechanical-energy}
and \eqref{eq:decay-correction-gradient} into \eqref{eq:lift-long-material-bound}
and \eqref{eq:lift-temperature-mean-identity}, we arrive at
\begin{equation}\label{eq:decay-correction-remainder}
\begin{aligned}
 \|\sqrt\rho\,\dot\phi\|_2^2
       +m_0\left|\frac{\dd}{\dd t}\langle\theta\rangle_\rho\right|^2
 \le CE_0^{1/4}\frac{\log^3(e+t)}{(1+t)^{3/8}}(1+Y(t))^2
      \bigl(\|\nabla\dot u\|_2^2+\|\sqrt\rho\,\dot\theta\|_2^2+\|P\|_4^4\bigr).
\end{aligned}
\end{equation}
Using \eqref{eq:lift-residual-elliptic}, \eqref{eq:lift-mean-rate} and Young's inequality
in \eqref{eq:lift-residual-energy}, we get
\[
\begin{aligned}
 &\frac{\kappa}{2}\frac{\dd}{\dd t}\|\nabla(\theta-\phi)\|_2^2
                       +\frac{3c_v}{4}\|\sqrt\rho\,\dot\theta\|_2^2
 \le c_v\|\sqrt\rho\,\dot\phi\|_2^2
       +\frac{3c_vm_0}{4}\left|\frac{\dd}{\dd t}\langle\theta\rangle_\rho\right|^2
       +C\|\nabla u\|_2\|\sqrt\rho\,\dot\theta\|_2^2.
\end{aligned}
\]
Multiplying this inequality by \(L^2\) and adding \eqref{eq:decay-mechanical-differential},
we deduce from \eqref{eq:combined-lower-energy}, \eqref{eq:decay-acceleration-anchor},
\eqref{eq:decay-mechanical-energy}, \eqref{eq:decay-correction-remainder} and
\(\sup_{t\ge0}(1+t)^{-3/8}\log^3(e+t)<\infty\) that
\begin{equation}\label{eq:decay-coupled-remainder}
\begin{aligned}
 &\frac{\dd}{\dd t}\left(\frac{Y(t)}{(1+t)^{3/2}}\right)
       +\left[2-CE_0^{1/4}(1+Y(t))^2\right]\times\left(\frac{\mu}{2}\|\nabla\dot u\|_2^2+\|\sqrt\rho\,\dot\theta\|_2^2
                      +\frac{L}{4(2\mu+\lambda)}\|P\|_4^4\right)\le0.
\end{aligned}
\end{equation}

By \eqref{eq:decay-acceleration-anchor}, \eqref{eq:lift-residual-elliptic},
\eqref{eq:lift-mean-rate} and Young's inequality, we have
\begin{equation}\label{eq:lift-time-weight-absorption}
\begin{aligned}
 \frac{3}{2(1+t)}\left(\|\sqrt\rho\,\dot u\|_2^2
                   +\frac{\kappa L^2}{2}\|\nabla(\theta-\phi)\|_2^2\right)
 \le\frac{\mu}{4}\|\nabla\dot u\|_2^2+\frac12\|\sqrt\rho\,\dot\theta\|_2^2
             +\frac{C\log(e+t)}{(1+t)^{7/2}}Y(t).
\end{aligned}
\end{equation}
It follows from H\"older's inequality, \eqref{eq:physical-total-energy} and Young's inequality that
\begin{equation}\label{eq:lift-time-weight-pressure}
\begin{aligned}
 \frac{3L}{2}(1+t)^{1/2}\|P\|_3^3
 \le&\frac{3L}{2}(1+t)^{1/2}(E_0/c_v)^{1/3}\|P\|_4^{8/3}\\
 \le&\frac{L}{8(2\mu+\lambda)}(1+t)^{3/2}\|P\|_4^4
        +\frac{CE_0}{(1+t)^{3/2}}.
\end{aligned}
\end{equation}
Multiplying \eqref{eq:decay-coupled-remainder} by \((1+t)^{3/2}\) and using
\eqref{eq:lift-time-weight-absorption}--\eqref{eq:lift-time-weight-pressure}, we conclude that
\begin{equation}\label{eq:lift-time-weighted-differential}
\begin{aligned}
 &\frac{\dd Y}{\dd t}
       +\left[1-CE_0^{1/4}(1+Y)^2\right](1+t)^{3/2}
        \left(\frac{\mu}{2}\|\nabla\dot u\|_2^2+\|\sqrt\rho\,\dot\theta\|_2^2
                       +\frac{L}{4(2\mu+\lambda)}\|P\|_4^4\right)\\
 \le&\frac{C}{(1+t)^{3/2}}(Y+E_0).
\end{aligned}
\end{equation}
Fixing a common constant \(C\) in \eqref{eq:decay-coupled-remainder} and
\eqref{eq:lift-time-weighted-differential}, we set
\[
 Z(t)=(Y(\tau)+E_0)\exp\left(C\int_\tau^t\frac{\dd s}{(1+s)^{3/2}}\right)-E_0.
\]
By \eqref{eq:short-weighted-lift} and \eqref{eq:decay-weighted-energy}, we have
\begin{equation}\label{eq:decay-comparison}
 0\le Z(t)\le e^{2C}(Y(\tau)+E_0)
       \le C(1+\tau)^{3/2}\tau^{-2}E_0^{1/12}\le CE_0^{1/12}.
\end{equation}
Choosing \(\eps_0\) so that \(C\eps_0^{1/12}\le1\) in \eqref{eq:decay-comparison} and
\(8C\eps_0^{1/4}\le1\) with \(C\) from \eqref{eq:lift-time-weighted-differential}, we get
\begin{equation}\label{eq:decay-comparison-smallness}
 CE_0^{1/4}(1+Z(t))^2\le4CE_0^{1/4}\le\tfrac12,
 \qquad \tau\le t\le T.
\end{equation}
Since \(\frac{\dd Z}{\dd t}=C(1+t)^{-3/2}(Z+E_0)\), it follows from
\eqref{eq:lift-time-weighted-differential} and \eqref{eq:decay-comparison-smallness} that
\[
\begin{aligned}
 \frac{\dd}{\dd t}(Y-Z)_+
 \le&C\bigg[(1+t)^{-3/2}+E_0^{1/4}(2+Y+Z)(1+t)^{3/2}\\
 &\qquad\times\bigl(\|\nabla\dot u\|_2^2+\|\sqrt\rho\,\dot\theta\|_2^2
                         +\|P\|_4^4\bigr)\bigg](Y-Z)_+.
\end{aligned}
\]
By \eqref{eq:regularity}, \eqref{eq:gn-lp}, \eqref{eq:density-anchor} and
\eqref{eq:decay-correction-gradient}, we have
\[
 \int_\tau^T\Bigl[(1+t)^{-3/2}+E_0^{1/4}(2+Y+Z)(1+t)^{3/2}
      \bigl(\|\nabla\dot u\|_2^2+\|\sqrt\rho\,\dot\theta\|_2^2+\|P\|_4^4\bigr)\Bigr]\dd t<\infty.
\]
Since \(Y(\tau)=Z(\tau)\), Gronwall's inequality gives
\[
 Y(t)\le Z(t),\qquad \tau\le t\le T.
\]
Integrating \eqref{eq:lift-time-weighted-differential} over \([\tau,t]\) and using
\eqref{eq:decay-comparison}--\eqref{eq:decay-comparison-smallness}, we arrive at
\[
 Y(t)+\int_\tau^t(1+s)^{3/2}
   \bigl(\|\nabla\dot u\|_2^2+\|\sqrt\rho\,\dot\theta\|_2^2+\|P\|_4^4\bigr)\dd s
 \le C Z(t)\le C E_0^{1/12},
\]
which proves \eqref{eq:weighted-decay-lemma} since \(E_0\le1\).
Substituting \(Y\le Z\) and \eqref{eq:decay-comparison-smallness} into
\eqref{eq:decay-coupled-remainder}, we also have, on \([\tau,T]\),
\begin{equation}\label{eq:lift-joint-differential}
\begin{aligned}
 &\frac{\dd}{\dd t}\left(\|\sqrt\rho\,\dot u\|_2^2+L\|P\|_3^3
                       +\frac{\kappa L^2}{2}\|\nabla(\theta-\phi)\|_2^2\right)+\frac{\mu}{2}\|\nabla\dot u\|_2^2+\|\sqrt\rho\,\dot\theta\|_2^2
                       +\frac{L}{4(2\mu+\lambda)}\|P\|_4^4\le0.
\end{aligned}
\end{equation}

Substituting \eqref{eq:weighted-decay-lemma} into \eqref{eq:lift-mechanical-interpolation}
yields \eqref{eq:decay-mechanical-second}.
By \eqref{eq:decay-correction-gradient}, \eqref{eq:weighted-decay-lemma} and
\eqref{eq:decay-mechanical-second}, we get
\[
 \|\nabla\phi(t)\|_2^2\le CE_0^{99/128}\log(e+t)(1+t)^{-9/8},
\]
which, combined with \eqref{eq:weighted-decay-lemma}, gives
\[
 \|\nabla\theta(t)\|_2^2
 \le 2\|\nabla(\theta-\phi)(t)\|_2^2+2\|\nabla\phi(t)\|_2^2
 \le CE_0^{1/32}(1+t)^{-3/2}+CE_0^{99/128}\log(e+t)(1+t)^{-9/8},
\]
and proves \eqref{eq:decay-thermal}.

For \(2\le p\le3\), H\"older's inequality gives
\begin{equation}\label{eq:decay-p-interpolation}
 \|P\|_p^p\le(\|P\|_2^2)^{3-p}(\|P\|_3^3)^{p-2},
\end{equation}
and taking \(2\le p\le3\) in \eqref{eq:div-curl-estimate} and \eqref{eq:flux-interpolation}
together with \eqref{eq:initial-bootstrap-assumption} yields
\begin{equation}\label{eq:decay-p-interpolation-second}
 \|\nabla u\|_p^p
 \le C\|\sqrt\rho\,\dot u\|_2^{p-2}
                 (\|\nabla u\|_2+\|P\|_2)^2+C\|P\|_p^p.
\end{equation}
Substituting \eqref{eq:weighted-decay-lemma} and \eqref{eq:decay-mechanical-second}
into \eqref{eq:decay-p-interpolation}--\eqref{eq:decay-p-interpolation-second}, we conclude that
\[
\begin{aligned}
 \|\nabla u\|_p^p+\|P\|_p^p
 \le& C\bigl(E_0^{1/2}+E_0^{(3-p)/2}\bigr)
       E_0^{(p-1)/64}(1+t)^{-3(p-1)/4}\\
 \le& CE_0^{(95-31p)/64}(1+t)^{-3(p-1)/4},
\end{aligned}
\]
which proves \eqref{eq:decay-p-family}.
\end{proof}

Based on Lemma~\ref{lem:weighted-decay}, we obtain the following improved
time-weighted decay estimates on $[\tau,T]$.
\begin{lemma}\label{lem:endpoint-weighted-decay}
Under the assumptions of Proposition~\ref{prop:lower-order-closure}, if \(T>\tau\), then we have
\begin{equation}\label{eq:endpoint-weighted-lemma}
 \sup_{\tau\le t\le T}(1+t)^2
  \bigl(\|\sqrt\rho\,\dot u\|_2^2+\|P\|_3^3+\|\nabla(\theta-\phi)\|_2^2\bigr)
                                                    \le CE_0^{1/32},
\end{equation}
\begin{equation}\label{eq:endpoint-weighted-lemma-second}
 \int_{\tau}^T\frac{(1+t)^2}{\log^2(e+t)}
  \bigl(\|\nabla\dot u\|_2^2+\|\sqrt\rho\,\dot\theta\|_2^2+
                    \|P\|_4^4\bigr)\dd t\le CE_0^{1/32},
\end{equation}
and, for \(\tau\le t\le T\),
\begin{equation}\label{eq:endpoint-weighted-lemma-third}
 \|\nabla\theta(t)\|_2^2
  \le CE_0^{1/32}(1+t)^{-2}
            +CE_0^{99/128}\log(e+t)(1+t)^{-3/2},
\end{equation}
where $C=C(m_0,\bar\rho,M,\mu,\lambda,c_v,\kappa)$.
For each $2\le p\le3$, we also have
\begin{equation}\label{eq:endpoint-p-lemma}
 \sup_{\tau\le t\le T}(1+t)^{p-1}
       \bigl(\|\nabla u\|_p^p+\|P\|_p^p\bigr)
             \le CE_0^{(95-31p)/64},
\end{equation}
and
\begin{equation}\label{eq:endpoint-p-lemma-second}
 \int_{\tau}^T\frac{(1+t)^{p-1}}{\log^2(e+t)}\|P\|_{p+1}^{p+1}\dd t
             \le CE_0^{(96-31p)/96},
\end{equation}
where $C=C(m_0,\bar\rho,M,\mu,\lambda,c_v,\kappa,p)$.
\end{lemma}
\begin{proof}
By \eqref{eq:weighted-decay-lemma}, $E_0\le E_0^{1/12}$, and
\(\sup_{t\ge0}\log(e+t) (1+t)^{-3/4}<\infty\), we have
\[
 \log(e+t)\bigl(\|\sqrt\rho\,\dot u\|_2+\|P\|_3^{3/2}
                    +\|\nabla(\theta-\phi)\|_2\bigr)
                    +E_0^{1/2}\le CE_0^{1/64}.
\]
Taking \(w=\dot u\) and \(p=2\) in \eqref{eq:lift-log-anchor} and using
\eqref{eq:mean-material-acceleration}, \eqref{eq:lift-residual-elliptic},
\eqref{eq:lift-mean-rate} and \eqref{eq:lift-time-weight-pressure}, we obtain
\begin{equation}\label{eq:endpoint-scalar-coercivity}
\begin{aligned}
 &\left(\|\sqrt\rho\,\dot u\|_2^2+L\|P\|_3^3
          +\frac{L^2\kappa}{2}\|\nabla(\theta-\phi)\|_2^2\right)^{3/2}\\
 \le& C\left[\left(\|\sqrt\rho\,\dot u\|_2^2
          +\frac{L^2\kappa}{2}\|\nabla(\theta-\phi)\|_2^2\right)^{3/2}
          +(L\|P\|_3^3)^{3/2}\right]\\
 \le& CE_0^{1/64}\left(\frac\mu2\|\nabla\dot u\|_2^2
          +\|\sqrt\rho\,\dot\theta\|_2^2
          +\frac{L}{4(2\mu+\lambda)}\|P\|_4^4\right).
\end{aligned}
\end{equation}
Combining \eqref{eq:lift-joint-differential} and
\eqref{eq:endpoint-scalar-coercivity} gives
\[
\begin{aligned}
 &\frac{\dd}{\dd t}\left(\|\sqrt\rho\,\dot u\|_2^2+L\|P\|_3^3
                         +\frac{L^2\kappa}{2}\|\nabla(\theta-\phi)\|_2^2\right)\\
 &+cE_0^{-1/64}\left(\|\sqrt\rho\,\dot u\|_2^2+L\|P\|_3^3
                         +\frac{L^2\kappa}{2}\|\nabla(\theta-\phi)\|_2^2\right)^{3/2}
 \le0.
\end{aligned}
\]
Integrating over \([\tau,t]\) and using \eqref{eq:weighted-decay-lemma}, we conclude that
\begin{equation}\label{eq:endpoint-scalar-decay}
\begin{aligned}
 &\|\sqrt\rho\,\dot u(t)\|_2^2+L\|P(t)\|_3^3
                         +\frac{L^2\kappa}{2}\|\nabla(\theta-\phi)(t)\|_2^2\\
 \le& E_0^{1/32}
              \left[C^{-1/2}+\frac c2(t-\tau)\right]^{-2}\\
 \le& CE_0^{1/32}(1+t-\tau)^{-2}
 \le CE_0^{1/32}(1+t)^{-2},
\end{aligned}
\end{equation}
which yields \eqref{eq:endpoint-weighted-lemma}.

To get \eqref{eq:endpoint-weighted-lemma-second}, we observe that
\[
 0\le\frac{\dd}{\dd t}\log\frac{(1+t)^2}{\log^2(e+t)}
 =\frac{2}{1+t}-\frac2{(e+t)\log(e+t)}\le\frac{2}{1+t}.
\]
The first inequality in \eqref{eq:lift-time-weight-pressure} and Young's inequality give
\begin{equation}\label{eq:endpoint-pressure-remainder}
\begin{aligned}
 L\|P\|_3^3\,\frac{\mathrm d}{\mathrm d t}\left(\frac{(1+t)^2}{\log^2(e+t)}\right)
 \le&\frac{2L(1+t)}{\log^2(e+t)}(E_0/c_v)^{1/3}\|P\|_4^{8/3}\\
 \le&\frac{L(1+t)^2}{8(2\mu+\lambda)\log^2(e+t)}\|P\|_4^4+\frac{CE_0}{(1+t)\log^2(e+t)}.
\end{aligned}
\end{equation}
Multiplying \eqref{eq:lift-joint-differential} by \((1+t)^2/\log^2(e+t)\)
and applying Young's inequality, we deduce from \eqref{eq:lift-residual-elliptic},
\eqref{eq:lift-mean-rate}, \eqref{eq:decay-acceleration-anchor} and \eqref{eq:endpoint-pressure-remainder} that
\begin{equation}\label{eq:endpoint-weighted-differential}
\begin{aligned}
 &\frac{\dd}{\dd t}\left[\frac{(1+t)^2}{\log^2(e+t)}
       \left(\|\sqrt\rho\,\dot u\|_2^2+L\|P\|_3^3
                         +\frac{L^2\kappa}{2}\|\nabla(\theta-\phi)\|_2^2\right)\right]\\
 &+\frac{(1+t)^2}{\log^2(e+t)}
       \left(\frac\mu4\|\nabla\dot u\|_2^2
        +\frac12\|\sqrt\rho\,\dot\theta\|_2^2
        +\frac{L}{8(2\mu+\lambda)}\|P\|_4^4\right)\\
 \le&\frac{C\log(e+t)}{(1+t)^2}\frac{(1+t)^2}{\log^2(e+t)}
       \left(\|\sqrt\rho\,\dot u\|_2^2+L\|P\|_3^3
                         +\frac{L^2\kappa}{2}\|\nabla(\theta-\phi)\|_2^2\right)\\
 &+\frac{CE_0}{(1+t)\log^2(e+t)}.
\end{aligned}
\end{equation}
Notice that
\[
 \int_0^\infty\frac{\dd t}{(1+t)\log^2(e+t)}\le e,
 \qquad \frac{(1+\tau)^2}{\log^2(e+\tau)}\le4.
\]
Taking \(t=\tau\) in \eqref{eq:endpoint-scalar-decay} and applying Gronwall's inequality
to \eqref{eq:endpoint-weighted-differential} on \([\tau,T]\), we obtain \eqref{eq:endpoint-weighted-lemma-second}.

For \(\tau\le t\le T\), substituting \eqref{eq:endpoint-scalar-decay} into \eqref{eq:lift-mechanical-interpolation} gives
\begin{equation}\label{eq:endpoint-mechanical-l2}
 \|\nabla u(t)\|_2^2+\|P(t)\|_2^2
 \le CE_0^{33/64}(1+t)^{-1}.
\end{equation}
Substituting \eqref{eq:endpoint-scalar-decay} and \eqref{eq:endpoint-mechanical-l2} into
\eqref{eq:decay-p-interpolation}--\eqref{eq:decay-p-interpolation-second}, we get
\[
 \|\nabla u\|_p^p+\|P\|_p^p
 \le CE_0^{(95-31p)/64}(1+t)^{-(p-1)},
 \qquad 2\le p\le3,
\]
which proves \eqref{eq:endpoint-p-lemma}.
By \eqref{eq:decay-correction-gradient}, \eqref{eq:endpoint-scalar-decay} and \eqref{eq:endpoint-mechanical-l2}, we have
\begin{equation}\label{eq:endpoint-lift-gradient}
 \|\nabla\phi(t)\|_2^2\le CE_0^{99/128}\log(e+t)(1+t)^{-3/2}.
\end{equation}
Combining \eqref{eq:endpoint-scalar-decay} and \eqref{eq:endpoint-lift-gradient}, we conclude that
\[
 \|\nabla\theta(t)\|_2^2
 \le 2\|\nabla(\theta-\phi)(t)\|_2^2+2\|\nabla\phi(t)\|_2^2
 \le CE_0^{1/32}(1+t)^{-2}+CE_0^{99/128}\log(e+t)(1+t)^{-3/2},
\]
which gives \eqref{eq:endpoint-weighted-lemma-third}.

Finally, for \(2\le p<3\), H\"older's inequality and \eqref{eq:physical-total-energy} yield
\[
 \|P\|_{p+1}^{p+1}\le(E_0/c_v)^{(3-p)/3}(\|P\|_4^4)^{p/3},
\]
which implies that
\[
 \frac{(1+t)^{p-1}}{\log^2(e+t)}\|P\|_{p+1}^{p+1}
 \le(E_0/c_v)^{(3-p)/3}\left(\frac{(1+t)^2}{\log^2(e+t)}\|P\|_4^4\right)^{p/3}
  \left(\frac1{(1+t)\log^2(e+t)}\right)^{(3-p)/3}.
\]
Applying H\"older's inequality in time and using \eqref{eq:endpoint-weighted-lemma-second}, we arrive at
\begin{equation}\label{eq:endpoint-pressure-integral}
\begin{aligned}
 &\int_{\tau}^T\frac{(1+t)^{p-1}}{\log^2(e+t)}\|P\|_{p+1}^{p+1}\dd t\\
 \le&CE_0^{(3-p)/3}
       \left(\int_{\tau}^T\frac{(1+t)^2}{\log^2(e+t)}
                                  \|P\|_4^4\dd t\right)^{p/3}
       \left(\int_{\tau}^T\frac{\dd t}{(1+t)\log^2(e+t)}\right)^{(3-p)/3}\\
 \le& CE_0^{(96-31p)/96},
\end{aligned}
\end{equation}
which gives \eqref{eq:endpoint-p-lemma-second} for \(2\le p<3\). The case \(p=3\) follows directly from \eqref{eq:endpoint-weighted-lemma-second}.
\end{proof}

\begin{remark}\label{rem:full-pressure-decay}
It should be noted that, unlike the barotropic case in \cite[(3.27)--(3.38)]{LiXin2019},
it follows from \eqref{eq:perturbation-system}$_1$ and \eqref{eq:perturbation-system}$_3$ that
the pressure \(P=\rho\theta\) satisfies
\[
 \dot P+\left(1+\frac{1}{c_v}\right)P\Div u
 =\frac{\kappa}{c_v}\Delta\theta
  +\frac{1}{c_v}\bigl(2\mu|D(u)|^2+\lambda(\Div u)^2\bigr),
\]
so the barotropic iteration in \cite[(3.27)--(3.38)]{LiXin2019} cannot be directly applied here.
Instead, multiplying \(\dot P+P\Div u=\rho\dot\theta\) by \(pP^{p-1}\), integrating over \(\R^2\)
and using \eqref{eq:flux-definitions}, we have
\begin{equation}\label{eq:full-pressure-power-balance}
 \frac{\dd}{\dd t}\|P\|_p^p+\frac{p-1}{2\mu+\lambda}\|P\|_{p+1}^{p+1}
 =-\frac{p-1}{2\mu+\lambda}\int_{\R^2}P^pF\dd x
   +p\int_{\R^2}\rho P^{p-1}\dot\theta\dd x.
\end{equation}
When \(p=3\), applying the Cauchy--Schwarz and Young inequalities to the last term in
\eqref{eq:full-pressure-power-balance}, we obtain
\[
 3\left|\int_{\R^2}\rho P^2\dot\theta\dd x\right|
 \le C\|P\|_4^2\|\sqrt\rho\,\dot\theta\|_2
 \le \frac{1}{2\mu+\lambda}\|P\|_4^4+C\|\sqrt\rho\,\dot\theta\|_2^2,
\]
where \(\frac{1}{2\mu+\lambda}\|P\|_4^4\) can be absorbed by the left-hand side of
\eqref{eq:full-pressure-power-balance}, as in \eqref{eq:cubic-pressure-differential}.
However, for \(p>3\), by the Cauchy--Schwarz and Young inequalities, we only get
\[
 p\left|\int_{\R^2}\rho P^{p-1}\dot\theta\dd x\right|
 \le C\|P\|_{2p-2}^{p-1}\|\sqrt\rho\,\dot\theta\|_2
 \le C\|P\|_{2p-2}^{2p-2}+C\|\sqrt\rho\,\dot\theta\|_2^2,
\]
where \(2p-2>p+1\), so that \(\|P\|_{2p-2}^{2p-2}\) cannot be absorbed by \(\|P\|_{p+1}^{p+1}\).
\end{remark}

\begin{proof}[Proof of Proposition~\ref{prop:lower-order-closure}]
\phantomsection\label{proof:lower-order-closure}
On \([0,\min\{T,\tau\}]\), the first inequality in \eqref{eq:initial-bootstrap-improvement}
follows from Lemma~\ref{lem:initial-density-control} and the second from
\eqref{eq:initial-weighted-temperature-gradient}, since \(\sigma(t)=t\).

Suppose now that \(T>\tau\). By \eqref{eq:weighted-temperature-mean},
\eqref{eq:decay-mechanical-second} and \eqref{eq:decay-thermal}, we obtain
\begin{equation}\label{eq:lower-rough-improvement}
 \sup_{\tau\le t\le T}\bigl(\|\nabla u\|_2^2+\|\sqrt\rho\,\theta\|_2^2+\|\nabla\theta\|_2^2\bigr)
 \le CE_0^{1/32}.
\end{equation}
Using \eqref{eq:weighted-decay-lemma}, \eqref{eq:decay-thermal} and
\(\log(e+t)(1+t)^{-9/8}\le C(1+t)^{-17/16}\), we arrive at
\begin{equation}\label{eq:lower-rough-improvement-second}
\begin{aligned}
 \int_{\tau}^T\bigl(\|\sqrt\rho\,\dot u\|_2^2+\|\nabla\theta\|_2^2\bigr)\dd t
 \le& CE_0^{1/32}\int_{\tau}^T\Bigl((1+t)^{-3/2}+(1+t)^{-17/16}\Bigr)\dd t\\
 \le& CE_0^{1/32}.
\end{aligned}
\end{equation}
Next, taking \(p=2\) and \(w=\dot u\) in \eqref{eq:lift-log-anchor}, and using
\eqref{eq:mean-material-acceleration} together with \eqref{eq:lift-flux-supremum}, we have
\[
 \|F\|_\infty
 \le C\log^{1/6}(e+t)\|\sqrt\rho\,\dot u\|_2^{2/3}\|\nabla\dot u\|_2^{1/3}
 \le C\sqrt{\log(e+t)}\,\|\nabla\dot u\|_2.
\]
Applying the Cauchy--Schwarz inequality and \eqref{eq:endpoint-weighted-lemma-second}, we get
\begin{equation}\label{eq:late-log-weighted-flux}
\begin{aligned}
 \int_{\tau}^T\|F\|_\infty\dd t
 \le& C\left(\int_{\tau}^T\frac{(1+t)^2}{\log^2(e+t)}\|\nabla\dot u\|_2^2\dd t\right)^{1/2}
      \left(\int_0^\infty\frac{\log^3(e+t)}{(1+t)^2}\dd t\right)^{1/2}\\
 \le& CE_0^{1/64}.
\end{aligned}
\end{equation}
Solving \eqref{eq:perturbation-system}$_1$ along the flow \eqref{eq:particle-flow} on \([\tau,t]\),
and using \(P\ge0\), \eqref{eq:initial-data-bound}, \eqref{eq:flux-definitions},
\eqref{eq:initial-density-gain}, \eqref{eq:late-log-weighted-flux} and \(E_0^{1/8}\le E_0^{1/64}\le E_0^{1/96}\),
we conclude that
\begin{equation}\label{eq:density-characteristics}
\begin{aligned}
 \|\rho(t)\|_\infty
 \le&\|\rho(\tau)\|_\infty
       \exp\left\{\frac{1}{2\mu+\lambda}\int_{\tau}^t\|F\|_\infty\dd s\right\}\\
 \le&\bar\rho\exp\bigl(CE_0^{1/96}\bigr).
\end{aligned}
\end{equation}
Let \(C_1=C_1(m_0,\bar\rho,M,\mu,\lambda,c_v,\kappa)\) be the largest of the constants \(C\)
appearing in \eqref{eq:initial-density-gain}, the proof of Lemma~\ref{lem:mechanical-and-pressure-energy},
\eqref{eq:lift-time-weighted-differential}, \eqref{eq:decay-comparison},
\eqref{eq:lower-rough-improvement} and \eqref{eq:density-characteristics}.
We choose \(\eps_0\in(0,1]\) such that
\begin{equation}\label{eq:original-smallness-threshold}
 C_1\eps_0^{1/96}\le\min\left\{\frac{1}{8},\,K,\,\log\frac{3}{2}\right\}.
\end{equation}
Since \(0<E_0\le\eps_0\le1\), \eqref{eq:original-smallness-threshold} implies that
\[
 C_1\eps_0^{1/8}\le\log\frac{3}{2},\qquad
 C_1\eps_0\le\frac{1}{4},\qquad
 C_1\eps_0^{1/12}\le1,\qquad
 8C_1\eps_0^{1/4}\le1,
\]
and it follows from \eqref{eq:lower-rough-improvement}, \eqref{eq:density-characteristics},
\eqref{eq:original-smallness-threshold} and \(\sigma(t)\le1\) that
\[
 \sup_{\tau\le t\le T}\|\rho(t)\|_\infty\le\frac{3}{2}\bar\rho,
\]
and
\[
 \sup_{\tau\le t\le T}\bigl(\|\nabla u(t)\|_2^2+\|\sqrt\rho\,\theta(t)\|_2^2
             +\sigma(t)\|\nabla\theta(t)\|_2^2\bigr)
 \le C_1E_0^{1/32}\le C_1\eps_0^{1/96}\le K.
\]
This, together with the estimates on \([0,\min\{T,\tau\}]\), completes the proof of
Proposition~\ref{prop:lower-order-closure}.
\end{proof}

Combining the estimates on $[0,\min\{T,\tau\}]$ and $[\tau,T]$ yields the
following corollary.
\begin{corollary}\label{cor:global-lower-estimates}
Under the assumptions of Proposition~\ref{prop:lower-order-closure}, we have
\begin{equation}\label{eq:global-rough-budget}
 \int_0^T\bigl(\|\sqrt\rho\,\dot u\|_2^2+\|\nabla\theta\|_2^2\bigr)\dd t\le C,
\end{equation}
where \(C=C(m_0,\bar\rho,M,\mu,\lambda,c_v,\kappa)\).
Moreover,
\begin{equation}\label{eq:global-first-time-weight}
\begin{aligned}
 &\sup_{0\le t\le T}\sigma(t)
     \bigl(\|\nabla u(t)\|_2^2+\|\sqrt\rho\,\theta(t)\|_2^2+\|P(t)\|_2^2\bigr)+\int_0^T\sigma(t)
     \bigl(\|\sqrt\rho\,\dot u\|_2^2+\|\nabla\theta\|_2^2\bigr)\dd t
 \le CE_0^{1/32},
\end{aligned}
\end{equation}
and
\begin{equation}\label{eq:global-second-time-weight}
\begin{aligned}
 &\sup_{0\le t\le T}\sigma^2(t)
     \bigl(\|\sqrt\rho\,\dot u(t)\|_2^2+\|\nabla\theta(t)\|_2^2\bigr)+\int_0^T\sigma^2(t)
     \bigl(\|\nabla\dot u\|_2^2+\|\sqrt\rho\,\dot\theta\|_2^2\bigr)\dd t
 \le CE_0^{1/32},
\end{aligned}
\end{equation}
where \(C=C(m_0,\bar\rho,M,\mu,\lambda,c_v,\kappa)\).
\end{corollary}
\begin{proof}
On \([0,\min\{T,\tau\}]\), we have \(\sigma(t)=t\) and \(E_0^{1/12}\le E_0^{1/32}\).
By \eqref{eq:initial-rough-lemma}--\eqref{eq:initial-rough-lemma-second}, we obtain
\[
 \int_0^{\min\{T,\tau\}}\bigl(\|\sqrt\rho\,\dot u\|_2^2+\|\nabla\theta\|_2^2\bigr)\dd t\le C,
\]
and by \eqref{eq:initial-weighted-lemma-second} and \eqref{eq:short-squared-weight}, we have
\[
 \sup_{0\le t\le\min\{T,\tau\}}\sigma(t)
     \bigl(\|\nabla u(t)\|_2^2+\|\sqrt\rho\,\theta(t)\|_2^2+\|P(t)\|_2^2\bigr)
 +\int_0^{\min\{T,\tau\}}\sigma(t)
     \bigl(\|\sqrt\rho\,\dot u\|_2^2+\|\nabla\theta\|_2^2\bigr)\dd t
 \le CE_0^{1/32},
\]
and
\[
 \sup_{0\le t\le\min\{T,\tau\}}\sigma^2(t)
     \bigl(\|\sqrt\rho\,\dot u(t)\|_2^2+\|\nabla\theta(t)\|_2^2\bigr)
 +\int_0^{\min\{T,\tau\}}\sigma^2(t)
     \bigl(\|\nabla\dot u\|_2^2+\|\sqrt\rho\,\dot\theta\|_2^2\bigr)\dd t
 \le CE_0^{1/32}.
\]
If \(T>\tau\), using \(\sigma(t)\le1\) and \(\|P\|_2^2\le2\bar\rho\|\sqrt\rho\,\theta\|_2^2\),
we deduce from \eqref{eq:lower-rough-improvement}--\eqref{eq:lower-rough-improvement-second} that
\[
\begin{aligned}
 &\sup_{\tau\le t\le T}\sigma(t)
     \bigl(\|\nabla u(t)\|_2^2+\|\sqrt\rho\,\theta(t)\|_2^2+\|P(t)\|_2^2\bigr)
 +\int_\tau^T\bigl(\|\sqrt\rho\,\dot u\|_2^2+\|\nabla\theta\|_2^2\bigr)\dd t
 \le CE_0^{1/32}.
\end{aligned}
\]
Furthermore, using \(\sigma^2(t)\le(1+t)^{3/2}\) and \(\log(e+t)(1+t)^{-9/8}\le C\),
we deduce from \eqref{eq:weighted-decay-lemma} and \eqref{eq:decay-thermal} that
\[
\begin{aligned}
 &\sup_{\tau\le t\le T}\sigma^2(t)
     \bigl(\|\sqrt\rho\,\dot u(t)\|_2^2+\|\nabla\theta(t)\|_2^2\bigr)
 +\int_\tau^T\sigma^2(t)
     \bigl(\|\nabla\dot u\|_2^2+\|\sqrt\rho\,\dot\theta\|_2^2\bigr)\dd t
 \le CE_0^{1/32}.
\end{aligned}
\]
Adding the estimates on \([0,\min\{T,\tau\}]\) and \([\tau,T]\) yields
\eqref{eq:global-rough-budget}--\eqref{eq:global-second-time-weight}.
\end{proof}

\section{Higher-order estimates}\label{sec:four}
In this section, let \((\rho\ge0,u,\theta\ge0)\) be a smooth solution of
\eqref{eq:perturbation-system} on \(\R^2\times[0,T]\) satisfying the assumptions of
Theorem~\ref{thm:target}. Throughout this section, \(C\) denotes a generic
positive constant depending only on \(m_0,\bar\rho,M,a,q,\eta_0,\mu,\lambda,c_v,\kappa,T\),
\(\|\bar x^a\rho_0\|_{H^1\cap W^{1,q}}\), \(\|\nabla u_0\|_{H^1}\), \(\|\nabla\theta_0\|_2\)
and \(\|g\|_2\).

We first establish the following estimates for \(\sqrt\rho\,\dot u\) and \(\nabla\theta\).
\begin{lemma}\label{lem:higher-acceleration-temperature}\label{prop:higher-estimates}
Under the assumptions of Theorem~\ref{thm:target}, we have
\begin{equation}\label{eq:higher-final-bounds}
\begin{aligned}
 &\sup_{0\le t\le T}
       \bigl(\|\sqrt\rho\dot u(t)\|_2^2+\|\nabla\theta(t)\|_2^2\bigr)
 +\int_0^T
       \bigl(\|\nabla\dot u\|_2^2+\|\sqrt\rho\dot\theta\|_2^2\bigr)\dd t
 \le C,
\end{aligned}
\end{equation}
and
\begin{equation}\label{eq:higher-final-bounds-second}
\begin{aligned}
 &\int_0^T\bigl(\|\nabla^2\theta\|_2^2+\|\theta\|_\infty^2\bigr)\dd t
 \le C,
\end{aligned}
\end{equation}
where \(C=C\bigl(m_0,\bar\rho,M,\mu,\lambda,c_v,\kappa,T,\|\nabla\theta_0\|_2,\|g\|_2\bigr)\) in both estimates.
\end{lemma}
\begin{proof}
Taking \(p=3,4\) in \eqref{eq:density-anchor} and using
\eqref{eq:physical-total-energy} and \eqref{eq:physical-first-moment}, we obtain
\[
 \|P\|_p\le(2\bar\rho)^{1-1/p}\|\rho^{1/p}\theta\|_p
 \le C(1+\|\nabla\theta\|_2).
\]
It follows from \eqref{eq:physical-total-energy},
\eqref{eq:initial-bootstrap-improvement} and
\eqref{eq:lift-source-moment}--\eqref{eq:lift-source-moment-second} that
\begin{equation}\label{eq:higher-source-weight}
\begin{aligned}
 &\bigl\||\nabla u|^2+P^2\bigr\|_1
       +\bigl\||\nabla u|^2+P^2\bigr\|_{3/2}
       +\int_{\R^2}(1+|x|^2)^{1/8}(|\nabla u|^2+P^2)\dd x\\
 \le& C\bigl(1+\|\sqrt\rho\,\dot u\|_2^2+\|\nabla\theta\|_2^2\bigr),
\end{aligned}
\end{equation}
and by \eqref{eq:quartic-velocity-pressure}, we also have
\begin{equation}\label{eq:higher-source-weight-second}
 \|\nabla u\|_4^4+\|P\nabla u\|_2^2+\|P\|_4^4
 \le C\bigl(1+\|\sqrt\rho\,\dot u\|_2^4+\|\nabla\theta\|_2^4\bigr).
\end{equation}

Taking \(w=|\nabla u|^2\) and \(v=\theta\) in \eqref{eq:lift-two-measure}, we deduce from \eqref{eq:physical-total-energy},
\eqref{eq:higher-source-weight} and the boundedness of \(s\log(e+1/s)\) on \((0,1]\) that
\begin{equation}\label{eq:higher-two-measure}
\begin{aligned}
 \int_{\R^2}\theta^2|\nabla u|^2\dd x
 \le& 2\int_{\R^2}|\theta-\langle\theta\rangle_\rho|^2|\nabla u|^2\dd x
       +2\langle\theta\rangle_\rho^2\|\nabla u\|_2^2\\
 \le& C\left[1+\|\nabla\theta\|_2^2
       \log\bigl(e+\|\sqrt\rho\,\dot u\|_2^2+
                         \|\nabla\theta\|_2^2\bigr)\right].
\end{aligned}
\end{equation}
By the Cauchy--Schwarz and Young inequalities, \eqref{eq:initial-bootstrap-improvement},
\eqref{eq:higher-two-measure} and \(\log(e+s)\le\delta s+C_\delta\) (\(\delta>0\)), we get,
for any \(\zeta>0\),
\begin{equation}\label{eq:higher-cross-coercivity}
\begin{aligned}
 0\le\int_{\R^2}\bigl(2\mu|D(u)|^2+\lambda(\Div u)^2\bigr)\theta\dd x
 \le& C\|\nabla u\|_2\left(\int_{\R^2}\theta^2|\nabla u|^2\dd x\right)^{1/2}\\
 \le& \zeta\bigl(\|\sqrt\rho\,\dot u\|_2^2+\|\nabla\theta\|_2^2\bigr)+C_\zeta.
\end{aligned}
\end{equation}
From \eqref{eq:perturbation-system}$_3$, \eqref{eq:initial-bootstrap-improvement}
and \eqref{eq:higher-source-weight-second}, we also obtain
\begin{equation}\label{eq:higher-thermal-hessian-first}
 \|\nabla^2\theta\|_2
 \le C\bigl(\|\sqrt\rho\,\dot\theta\|_2+1+
                \|\sqrt\rho\,\dot u\|_2^2+\|\nabla\theta\|_2^2\bigr).
\end{equation}

Substituting \eqref{eq:higher-source-weight-second} and \eqref{eq:higher-two-measure}
into \eqref{eq:heat-source-products} and applying \eqref{eq:higher-source-weight-second}
to \eqref{eq:material-acceleration-products}, we arrive at
\begin{equation}\label{eq:higher-thermal-cross-bound}
\begin{aligned}
 &\frac{\dd}{\dd t}\left(\frac\kappa2\|\nabla\theta\|_2^2
                         -\int_{\R^2} \bigl(2\mu|D(u)|^2+\lambda(\Div u)^2\bigr)\theta\dd x\right)
       +\frac{c_v}{2}\|\sqrt\rho\,\dot\theta\|_2^2\\
 \le&\eta\|\nabla\dot u\|_2^2
       +C_\eta\bigl(1+\|\sqrt\rho\,\dot u\|_2^4+\|\nabla\theta\|_2^4\bigr)+C_\eta\|\nabla\theta\|_2^2
       \log\bigl(e+\|\sqrt\rho\,\dot u\|_2^2+\|\nabla\theta\|_2^2\bigr),
\end{aligned}
\end{equation}
and
\begin{equation}\label{eq:higher-acceleration-bound}
\begin{aligned}
 &\frac{\dd}{\dd t}\|\sqrt\rho\,\dot u\|_2^2
       +\mu\|\nabla\dot u\|_2^2
 \le C_0\|\sqrt\rho\,\dot\theta\|_2^2
        +C\bigl(1+\|\sqrt\rho\,\dot u\|_2^4+\|\nabla\theta\|_2^4\bigr).
\end{aligned}
\end{equation}
Choose \(\eta=\frac{\mu}{2(1+4C_0/c_v)}\). Adding \(\bigl(1+\frac{4C_0}{c_v}\bigr)\) times
\eqref{eq:higher-thermal-cross-bound} to \eqref{eq:higher-acceleration-bound} and using
\(\log(e+s)\le1+s\) for \(s\ge0\), we obtain
\begin{equation}\label{eq:higher-combined}
\begin{aligned}
 &\frac{\dd}{\dd t}\Bigg[\|\sqrt\rho\,\dot u\|_2^2
       +\left(1+\frac{4C_0}{c_v}\right)
          \bigg(\frac\kappa2\|\nabla\theta\|_2^2
 -\int_{\R^2}\bigl(2\mu|D(u)|^2+\lambda(\Div u)^2\bigr)
                                      \theta\dd x\bigg)\Bigg]\\
 &+\frac\mu2\|\nabla\dot u\|_2^2
       +\frac{c_v+4C_0}{4}\|\sqrt\rho\,\dot\theta\|_2^2\\
 \le& C\bigl(1+\|\sqrt\rho\,\dot u\|_2^2+\|\nabla\theta\|_2^2\bigr)^2.
\end{aligned}
\end{equation}
Taking \(\zeta=\min\bigl\{\frac{1}{4(1+4C_0/c_v)},\frac{\kappa}{8}\bigr\}\) in
\eqref{eq:higher-cross-coercivity}, we set
\[
 \Phi(t)=\|\sqrt\rho\,\dot u(t)\|_2^2
 +\left(1+\frac{4C_0}{c_v}\right)
 \left[\frac{\kappa}{2}\|\nabla\theta(t)\|_2^2
       -\int_{\R^2}\bigl(2\mu|D(u)|^2+\lambda(\Div u)^2\bigr)\theta\dd x
       +C_\zeta\right]+1,
\]
and obtain
\begin{equation}\label{eq:higher-energy-equivalence}
 C^{-1}\bigl(1+\|\sqrt\rho\,\dot u(t)\|_2^2+\|\nabla\theta(t)\|_2^2\bigr)
 \le\Phi(t)
 \le C\bigl(1+\|\sqrt\rho\,\dot u(t)\|_2^2+\|\nabla\theta(t)\|_2^2\bigr).
\end{equation}
Integrating \eqref{eq:higher-combined} over \([0,t]\) and using
\eqref{eq:compatibility}, \eqref{eq:perturbation-system}$_2$ and \eqref{eq:higher-energy-equivalence}, we obtain
\begin{equation}\label{eq:higher-combined-integral}
\begin{aligned}
 &\|\sqrt\rho\,\dot u(t)\|_2^2+\|\nabla\theta(t)\|_2^2
       +\int_0^t\bigl(\|\nabla\dot u\|_2^2+\|\sqrt\rho\,\dot\theta\|_2^2\bigr)\dd r\\
 \le& C\bigl(1+\|g\|_2^2+\|\nabla\theta_0\|_2^2\bigr)
       +C\int_0^t\bigl(1+\|\sqrt\rho\,\dot u\|_2^2+\|\nabla\theta\|_2^2\bigr)^2\dd r.
\end{aligned}
\end{equation}
Applying Gronwall's inequality to \eqref{eq:higher-combined-integral} and using
\eqref{eq:global-rough-budget} yields \eqref{eq:higher-final-bounds}.

Next, integrating \eqref{eq:perturbation-system}$_3$ over \(\R^2\) and using \eqref{eq:perturbation-system}$_1$,
\eqref{eq:physical-thermal-energy} and \eqref{eq:viscous-dissipation-identity}, we arrive at
\begin{equation}\label{eq:thermal-laplacian-source}
 \int_{\R^2}\Delta\theta\dd x=0.
\end{equation}
By \eqref{eq:perturbation-system}$_3$, we have
\[
 |\Delta\theta|\le C\bigl(\rho|\dot\theta|+|\nabla u|^2+P^2\bigr),
\]
and by H\"older's inequality, \eqref{eq:physical-first-moment} and \eqref{eq:initial-bootstrap-improvement}, we get
\[
 \|\rho\dot\theta\|_{4/3}
 +\|(1+|x|^2)^{1/8}\rho\dot\theta\|_1
 \le C\|\sqrt\rho\,\dot\theta\|_2.
\]
Using H\"older's inequality, \eqref{eq:higher-source-weight}, \eqref{eq:higher-final-bounds}
and the two preceding inequalities, we conclude that
\begin{equation}\label{eq:thermal-laplacian-source-second}
 \|\Delta\theta\|_{4/3}
 +\int_{\R^2}(1+|x|^2)^{1/8}|\Delta\theta|\dd x
 \le C(1+\|\sqrt\rho\,\dot\theta\|_2).
\end{equation}
By \eqref{eq:higher-final-bounds}, \eqref{eq:thermal-laplacian-source}--\eqref{eq:thermal-laplacian-source-second}
and Lemma~\ref{lem:continuation-newton} with \(v=\theta\) and \(f=\Delta\theta\), there exists a constant
\(c=c(t)\) such that \(\theta-\Delta^{-1}\Delta\theta=c\). Multiplying this identity by \(\frac{\rho}{m_0}\)
and integrating over \(\R^2\), we deduce from \eqref{eq:mass} that
\begin{equation}\label{eq:continuation-newton}
 \theta=\Delta^{-1}\Delta\theta+\langle\theta\rangle_\rho
               -\frac{1}{m_0}\int_{\R^2}\rho\Delta^{-1}\Delta\theta\dd x.
\end{equation}
By \eqref{eq:mass}, \eqref{eq:physical-total-energy}, \eqref{eq:newton-potential-bound-second},
\eqref{eq:thermal-laplacian-source-second} and \eqref{eq:continuation-newton}, we arrive at
\begin{equation}\label{eq:continuation-newton-second}
 \|\theta\|_\infty
 \le 2\|\Delta^{-1}\Delta\theta\|_\infty+\langle\theta\rangle_\rho
 \le C(1+\|\sqrt\rho\,\dot\theta\|_2).
\end{equation}
Combining \eqref{eq:higher-thermal-hessian-first} and \eqref{eq:continuation-newton-second},
integrating over \([0,T]\) and using \eqref{eq:higher-final-bounds}, we obtain
\[
 \int_0^T\bigl(\|\nabla^2\theta\|_2^2+\|\theta\|_\infty^2\bigr)\dd t
 \le C\int_0^T\bigl(1+\|\sqrt\rho\,\dot\theta\|_2^2
      +\|\sqrt\rho\,\dot u\|_2^4+\|\nabla\theta\|_2^4\bigr)\dd t
 \le C,
\]
which gives \eqref{eq:higher-final-bounds-second}.
\end{proof}

Using Lemma~\ref{lem:higher-acceleration-temperature}, we now derive the
following estimates for $\rho$ and $\nabla u$.
\begin{lemma}\label{lem:higher-density}
Under the assumptions of Theorem~\ref{thm:target}, we have
\begin{equation}\label{eq:higher-density-bounds}
 \sup_{0\le t\le T}\|\rho(t)\|_{H^1\cap W^{1,q}}
 +\int_0^T\bigl(\|\nabla^2u\|_q^2+\|\nabla u\|_\infty^2\bigr)\dd t
 \le C,
\end{equation}
where \(C=C\bigl(m_0,\bar\rho,M,q,\mu,\lambda,c_v,\kappa,T,
\|\nabla\rho_0\|_{L^2\cap L^q},\|\nabla\theta_0\|_2,\|g\|_2\bigr)\).
\end{lemma}
\begin{proof}
Taking \(w=\dot u\) in \eqref{eq:lift-log-anchor} and using \eqref{eq:mean-material-acceleration},
\eqref{eq:flux-gradient-lp} and \eqref{eq:initial-bootstrap-improvement}, we obtain, for \(2\le p<\infty\),
\begin{equation}\label{eq:higher-flux-grad}
 \|\nabla F\|_p+\|\nabla\omega\|_p
 \le C\|\rho\dot u\|_p
 \le C\|\nabla\dot u\|_2.
\end{equation}
Substituting \eqref{eq:lift-lame-gradient-representation} into \eqref{eq:flux-definitions},
and applying \eqref{eq:gn-infty} with \(q=r=4\), \eqref{eq:riesz-potential-lp} with \(p=4/3\),
\eqref{eq:mass}, H\"older's inequality and \eqref{eq:higher-flux-grad}, we arrive at
\[
\begin{aligned}
 \|F\|_\infty+\|\omega\|_\infty
 \le& C\|I_1(|\rho\dot u|)\|_4^{1/2}\|\rho\dot u\|_4^{1/2}\\
 \le& C\|\rho\dot u\|_{4/3}^{1/2}\|\rho\dot u\|_4^{1/2}\\
 \le& C\|\nabla\dot u\|_2.
\end{aligned}
\]
Applying the div--curl estimate to \(\nabla u\), and using \eqref{eq:flux-definitions},
\eqref{eq:initial-bootstrap-improvement} and \eqref{eq:higher-flux-grad}, we have, for \(p=2,q\),
\begin{equation}\label{eq:higher-density-differential-second}
\begin{aligned}
 \|\nabla^2u\|_p
 \le& C\bigl(\|\nabla\Div u\|_p+\|\nabla\omega\|_p\bigr)\\
 \le& C\bigl(\|\nabla F\|_p+\|\nabla P\|_p+\|\nabla\omega\|_p\bigr)\\
 \le& C\bigl(\|\nabla\dot u\|_2+\|\theta\|_\infty\|\nabla\rho\|_p+\|\nabla\theta\|_p\bigr).
\end{aligned}
\end{equation}
Applying \(\nabla\) to \eqref{eq:perturbation-system}$_1$ gives
\begin{equation}\label{eq:higher-density-gradient-eq}
 (\partial_t+u\cdot\nabla)\nabla\rho
 +\bigl((\nabla u)^{\mathsf T}+(\Div u)I\bigr)\nabla\rho
 +\rho\nabla\Div u=0.
\end{equation}
Multiplying \eqref{eq:higher-density-gradient-eq} by \(p|\nabla\rho|^{p-2}\nabla\rho\)
(\(p=2,q\)) and integrating by parts, we deduce from H\"older's inequality,
\eqref{eq:initial-bootstrap-improvement} and \eqref{eq:higher-density-differential-second} that
\begin{equation}\label{eq:higher-density-differential}
\begin{aligned}
 \frac{\dd}{\dd t}\|\nabla\rho\|_p
 \le& C\|\nabla u\|_\infty\|\nabla\rho\|_p+C\|\nabla\Div u\|_p\\
 \le& C(\|\nabla u\|_\infty+\|\theta\|_\infty)\|\nabla\rho\|_p
       +C(\|\nabla\dot u\|_2+\|\nabla\theta\|_p).
\end{aligned}
\end{equation}
Combining \eqref{eq:flux-definitions}, \eqref{eq:higher-bkm}, \eqref{eq:continuation-newton-second},
\eqref{eq:higher-density-differential-second} and \eqref{eq:higher-density-differential}, we conclude that
\begin{equation}\label{eq:density-log-derivative}
\begin{aligned}
 \frac{\dd}{\dd t}
       \log(e+\|\nabla\rho\|_2+\|\nabla\rho\|_q)
 \le& C(1+\|\nabla\dot u\|_2+\|\sqrt\rho\,\dot\theta\|_2+\|\nabla\theta\|_q)
 \times\bigl[\log(e+\|\nabla\rho\|_2+\|\nabla\rho\|_q)\\
       &+\log(e+\|\nabla\dot u\|_2
                       +\|\sqrt\rho\,\dot\theta\|_2+\|\nabla\theta\|_q)\bigr].
\end{aligned}
\end{equation}
Applying \eqref{eq:gn-lp} to \(\nabla\theta\) and using
\eqref{eq:higher-final-bounds}--\eqref{eq:higher-final-bounds-second}, we get
\begin{equation}\label{eq:temperature-gradient-lq-time}
 \int_0^T\|\nabla\theta\|_q^2\dd t
 \le C\int_0^T\|\nabla\theta\|_2^{4/q}
                 \|\nabla^2\theta\|_2^{2-4/q}\dd t\le C.
\end{equation}
Combining \eqref{eq:higher-final-bounds} and \eqref{eq:temperature-gradient-lq-time} gives
\begin{equation}\label{eq:density-log-coefficient}
 \int_0^T
 \bigl(1+\|\nabla\dot u\|_2+\|\sqrt\rho\,\dot\theta\|_2
                         +\|\nabla\theta\|_q\bigr)^2\dd t\le C.
\end{equation}
Using \(r\log(e+r)\le C(1+r^2)\) in \eqref{eq:density-log-derivative} and applying Gronwall's inequality
with \eqref{eq:density-log-coefficient}, we obtain from \eqref{eq:higher-bkm} and
\eqref{eq:higher-density-differential-second}--\eqref{eq:higher-density-differential} that
\[
 \sup_{0\le t\le T}
       \bigl(\|\nabla\rho\|_2+\|\nabla\rho\|_q\bigr)
 +\int_0^T\bigl(\|\nabla u\|_\infty+\|\nabla^2u\|_q^2\bigr)\dd t
 \le C.
\]
Since \(q>2\), we conclude from \eqref{eq:gn-infty} that
\[
 \int_0^T\|\nabla u\|_\infty^2\dd t
 \le C\int_0^T\|\nabla u\|_2^{(q-2)/(q-1)}
                  \|\nabla^2u\|_q^{q/(q-1)}\dd t\le C.
\]
Finally, by \eqref{eq:mass} and \eqref{eq:initial-bootstrap-improvement}, we have
\[
 \sup_{0\le t\le T}\bigl(\|\rho(t)\|_2+\|\rho(t)\|_q\bigr)\le C,
\]
which completes the proof of \eqref{eq:higher-density-bounds}.
\end{proof}

The following lemma gives the spatially weighted estimates for $\rho$.
\begin{lemma}\label{lem:higher-weighted-density}
Under the assumptions of Theorem~\ref{thm:target}, we have
\begin{equation}\label{eq:higher-weighted-density-bounds}
 \sup_{0\le t\le T}\|\bar x^a\rho(t)\|_{L^1\cap H^1\cap W^{1,q}}
 \le C,
\end{equation}
where \(C=C\bigl(m_0,\bar\rho,M,a,q,\eta_0,\mu,\lambda,c_v,\kappa,T,
\|\bar x^a\rho_0\|_{H^1\cap W^{1,q}},\|\nabla\theta_0\|_2,\|g\|_2\bigr)\).
\end{lemma}
\begin{proof}
By \eqref{eq:gn-infty}, \eqref{eq:spatial-weighted-lp}, \eqref{eq:physical-first-moment},
\eqref{eq:higher-source-weight-second} and \eqref{eq:higher-final-bounds}, we arrive at
\[
 \left\|\frac{u}{(1+|x|^2)^{1/2}}\right\|_\infty\le C.
\]
Multiplying \eqref{eq:perturbation-system}$_1$ by \(\bar x^a\) gives
\begin{equation}\label{eq:weighted-density-transport}
 (\partial_t+u\cdot\nabla)(\bar x^a\rho)
       =(au\cdot\nabla\log\bar x-\Div u)\bar x^a\rho.
\end{equation}
By \eqref{eq:higher-density-bounds}, \eqref{eq:weighted-density-transport} and Gronwall's inequality, we have
\[
 \sup_{0\le t\le T}\|\bar x^a\rho(t)\|_{L^1\cap L^\infty}\le C.
\]
Differentiating \eqref{eq:weighted-density-transport} in space yields
\begin{equation}\label{eq:higher-weighted-transport}
\begin{aligned}
 (\partial_t+u\cdot\nabla)\nabla(\bar x^a\rho)
 =&\bigl[(au\cdot\nabla\log\bar x-\Div u)I
                         -(\nabla u)^{\mathsf T}\bigr]\nabla(\bar x^a\rho)\\
 &+a\bar x^a\rho\bigl[(\nabla u)^{\mathsf T}\nabla\log\bar x
                              +(\nabla^2\log\bar x)u\bigr]
       -\bar x^a\rho\nabla\Div u,
\end{aligned}
\end{equation}
where, by \eqref{eq:flux-definitions},
\begin{equation}\label{eq:higher-weighted-transport-second}
\begin{aligned}
 -\bar x^a\rho\nabla\Div u
 =&-\frac1{2\mu+\lambda}\bigl[\bar x^a\rho\nabla F
                                +P\nabla(\bar x^a\rho)
                                -a\bar x^a\rho P\nabla\log\bar x
                                +\bar x^a\rho^2\nabla\theta\bigr].
\end{aligned}
\end{equation}
Taking \(L^p\) norms in \eqref{eq:higher-weighted-transport}--\eqref{eq:higher-weighted-transport-second} for \(p=2,q\)
and using \eqref{eq:flux-gradient-lp}, \eqref{eq:lift-log-anchor} with \(w=\dot u\)
and \eqref{eq:mean-material-acceleration}, we conclude that
\begin{equation}\label{eq:weighted-density-gradient}
\begin{aligned}
 \frac{\dd}{\dd t}\|\nabla(\bar x^a\rho)\|_p
 \le& C(1+\|\nabla u\|_\infty+\|\theta\|_\infty)
                                  \|\nabla(\bar x^a\rho)\|_p\\
 &+C(1+\|\nabla\dot u\|_2+\|\nabla\theta\|_p+\|\theta\|_\infty).
\end{aligned}
\end{equation}
By \eqref{eq:higher-final-bounds}--\eqref{eq:higher-final-bounds-second}, \eqref{eq:higher-density-bounds}
and \eqref{eq:temperature-gradient-lq-time}, we have, for \(p=2,q\),
\begin{equation}\label{eq:weighted-density-coefficient}
 \int_0^T
 \bigl(1+\|\nabla u\|_\infty+\|\theta\|_\infty
           +\|\nabla\dot u\|_2+\|\nabla\theta\|_p\bigr)\dd t\le C.
\end{equation}
Applying Gronwall's inequality to \eqref{eq:weighted-density-gradient} with
\eqref{eq:weighted-density-coefficient} yields \eqref{eq:higher-weighted-density-bounds}.
\end{proof}

From Lemma~\ref{lem:higher-weighted-density}, we obtain the higher-order
estimates for $u$ and its far-field limit.
\begin{lemma}\label{lem:higher-velocity}
Under the assumptions of Theorem~\ref{thm:target}, we have
\begin{equation}\label{eq:higher-velocity-bounds}
\begin{aligned}
 &\sup_{0\le t\le T}
       \bigl(\|\nabla^2u\|_2^2+\|u\|_\infty^2
               +\|\sqrt\rho u_t\|_2^2
               +\|\bar x^{\alpha/2}\nabla u\|_2^2\bigr)+\int_0^T\|\nabla u_t\|_2^2\dd t\le C,
\end{aligned}
\end{equation}
where \(C=C\bigl(m_0,\bar\rho,M,a,q,\eta_0,\mu,\lambda,c_v,\kappa,T,
\|\bar x^a\rho_0\|_{H^1\cap W^{1,q}},\|\nabla u_0\|_{H^1},\|\nabla\theta_0\|_2,\|g\|_2\bigr)\).
For $0\le t\le T$, the velocity has the spatial limit
\begin{equation}\label{eq:moving-velocity-limit}
 \lim_{|x|\to\infty}u(x,t)
 =\langle u\rangle_\rho+\frac1{m_0}
       \int_{\R^2}(\nabla\Delta^{-1}\rho\cdot\nabla)u\dd x .
\end{equation}
\end{lemma}
\begin{proof}
By \eqref{eq:spatial-weighted-lp} and \eqref{eq:higher-weighted-density-bounds}, we have
\[
 \|(1+|x|^2)^{(q-1)/(2q)}\nabla\rho\|_q\le C,\qquad
 \|(1+|x|^2)^{-(q-1)/(2q)}\theta\|_{2q/(q-2)}\le C.
\]
Using H\"older's inequality, we arrive at
\begin{equation}\label{eq:pressure-density-gradient}
 \|\theta\nabla\rho\|_2
 \le\|(1+|x|^2)^{(q-1)/(2q)}\nabla\rho\|_q
       \|(1+|x|^2)^{-(q-1)/(2q)}\theta\|_{2q/(q-2)}\le C.
\end{equation}
For any fixed \(t\in[0,T]\), applying Plancherel's identity to \eqref{eq:perturbation-system}$_2$
and using \eqref{eq:physical}, \eqref{eq:initial-bootstrap-improvement}, \eqref{eq:higher-final-bounds}
and \eqref{eq:pressure-density-gradient}, we obtain
\begin{equation}\label{eq:velocity-gradient-h1}
\begin{aligned}
 &\|\nabla u(t)\|_{H^1}^2
 \le K+\frac{1}{\mu^2}\bigl(\|\rho\dot u\|_2+\|\theta\nabla\rho\|_2
                              +\|\rho\nabla\theta\|_2\bigr)^2\le C.
\end{aligned}
\end{equation}
Applying \eqref{eq:gn-lp} with \(p=16\) to \(\nabla u(t)\) and using
\eqref{eq:velocity-gradient-h1}, we have
\[
 \|\nabla u(t)\|_{16}
 \le C\|\nabla u(t)\|_2^{1/8}\|\nabla^2u(t)\|_2^{7/8}\le C.
\]
Moreover, it follows from \eqref{eq:higher-final-bounds} and \eqref{eq:higher-source-weight} that
\[
 \||x|^{1/8}\nabla u(t)\|_2^2
 \le\int_{\R^2}(1+|x|^2)^{1/8}|\nabla u(t)|^2\dd x\le C.
\]
Let \(0\le\eta_\varepsilon\in C_c^\infty(B_\varepsilon)\) be a standard mollifier with
\(\int_{\R^2}\eta_\varepsilon\dd x=1\). Using \(|x|^{1/8}\le|x-y|^{1/8}+|y|^{1/8}\)
and Minkowski's inequality, we get
\[
 \||x|^{1/8}\nabla(\eta_\varepsilon*u)\|_2
 \le\||x|^{1/8}\nabla u\|_2+\varepsilon^{1/8}\|\nabla u\|_2.
\]
By \cite[Lemma~2.1]{FanLuWang2026}, there exists \(c_\varepsilon\in\R^2\) such that
\[
 \|\eta_\varepsilon*u-c_\varepsilon\|_{16}
 \le C\||x|^{1/8}\nabla(\eta_\varepsilon*u)\|_2
 \le C\bigl(\||x|^{1/8}\nabla u\|_2+\varepsilon^{1/8}\|\nabla u\|_2\bigr).
\]
By H\"older's inequality on \(B_1\), we deduce that
\[
\begin{aligned}
 |c_\varepsilon|
 =&\frac{1}{\pi}\int_{B_1}|c_\varepsilon|\dd x\\
 \le&\frac{1}{\pi}\int_{B_1}\bigl(|\eta_\varepsilon*u|+|\eta_\varepsilon*u-c_\varepsilon|\bigr)\dd x\\
 \le&\frac{1}{\pi}\|\eta_\varepsilon*u\|_{L^1(B_1)}
      +\pi^{-1/16}\|\eta_\varepsilon*u-c_\varepsilon\|_{16}.
\end{aligned}
\]
Since \(\eta_\varepsilon*u\to u\) in \(H^1_{\rm loc}\) as \(\varepsilon\to0^+\),
 \(\{c_\varepsilon\}\) is bounded, which implies that
\(c_\varepsilon\to c\in\R^2\), and \(\eta_\varepsilon*u\to u\) almost everywhere in \(\R^2\).
Passing to the limit by Fatou's lemma, we arrive at
\begin{equation}\label{eq:weighted-sobolev-16}
 \|u-c\|_{16}\le C\||x|^{1/8}\nabla u\|_2.
\end{equation}
By \eqref{eq:weighted-sobolev-16} and \(\nabla u(t)\in L^{16}\), we have
\(u(t)-c\in W^{1,16}(\R^2)\hookrightarrow C_0(\R^2)\),
so \(c=\lim_{|x|\to\infty}u(x,t)\) and
\begin{equation}\label{eq:velocity-spatial-remainder}
 \left\|u(t)-\lim_{|y|\to\infty}u(y,t)\right\|_{W^{1,16}}
 \le C\bigl(\||x|^{1/8}\nabla u(t)\|_2+\|\nabla u(t)\|_{16}\bigr)\le C.
\end{equation}
Multiplying \(u(x,t)=u(y,t)+(u(x,t)-c)-(u(y,t)-c)\) by \(\frac{\rho(y,t)}{m_0}\) and
integrating over \(\R^2\) with respect to \(y\), we deduce from \eqref{eq:mass},
\eqref{eq:physical-total-energy}, \eqref{eq:velocity-spatial-remainder} and the Cauchy--Schwarz
inequality that
\[
 \|u(t)\|_\infty\le m_0^{-1/2}\|\sqrt\rho\,u(t)\|_2
     +2\left\|u(t)-\lim_{|y|\to\infty}u(y,t)\right\|_\infty\le C.
\]
Applying
Lemmas~\ref{lem:higher-acceleration-temperature} and \ref{lem:higher-density}, we obtain
\begin{equation}\label{eq:higher-velocity-time}
 \sup_{0\le t\le T}\|\sqrt\rho\,u_t\|_2\le C,\qquad
 \int_{0}^{T}\|\nabla u_t\|_2^2\dd t\le C.
\end{equation}

Next, for \(\alpha/2<\beta<\gamma<\min\{a/2,1\}\), we choose
\[
 1<p<\min\left\{\frac2{2-\gamma+\beta},
                         \frac{2a}{a+2\gamma}\right\},
\]
such that
\[
 |x|^{\beta-\gamma}\in L^{p/(p-1)}(\{|x|>1\}),\qquad
 \frac{2\gamma p}{2-p}<a.
\]
By H\"older's inequality and \eqref{eq:higher-weighted-density-bounds}, we have
\begin{equation}\label{eq:weighted-acceleration-lp}
 \||x|^\gamma\rho\dot u\|_p
 \le\|\sqrt\rho\,\dot u\|_2
 \left(\int_{\R^2} |x|^{2\gamma p/(2-p)}\rho^{p/(2-p)}\dd x\right)^{(2-p)/(2p)}
 \le C\|\sqrt\rho\,\dot u\|_2.
\end{equation}
Fundamental matrix \(\Gamma\) is defined by \eqref{eq:lift-lame-fundamental-solution}. For \(|x|>2|y|\), the mean value theorem
applied to \(\partial_i\Gamma\) gives
\[
 |(\partial_i\Gamma)(x-y)-(\partial_i\Gamma)(x)|\le C|y|^\gamma|x|^{-1-\gamma}.
\]
By \eqref{eq:mean-material-acceleration}, we get, for \(|x|>1\),
\[
 |( (\partial_i\Gamma)*(\rho\dot u) )(x)|
 \le C|x|^{-1-\gamma}\||x|^\gamma\rho\dot u\|_1
     +C|x|^{-\gamma} I_1(|\cdot|^\gamma|\rho\dot u|)(x).
\]
Applying \eqref{eq:riesz-potential-lp} with \(f=|\rho\dot u|\) on \(B_1\) and
\(f=|\cdot|^\gamma|\rho\dot u|\) on \(\R^2\setminus B_1\), and using H\"older's inequality
and \eqref{eq:weighted-riesz-pressure} in \eqref{eq:lift-lame-gradient-representation}, we arrive at
\begin{equation}\label{eq:weighted-lame-gradient}
\begin{aligned}
 \|(1+|x|^2)^{\beta/2}\nabla u\|_2
 \le& C\bigl(\|\rho\dot u\|_p+\||x|^\gamma\rho\dot u\|_p\bigr)\\
 &+C\left(\int_{\R^2}(1+|x|^2)^\gamma\rho\dd x\right)^{1/2}
                         \|\sqrt\rho\,\dot u\|_2
       +C\|(1+|x|^2)^{\beta/2}P\|_2.
\end{aligned}
\end{equation}
By \eqref{eq:higher-weighted-density-bounds} and \(2\gamma<a\), we have
\begin{equation}\label{eq:weighted-density-pressure}
 \int_{\R^2}(1+|x|^2)^{\gamma}\rho\dd x\le C,
\end{equation}
and
\begin{equation}\label{eq:weighted-density-pressure-second}
 \int_{\R^2}(1+|x|^2)^{\beta}P^2\dd x
 \le\|(1+|x|^2)^{\beta}\rho\|_\infty
           \|\sqrt\rho\,\theta\|_2^2\le C.
\end{equation}
Substituting \eqref{eq:weighted-acceleration-lp} and
\eqref{eq:weighted-density-pressure}--\eqref{eq:weighted-density-pressure-second} into
\eqref{eq:weighted-lame-gradient}, and using \(\bar x^\alpha\le C(1+|x|^2)^\beta\) and
\eqref{eq:higher-final-bounds}, we conclude that
\[
 \sup_{0\le t\le T}\|\bar x^{\alpha/2}\nabla u(t)\|_2\le C.
\]

Finally, since \(\alpha>1/2\), H\"older's inequality gives
\[
 \|\nabla u\|_{16/9}
 \le\|\bar x^{-\alpha/2}\|_{16}\|\bar x^{\alpha/2}\nabla u\|_2<\infty,
\]
and \eqref{eq:riesz-potential-lp} with \(p=16/15\) gives
\[
 \|\nabla\Delta^{-1}\rho\|_{16/7}\le C\|\rho\|_{16/15}<\infty.
\]
Hence, by H\"older's inequality, we have
\[
 \int_{\R^2}|(\nabla\Delta^{-1}\rho\cdot\nabla)u|\dd x
 \le\|\nabla\Delta^{-1}\rho\|_{16/7}\|\nabla u\|_{16/9}<\infty,
\]
and
\[
 \int_{\R^2}\rho\left|u-\lim_{|y|\to\infty}u(y,t)\right|\dd x
 \le\|\rho\|_{16/15}
      \left\|u-\lim_{|y|\to\infty}u(y,t)\right\|_{16}<\infty.
\]
Multiplying \(\Div(\nabla\Delta^{-1}\rho)=\rho\) by
\(u-\lim_{|y|\to\infty}u(y,t)\) and integrating by parts, we deduce from
\eqref{eq:mass} that
\begin{equation}\label{eq:velocity-mean-limit-identity}
\begin{aligned}
 \int_{\R^2}(\nabla\Delta^{-1}\rho\cdot\nabla)u\dd x
 =&-\int_{\R^2}\rho\left(u-\lim_{|y|\to\infty}u(y,t)\right)\dd x\\
 =&m_0\left(\lim_{|y|\to\infty}u(y,t)-\langle u\rangle_\rho\right),
\end{aligned}
\end{equation}
which proves \eqref{eq:moving-velocity-limit}.
\end{proof}

Finally, we establish the higher-order time-weighted estimates for $\theta$.
\begin{lemma}\label{lem:higher-thermal-estimates}
Under the assumptions of Theorem~\ref{thm:target}, we have
\begin{equation}\label{eq:higher-temperature-bounds}
\begin{aligned}
 &\sup_{0\le t\le T}
       \bigl(\|\sqrt{t\rho}\,\dot\theta(t)\|_2^2
                +\|\sqrt t\,\nabla^2\theta(t)\|_2^2+\|\sqrt t\,\theta(t)\|_\infty^2\bigr)\\
 &\qquad+\int_0^Tt
       \bigl(\|\nabla\dot\theta\|_2^2+\|\nabla\theta_t\|_2^2
                                 +\|\nabla^3\theta\|_2^2\bigr)\dd t\le C,
\end{aligned}
\end{equation}
where \(C=C\bigl(m_0,\bar\rho,M,a,q,\eta_0,\mu,\lambda,c_v,\kappa,T,
\|\bar x^a\rho_0\|_{H^1\cap W^{1,q}},\|\nabla u_0\|_{H^1},\|\nabla\theta_0\|_2,\|g\|_2\bigr)\).
\end{lemma}
\begin{proof}
Applying \(\partial_t+\Div(u\,\cdot)\) to \eqref{eq:perturbation-system}$_3$ and using \eqref{eq:perturbation-system}$_1$, we obtain
\begin{equation}\label{eq:higher-material-heat}
\begin{aligned}
 &c_v\rho (\partial_t+u\cdot\nabla)\dot\theta-\kappa\Delta\dot\theta\\
 =&\kappa\bigl[(\Div u)\Delta\theta
       -\partial_i(\partial_iu\cdot\nabla\theta)
       -\partial_iu\cdot\nabla\partial_i\theta\bigr]\\
 &+(\partial_t+u\cdot\nabla+\Div u)\bigl(2\mu|D(u)|^2+\lambda(\Div u)^2\bigr)
       -\rho\dot\theta\Div u\\
 &-P\Div\dot u+P\partial_iu_j\partial_ju_i.
\end{aligned}
\end{equation}

By \eqref{eq:regularity}, \eqref{eq:higher-final-bounds} and \eqref{eq:higher-velocity-bounds},
we have \(\dot\theta-\langle\dot\theta\rangle_\rho\in\mathcal H_\rho\) for almost every \(t\in(0,T)\). Applying \eqref{eq:lift-two-measure} with \(w=|\nabla u|^2+P^2\) and \(v=\dot\theta\),
and using \eqref{eq:higher-source-weight}, \eqref{eq:higher-final-bounds} and
\(m_0|\langle\dot\theta\rangle_\rho|\le\sqrt{m_0}\|\sqrt\rho\,\dot\theta\|_2\), we arrive at
\begin{equation}\label{eq:higher-material-products}
 \|\dot\theta\nabla u\|_2+\|P\dot\theta\|_2
 \le C\bigl(\|\sqrt\rho\,\dot\theta\|_2+\|\nabla\dot\theta\|_2\bigr).
\end{equation}
Multiplying \eqref{eq:higher-material-heat} by \(\dot\theta\), integrating by parts
and using H\"older's inequality together with \eqref{eq:gn-lp} for \(\nabla\theta\), we deduce that
\begin{equation}\label{eq:thermal-mixed-products}
\begin{aligned}
 &\frac{c_v}{2}\frac{\dd}{\dd t}\|\sqrt\rho\,\dot\theta\|_2^2+\kappa\|\nabla\dot\theta\|_2^2\\
 \le& C\|\Div u\|_\infty\|\sqrt\rho\,\dot\theta\|_2^2
      +C\|\nabla u\|_4\|\nabla\theta\|_2^{1/2}\|\nabla^2\theta\|_2^{1/2}\|\nabla\dot\theta\|_2+C\|\dot\theta\nabla u\|_2\bigl(\|\nabla^2\theta\|_2\\
      &+\|\nabla\dot u\|_2+\|\nabla u\|_4^2\bigr)
      +C\|P\dot\theta\|_2\bigl(\|\nabla\dot u\|_2+\|\nabla u\|_4^2\bigr).
\end{aligned}
\end{equation}
Since \(\|\nabla u\|_4\le C\) by \eqref{eq:higher-source-weight-second} and \eqref{eq:higher-final-bounds},
applying \eqref{eq:higher-material-products} and Young's inequality in \eqref{eq:thermal-mixed-products} yields
\begin{equation}\label{eq:higher-thermal-material-budget}
\begin{aligned}
 \frac{c_v}{2}\frac{\dd}{\dd t}\|\sqrt\rho\,\dot\theta\|_2^2
       +\frac\kappa2\|\nabla\dot\theta\|_2^2
 \le& C(1+\|\Div u\|_\infty)\|\sqrt\rho\,\dot\theta\|_2^2\\
 &+C(1+\|\nabla\dot u\|_2^2+\|\nabla^2\theta\|_2^2).
\end{aligned}
\end{equation}
Multiplying \eqref{eq:higher-thermal-material-budget} by \(s\), integrating over \([0,t]\)
and applying Gronwall's inequality, we deduce from \eqref{eq:higher-final-bounds}--\eqref{eq:higher-final-bounds-second}
and \eqref{eq:higher-density-bounds} that
\begin{equation}\label{eq:weighted-temperature-derivative}
\begin{aligned}
 \|\sqrt{t\rho}\,\dot\theta(t)\|_2^2+\int_0^t s\|\nabla\dot\theta(s)\|_2^2\dd s
 \le& C\int_0^t\|\sqrt\rho\,\dot\theta(s)\|_2^2\dd s
       +C\int_0^t s\bigl(1+\|\nabla\dot u(s)\|_2^2+\|\nabla^2\theta(s)\|_2^2\bigr)\dd s\\
       \le& C.
\end{aligned}
\end{equation}
By \eqref{eq:higher-thermal-hessian-first}, \eqref{eq:higher-final-bounds},
\eqref{eq:continuation-newton-second} and \eqref{eq:weighted-temperature-derivative}, we have
\[
 \|\sqrt t\,\nabla^2\theta(t)\|_2^2+\|\sqrt t\,\theta(t)\|_\infty^2
 \le C\bigl(t+\|\sqrt{t\rho}\,\dot\theta(t)\|_2^2\bigr)\le C.
\]
Noting that \(\nabla\theta_t=\nabla\dot\theta-(\nabla u)^{\mathsf T}\nabla\theta-u\cdot\nabla^2\theta\)
and using \eqref{eq:higher-final-bounds}--\eqref{eq:higher-final-bounds-second},
\eqref{eq:higher-density-bounds}, \eqref{eq:higher-velocity-bounds} and
\eqref{eq:weighted-temperature-derivative}, we conclude that
\[
 \int_0^Tt\|\nabla\theta_t\|_2^2\dd t
 \le 3\int_0^Tt\bigl(\|\nabla\dot\theta\|_2^2+
       \|\nabla u\|_\infty^2\|\nabla\theta\|_2^2+
       \|u\|_\infty^2\|\nabla^2\theta\|_2^2\bigr)\dd t
 \le C.
\]

Finally, differentiating \eqref{eq:perturbation-system}$_3$ in space gives
\begin{equation}\label{eq:differentiated-heat-equation}
 \kappa\nabla\Delta\theta
 =c_v\nabla\rho\,\dot\theta+c_v\rho\nabla\dot\theta
       -\nabla \bigl(2\mu|D(u)|^2+\lambda(\Div u)^2\bigr)+\nabla P\Div u+P\nabla\Div u.
\end{equation}
By \eqref{eq:spatial-weighted-lp}, \eqref{eq:higher-weighted-density-bounds} and H\"older's inequality, we have
\begin{equation}\label{eq:thermal-density-derivative}
\begin{aligned}
 \|\dot\theta\nabla\rho\|_2
 \le&\|(1+|x|^2)^{(q-1)/(2q)}\nabla\rho\|_q
       \|(1+|x|^2)^{-(q-1)/(2q)}\dot\theta\|_{2q/(q-2)}\\
 \le& C\bigl(\|\sqrt\rho\,\dot\theta\|_2+\|\nabla\dot\theta\|_2\bigr),
\end{aligned}
\end{equation}
and by \eqref{eq:higher-final-bounds} and \eqref{eq:pressure-density-gradient}, we also have
\[
 \|\nabla P\|_2\le\|\theta\nabla\rho\|_2+\|\rho\nabla\theta\|_2\le C.
\]
Applying the \(L^2\) elliptic estimate to \eqref{eq:differentiated-heat-equation} and using
\eqref{eq:higher-final-bounds}--\eqref{eq:higher-final-bounds-second}, \eqref{eq:higher-density-bounds},
\eqref{eq:higher-velocity-bounds}, \eqref{eq:weighted-temperature-derivative} and
\eqref{eq:thermal-density-derivative}, we arrive at
\[
 \int_0^Tt\|\nabla^3\theta\|_2^2\dd t
 \le C\int_0^Tt\bigl(1+\|\sqrt\rho\,\dot\theta\|_2^2+
       \|\nabla\dot\theta\|_2^2+\|\nabla u\|_\infty^2+
       \|\theta\|_\infty^2\bigr)\dd t\le C,
\]
which completes the proof of \eqref{eq:higher-temperature-bounds}.
\end{proof}

\section{Proofs of Theorems~\ref{thm:target} and~\ref{thm:decay}}\label{sec:five}\label{subsec:decay}
With all the a priori estimates established in Sections~\ref{sec:three} and~\ref{sec:four},
we are now ready to prove Theorems~\ref{thm:target} and~\ref{thm:decay}.

\begin{proof}[Proof of Theorem~\ref{thm:target}]
Using the Galilean transformation \eqref{eq:moving-frame}, we first consider
\eqref{eq:perturbation-system}. By Remark~\ref{rem:local-existence}, there exists
\(T_0>0\) such that \eqref{eq:perturbation-system} with initial data \eqref{eq:initial}
admits a local strong solution \((\rho,u,\theta)\) on \(\R^2\times[0,T_0]\)
satisfying \eqref{eq:regularity}. The continuation argument assumes local uniqueness
in this class as specified in Remark~\ref{rem:local-existence}.

 By \eqref{eq:initial-data-bound}, \eqref{eq:bootstrap-K}
and \(\sigma(0)=0\), we have
\[
 \|\rho_0\|_\infty\le\bar\rho<2\bar\rho,\qquad
 \|\nabla u_0\|_2^2+\|\sqrt{\rho_0}\,\theta_0\|_2^2
 +\sigma(0)\|\nabla\theta_0\|_2^2\le M^2<K<2K.
\]
By the regularity \eqref{eq:regularity} of the local strong solution, there exists
\(T_1\in(0,T_0]\) such that \eqref{eq:initial-bootstrap-assumption} holds on \([0,T_1]\).
Set
\begin{equation}\label{eq:bootstrap-maximal-time}
 T_*=\sup\bigl\{T>0:\ (\rho,u,\theta)\text{ exists on }\R^2\times[0,T]
               \text{ and satisfies }\eqref{eq:initial-bootstrap-assumption}\bigr\}.
\end{equation}
Clearly, \(T_*\ge T_1>0\), and Proposition~\ref{prop:lower-order-closure} implies that
\eqref{eq:initial-bootstrap-improvement} holds on \([0,T]\) for every \(0<T<T_*\).

Suppose that \(T_*<\infty\). By \eqref{eq:physical-total-energy},
\eqref{eq:initial-bootstrap-improvement}, \eqref{eq:higher-final-bounds},
\eqref{eq:higher-weighted-density-bounds} and \eqref{eq:higher-velocity-bounds}, we have
\begin{equation}\label{eq:restart-data-bound}
\begin{aligned}
 \operatorname*{ess\,sup}_{0<t<T_*}\bigl(&\|\bar x^a\rho(t)\|_{L^1\cap H^1\cap W^{1,q}}
       +\|\sqrt\rho\,u(t)\|_2+\|\sqrt\rho\,\theta(t)\|_2\\
       &+\|\nabla u(t)\|_{H^1}+\|\nabla\theta(t)\|_2+\|\sqrt\rho\,\dot u(t)\|_2\bigr)\le C,
\end{aligned}
\end{equation}
where \(C\) is uniform for \(t<T_*\). For almost every \(s\in(0,T_*)\),
\eqref{eq:perturbation-system}$_2$ gives
\begin{equation}\label{eq:restart-compatibility}
 -\mu\Delta u(s)-(\mu+\lambda)\nabla\Div u(s)+\nabla(\rho\theta)(s)
       =-\sqrt{\rho(s)}\,\bigl(\sqrt{\rho(s)}\,\dot u(s)\bigr).
\end{equation}
By \eqref{eq:mass} and \eqref{eq:physical-first-moment}, taking \(R>0\) sufficiently large
and independent of \(s\), we have
\[
 \int_{B_R}\rho(s)\dd x
 =m_0-\int_{|x|\ge R}\rho(s)\dd x
 \ge m_0-\frac{1}{R}\int_{\R^2}|x|\rho(s)\dd x
 \ge\frac{m_0}{2}.
\]
By \eqref{eq:restart-data-bound}--\eqref{eq:restart-compatibility} and
Remark~\ref{rem:local-existence}, there exists a uniform time \(T_0>0\) such that the
Cauchy problem with initial data \((\rho,u,\theta)(\cdot,s)\) admits a unique local strong
solution on \(\R^2\times[s,s+T_0]\). Taking \(s\in(0,T_*)\) sufficiently close to \(T_*\)
so that \(s+T_0>T_*\), we deduce from local uniqueness that the solution extends to
\([0,s+T_0]\) and satisfies \eqref{eq:regularity}. Together with \eqref{eq:initial-bootstrap-improvement} and continuity,
this implies that \eqref{eq:initial-bootstrap-assumption} holds on \([0,T]\) for some
\(T>T_*\), contradicting \eqref{eq:bootstrap-maximal-time}. Hence \(T_*=\infty\).
This completes the proof of Theorem~\ref{thm:target}.
\end{proof}

\begin{proof}[Proof of Theorem~\ref{thm:decay}]
Consider first the transformed solution under the Galilean transformation \eqref{eq:moving-frame}.
For all \(t\ge2>\tau\), it follows directly from \eqref{eq:weighted-decay-lemma},
\eqref{eq:decay-mechanical-second} and \eqref{eq:decay-thermal} that
\[
 \|\nabla u(t)\|_2+\|P(t)\|_2\le C(1+t)^{-3/8},
 \qquad \|P(t)\|_3\le C(1+t)^{-1/2},
\]
and
\[
 \|\nabla\theta(t)\|_2\le C\log^{1/2}(e+t)(1+t)^{-9/16}.
\]

By \eqref{eq:physical-total-energy}, we have \(\langle u(t)\rangle_\rho=u_*-u_\infty\).
Applying \eqref{eq:lift-log-anchor} and \eqref{eq:decay-mechanical-second}, we obtain
\begin{equation}\label{eq:velocity-variance-decay}
 \|\sqrt\rho(u-\langle u\rangle_\rho)(t)\|_2^2
 \le C\log(e+t)\|\nabla u(t)\|_2^2
 \le C\log(e+t)(1+t)^{-3/4}\longrightarrow0
 \qquad (t\to\infty).
\end{equation}
Using \(\langle u(t)\rangle_\rho=u_*-u_\infty\) and \eqref{eq:mass}, we get
\begin{equation}\label{eq:conserved-limits}
 \frac{1}{2}\int_{\R^2}\rho|u|^2\dd x
 =\frac{m_0}{2}|u_*-u_\infty|^2
       +\frac{1}{2}\|\sqrt\rho(u-\langle u\rangle_\rho)(t)\|_2^2.
\end{equation}
Combining \eqref{eq:conserved-limits} with \eqref{eq:physical-total-energy},
\eqref{eq:velocity-variance-decay} and the definition of \(\theta_*\), we arrive at
\begin{equation}\label{eq:conserved-limits-second}
 \langle\theta(t)\rangle_\rho
 =\theta_*-\frac{\|\sqrt\rho(u-\langle u\rangle_\rho)(t)\|_2^2}{2c_vm_0}
 \longrightarrow\theta_*
 \qquad (t\to\infty).
\end{equation}
Moreover, by \eqref{eq:lift-log-anchor}, \eqref{eq:decay-thermal},
\eqref{eq:velocity-variance-decay} and \eqref{eq:conserved-limits-second}, we deduce that
\begin{equation}\label{eq:thermal-variance-decay}
\begin{aligned}
 \int_{\R^2}\rho|\theta-\theta_*|^2\dd x
 =&\int_{\R^2}\rho|\theta-\langle\theta\rangle_\rho|^2\dd x
            +m_0|\langle\theta\rangle_\rho-\theta_*|^2\\
 \le& C\log(e+t)\|\nabla\theta(t)\|_2^2
                     +C\|\sqrt\rho(u-\langle u\rangle_\rho)(t)\|_2^4\\
 \le& C\log^2(e+t)(1+t)^{-9/8}\longrightarrow0
 \qquad (t\to\infty).
\end{aligned}
\end{equation}

Finally, for the solution \((\rho,u,\theta)\) of \eqref{eq:fns}--\eqref{eq:initial}, it follows from \eqref{eq:mass},
\eqref{eq:moving-frame} and \eqref{eq:physical-total-energy} that
\(\langle u(t)\rangle_\rho=u_*\) and
\begin{equation}\label{eq:unshifted-kinetic-identity}
 \frac{1}{2}\int_{\R^2}\rho|u|^2\dd x
 =\frac{m_0}{2}|u_*|^2+\frac{1}{2}\|\sqrt\rho(u-u_*)(t)\|_2^2.
\end{equation}
Since \(\|\sqrt\rho(u-u_*)(t)\|_2\) is invariant under \eqref{eq:moving-frame},
combining \eqref{eq:velocity-variance-decay}, \eqref{eq:conserved-limits-second},
\eqref{eq:thermal-variance-decay} and \eqref{eq:unshifted-kinetic-identity} yields
\eqref{eq:main-decay} and
\[
 \frac{1}{2}\int_{\R^2}\rho|u|^2\dd x\longrightarrow\frac{m_0}{2}|u_*|^2,\qquad
 \frac{1}{m_0}\int_{\R^2}\rho\theta\dd x\longrightarrow\theta_*
 \qquad (t\to\infty).
\]
This completes the proof of Theorem~\ref{thm:decay}.
\end{proof}

\appendix
\section{The far-field limit of the initial velocity}\label{app:initial-velocity}

\begin{proposition}\label{prop:initial-velocity-limit}
Under the initial regularity and compatibility assumptions of Theorem~\ref{thm:target}, we have
\[
 \lim_{j\to\infty}\sup_{k>j}|m_k-m_j|=0,
\]
where $A_j=B_{2^{j+1}}\setminus B_{2^j}$ and
$m_j=|A_j|^{-1}\int_{A_j}u_0(x)\dd x$ for $j=0,1,\ldots$.
The limit $u_\infty=\lim_{j\to\infty}m_j$ exists in $\R^2$, and we have
\[
 \|u_0-u_\infty\|_{L^{16}(\R^2)}
 \le C\|\nabla u_0\|_{L^{16/9}(\R^2)}<\infty,
\]
and
\begin{equation}\label{eq:initial-velocity-remainder}
 u_0-u_\infty\in W^{1,16}(\R^2),\qquad
 \lim_{|x|\to\infty}u_0(x)=u_\infty .
\end{equation}
\end{proposition}
\begin{proof}
Write $P_0=\rho_0\theta_0$. By \eqref{eq:mass} and \eqref{eq:initial-data-bound}, we have
\[
 \|\sqrt{\rho_0}g\|_1\le\sqrt{m_0}\,\|g\|_2,\qquad
 \|P_0\|_2\le\sqrt{\bar\rho}\,\|\sqrt{\rho_0}\theta_0\|_2.
\]
For $R>0$, we set $\chi_R(x)=\chi(x/R)$, where
\[
 \chi\in C_c^\infty(B_2),\qquad 0\le\chi\le1,\qquad
 \chi(x)=1\quad (|x|\le1).
\]
Integrating \eqref{eq:compatibility} against $\chi_R$ and using
$\|\nabla\chi_R\|_2=\|\nabla\chi\|_2$, we obtain
\[
 \left|\int_{\R^2}\chi_R\sqrt{\rho_0}g\dd x\right|
 \le C\bigl(\|\nabla u_0\|_{L^2(B_{2R}\setminus B_R)}
                 +\|P_0\|_{L^2(B_{2R}\setminus B_R)}\bigr)
 \longrightarrow0 \qquad (R\to\infty).
\]
By Lebesgue's theorem, we arrive at
\begin{equation}\label{eq:initial-force-zero-mean}
 \int_{\R^2}\sqrt{\rho_0}\,g\dd x=0.
\end{equation}

Fix $\alpha/2<\beta<\gamma<\alpha$, where $\alpha=\min\{a/2,1\}$.
Since $q>2$, \eqref{eq:data} gives $\bar x^a\rho_0\in L^1\cap L^\infty$.
Using \((1+|x|^2)^\gamma\le\bar x^a\) and Cauchy--Schwarz, we obtain
\begin{equation}\label{eq:weighted-initial-lame-data}
 \|(1+|x|^2)^{\gamma/2}\sqrt{\rho_0}\,g\|_1
 \le \|\bar x^a\rho_0\|_1^{1/2}\|g\|_2,
 \qquad
 \|(1+|x|^2)^{\gamma/2}\sqrt{\rho_0}\,g\|_2
 \le \|\bar x^a\rho_0\|_\infty^{1/2}\|g\|_2,
\end{equation}
and
\begin{equation}\label{eq:weighted-initial-lame-data-second}
 \|(1+|x|^2)^{\beta/2}P_0\|_2
 \le\|\bar x^a\rho_0\|_\infty^{1/2}\|\sqrt{\rho_0}\theta_0\|_2.
\end{equation}
By \eqref{eq:weighted-initial-lame-data}, we have $\sqrt{\rho_0}\,g\in L^1\cap L^2\subset L^{4/3}$.
Equations \eqref{eq:compatibility} and \eqref{eq:lift-lame-gradient-representation} give
\[
 \partial_i u_0
 =\frac{\partial_i\nabla\Delta^{-1}P_0}{2\mu+\lambda}
  -(\partial_i\Gamma)*(\sqrt{\rho_0}\,g),\qquad i=1,2.
\]
For the fundamental matrix $\Gamma=(\Gamma_{jk})_{j,k=1}^2$ in
\eqref{eq:lift-lame-fundamental-solution}, we have
\[
 \bigl(\mu\Delta+(\mu+\lambda)\nabla\Div\bigr)\Gamma
 =\delta_0 I_2\quad\text{in }\mathcal D'(\R^2),\qquad
 (\partial_i\Gamma)_{jk}=\partial_{x_i}\Gamma_{jk},\quad i=1,2,
\]
where $\delta_0$ is the Dirac function at the origin.
For $x\ne0$ and $|y|\le |x|/2$, the mean value theorem,
\eqref{eq:lift-lame-fundamental-solution} and $0<\gamma<1$ give
\[
 |\partial_i\Gamma(x-y)-\partial_i\Gamma(x)|
 \le C|y|\,|x|^{-2}\le C|y|^\gamma|x|^{-1-\gamma}.
\]
Using \eqref{eq:initial-force-zero-mean}, we obtain, for $|x|>1$,
\[
 |(\partial_i\Gamma)*(\sqrt{\rho_0}g)(x)|
 \le \frac{C}{|x|^{1+\gamma}}\big\||\,\cdot\,|^\gamma\sqrt{\rho_0}g\big\|_1
     +\frac{C}{|x|^\gamma}I_1\big(|\,\cdot\,|^\gamma\sqrt{\rho_0}|g|\big)(x).
\]
For $1<r<2/(2-\gamma+\beta)$, we have
\[
 (\gamma-\beta)\frac{r}{r-1}>2,\qquad
 \frac{r-1}{r}+\frac{2-r}{2r}=\frac12.
\]
By H\"older's inequality and \eqref{eq:riesz-potential-lp}, we obtain
\begin{equation}\label{eq:initial-force-exterior}
\begin{aligned}
 &\bigl\||x|^{\beta-\gamma}I_1(|\,\cdot\,|^\gamma\sqrt{\rho_0}|g|)
       \bigr\|_{L^2(\{|x|>1\})}\\
 \le{}&\bigl\||x|^{\beta-\gamma}\bigr\|_{L^{r/(r-1)}(\{|x|>1\})}
       \bigl\|I_1(|\,\cdot\,|^\gamma\sqrt{\rho_0}|g|)\bigr\|_{2r/(2-r)}\\
 \le{}&C\bigl\||x|^\gamma\sqrt{\rho_0}\,g\bigr\|_r.
\end{aligned}
\end{equation}
On $B_1$, \eqref{eq:riesz-potential-lp} and \eqref{eq:lift-lame-fundamental-solution} give
\begin{equation}\label{eq:initial-force-ball}
\begin{aligned}
 \bigl\|(\partial_i\Gamma)*(\sqrt{\rho_0}\,g)\bigr\|_{L^2(B_1)}\le&C\bigl\|I_1(\sqrt{\rho_0}|g|)\bigr\|_{L^2(B_1)}\\
 \le&C|B_1|^{(r-1)/r}\bigl\|I_1(\sqrt{\rho_0}|g|)\bigr\|_{2r/(2-r)}\\
 \le& C\|\sqrt{\rho_0}\,g\|_r.
\end{aligned}
\end{equation}
Interpolating the two estimates in \eqref{eq:weighted-initial-lame-data}, we obtain
\[
 \bigl\|(1+|x|^2)^{\gamma/2}\sqrt{\rho_0}\,g\bigr\|_r
 \le \|\bar x^a\rho_0\|_1^{1/r-1/2}
     \|\bar x^a\rho_0\|_\infty^{1-1/r}\|g\|_2.
\]
Applying \eqref{eq:weighted-riesz-pressure} and \eqref{eq:weighted-initial-lame-data-second}
to $P_0$ and \eqref{eq:initial-force-exterior}--\eqref{eq:initial-force-ball}
to $(\partial_i\Gamma)*(\sqrt{\rho_0}\,g)$, we obtain
\begin{equation}\label{eq:initial-farfield-gradient}
\begin{aligned}
 \|\bar x^{\alpha/2}\nabla u_0\|_2
 \le& C\|(1+|x|^2)^{\beta/2}\nabla u_0\|_2\\
 \le& C\bigl(\|\bar x^a\rho_0\|_1^{1/2}
             +\|\bar x^a\rho_0\|_\infty^{1/2}\bigr)\|g\|_2
        +C\|\bar x^a\rho_0\|_\infty^{1/2}\|\sqrt{\rho_0}\theta_0\|_2 .
\end{aligned}
\end{equation}
Since $8\alpha>2$, we have $\bar x^{-\alpha/2}\in L^{16}(\R^2)$.
H\"older's inequality and \eqref{eq:initial-farfield-gradient} yield
\begin{equation}\label{eq:initial-gradient-16-9}
 \|\nabla u_0\|_{16/9}
 \le\|\bar x^{-\alpha/2}\|_{16}\,\|\bar x^{\alpha/2}\nabla u_0\|_2<\infty .
\end{equation}

Since $u_0\in L^1_{\rm loc}(\R^2)$, each $m_j$ is finite and
\eqref{eq:initial-gradient-16-9} gives $u_0\in W^{1,1}(B)$ on every ball $B$.
Applying \cite[Lemma~7.16]{GilbargTrudinger2001} on \(B\) and
\eqref{eq:riesz-potential-lp} with \(p=16/9\), we obtain
\[
\begin{aligned}
 \|u_0-(u_0)_B\|_{L^{16}(B)}
 \le& \frac{2}{\pi}\bigl\|I_1(\mathbf{1}_B|\nabla u_0|)\bigr\|_{16}
 \le C\|\nabla u_0\|_{L^{16/9}(B)},
\end{aligned}
\]
where \((u_0)_B=|B|^{-1}\int_B u_0\dd x\) and \(C\) is independent of \(B\).
For \(B=B_{2^{j+1}}\), H\"older's inequality gives
\[
 |m_j-(u_0)_B|\le |A_j|^{-1/16}\|u_0-(u_0)_B\|_{L^{16}(B)}.
\]
Since \(|B|/|A_j|=4/3\), it follows that
\begin{equation}\label{eq:initial-annular-estimate}
 \|u_0-m_j\|_{L^{16}(B_{2^{j+1}})}
 \le \bigl(1+(4/3)^{1/16}\bigr)\|u_0-(u_0)_B\|_{L^{16}(B)}
 \le C\|\nabla u_0\|_{16/9},
\end{equation}
where \(C\) is independent of \(j\).
For any integers \(k>j\), we have \(A_j\subset B_{2^{k+1}}\) and \(|A_j|=3\pi\,2^{2j}\).
Using \eqref{eq:initial-gradient-16-9}
and \eqref{eq:initial-annular-estimate} at index \(k\), we obtain
\begin{equation}\label{eq:initial-annular-cauchy}
\begin{aligned}
 \sup_{k>j}|m_j-m_k|
 =&\sup_{k>j}\left|\frac1{|A_j|}\int_{A_j}(u_0-m_k)\dd x\right|\\
 \le& |A_j|^{-1/16}\sup_{k>j}\|u_0-m_k\|_{L^{16}(B_{2^{k+1}})}\\
 \le& \frac{C\,2^{-j/8}}{(3\pi)^{1/16}}\|\nabla u_0\|_{16/9}
 \longrightarrow0\qquad(j\to\infty).
\end{aligned}
\end{equation}
Thus $\{m_j\}$ is Cauchy in $\R^2$. By completeness its limit exists and we denote it by $u_\infty$.
Fatou's lemma and \eqref{eq:initial-annular-estimate} give
\[
 \|u_0-u_\infty\|_{16}
 \le\liminf_{j\to\infty}\|u_0-m_j\|_{L^{16}(B_{2^{j+1}})}
 \le C\|\nabla u_0\|_{16/9}.
\]

Finally, $\nabla u_0\in H^1(\R^2)\subset L^{16}(\R^2)$ by \eqref{eq:data-second}.
Thus $u_0-u_\infty\in W^{1,16}(\R^2)$ has a continuous representative.
Applying the embedding $W^{1,16}(B_1(x))\hookrightarrow C(\overline{B_1(x)})$
for $|x|\ge R+1$, we obtain
\[
 \sup_{|x|\ge R+1}|u_0(x)-u_\infty|
 \le C\bigl(\|u_0-u_\infty\|_{L^{16}(\R^2\setminus B_R)}
             +\|\nabla u_0\|_{L^{16}(\R^2\setminus B_R)}\bigr)
 \longrightarrow0\qquad(R\to\infty). \qedhere
\]
\end{proof}

\section*{Statements and Declarations}
{\bf Conflict of interest statement}
On behalf of all authors, the corresponding author states that there is no conflict of interest. \\
{\bf  Data availability statement} 
No data were used for this research.

\section*{Acknowledgements} 
Xinyu Fan is supported by  the Fundamental Research Funds for the Central Universities, Nankai University.  

\section*{AI assistance}
The authors acknowledge the use of AI to explore proof strategies, check calculations, and polish
the manuscript.  The authors take full responsibility for the mathematical arguments and
the final manuscript.

\end{document}